\documentclass[a4paper,11pt]{article}
\usepackage{graphicx}
\usepackage{a4wide}
\usepackage[top=80pt,bottom=70pt]{geometry}
\usepackage{amssymb,amsmath,amsthm}
\usepackage{siunitx}
\usepackage{subfigure,color,colordvi,graphicx,xcolor}
\usepackage[only,llbracket,rrbracket,interleave]{stmaryrd}
\usepackage[superscript,sort,compress]{cite}

\makeatletter
\newcommand*\wt[2][0.2ex]{%
        \begingroup
        \mathchoice{\wt@helper{#1}{#2}{\displaystyle}{\textfont}}
                   {\wt@helper{#1}{#2}{\textstyle}{\textfont}}
                   {\wt@helper{#1}{#2}{\scriptstyle}{\scriptfont}}
                   {\wt@helper{#1}{#2}{\scriptscriptstyle}{\scriptscriptfont}}%
        \endgroup
        #2%
}
\newcommand*\wt@helper[4]{%
        \def\currentfont{\the#41}%
        \def\currentskewchar{\char\the\skewchar\currentfont}%
        \setbox\tw@\hbox{\currentfont#2\currentskewchar}%
        \dimen@ii\wd\tw@
        \setbox\tw@\hbox{\currentfont#2{}\currentskewchar}%
        \advance\dimen@ii-\wd\tw@
        \rlap{\raisebox{-#1}{$\m@th#3\kern\dimen@ii\widetilde{\phantom{#2}}$}}%
}
\makeatother

\newcommand{\bm}[1]{\text{\boldmath $#1$\unboldmath}}
\newcommand{\grad}{\bm{\nabla}}
\newcommand{\tr}{\operatorname{tr}}
\newcommand{\bn}{\bm{n}}
\newcommand{\bq}{\bm{q}}
\newcommand{\bu}{\bm{u}}
\newcommand{\bv}{\bm{v}}
\newcommand{\bw}{\bm{w}}
\newcommand{\bhu}{\widehat{\bu}}
\newcommand{\bhv}{\widehat{\bv}}
\newcommand{\bL}{\bm{L}}
\newcommand{\hu}{\hat{u}}
\newcommand{\hv}{\hat{v}}
\newcommand{\Insd}{\mat{I}_{\nsd}}
\newcommand{\stress}{\bm{\sigma}}
\newcommand{\defo}[1]{\bm{\varepsilon}#1}
\newcommand{\gradS}{\bm{\nabla}_{\!\texttt{S}}}
\newcommand{\gradW}{\bm{\nabla}_{\!\texttt{W}}}
\newcommand{\bD}{\mat{D}}
\newcommand{\bDHalf}{\bD^{1/2}}

\newcommand{\bN}{\mat{N}}
\newcommand{\bT}{\mat{T}}

\newcommand{\stressV}{\bm{\sigma}_{\texttt{V}}}
\newcommand{\strainV}{\bm{\varepsilon}_{\texttt{V}}}
\newcommand{\numel}{\ensuremath{\texttt{n}_{\texttt{el}}}}
\newcommand{\nsd}{\ensuremath{\texttt{n}_{\texttt{sd}}}}
\newcommand{\msd}{\ensuremath{\texttt{m}_{\texttt{sd}}}}
\newcommand{\nrr}  {\ensuremath{\texttt{n}_{\texttt{rr}}}}
\newcommand{\RR}{\mathbb{R}}

\newcommand{\eltwo}{\ensuremath{\mathcal{L}_2}}

\newcommand{\Vh}{\ensuremath{\mathcal{V}^h}}
\newcommand{\Wh}{\ensuremath{\mathcal{W}^h}}
\newcommand{\HWh}{\ensuremath{\reallywidehat{\mathcal{W}}^h}}
\newcommand{\Czero}{\ensuremath{\mathcal{C}^0}}
\newcommand{\Poly}[1][1]{\ensuremath{\mathcal{P}^{#1}}}
\newcommand{\jump}[1]{\llbracket #1\rrbracket}

\newcommand{\vect}[1]{\mathbf{#1}}
\newcommand{\mat}[1]{\mathbf{#1}}

\newcommand{\HDG}[1]{{#1}_{\texttt{HDG}}}
\newcommand{\CG}[1]{{#1}_{\texttt{CG}}}
\newcommand{\KInt}{\mat{K}_{\texttt{I}}}
\newcommand{\Kdeg}[1]{{\text{#1}}_{k}}
\newcommand{\KPdeg}[1]{{\text{#1}}_{k{+}1}}

\usepackage{scalerel}
\usepackage{stackengine}
\stackMath
\def\hatgap{0pt}
\def\subdown{-2pt}
\newcommand\reallywidehat[2][]{
\renewcommand\stackalignment{l}
\stackon[\hatgap]{#2}{
\stretchto{\scalerel*[\widthof{$#2$}]{\kern-.6pt\bigwedge\kern-.6pt}{\rule[-\textheight/2]{1ex}{\textheight}}}{0.5ex}_{\smash{\belowbaseline[\subdown]{\scriptstyle#1}}}}}

\newenvironment{keywords}{\begin{quote}\emph{\textbf{Keywords:}}}{\end{quote}}

\theoremstyle{definition}

\RequirePackage{algorithm}[0.1] 
\usepackage{algorithmic} 

\graphicspath{{figsPDF/}{figsPNG/}{fig_arXiv/}}

\begin{document}
\title{Hybrid coupling of CG and HDG discretizations based on Nitsche's method}

\author{
\renewcommand{\thefootnote}{\arabic{footnote}}
			  Andrea La Spina\footnotemark[1]\textsuperscript{ \ ,}\footnotemark[2] \ ,
			  Matteo Giacomini\footnotemark[2], \ and
             Antonio Huerta\footnotemark[2]
}

\date{\today}
\maketitle

\renewcommand{\thefootnote}{\arabic{footnote}}

\footnotetext[1]{Lehrstuhl f\"ur Numerische Mechanik (LNM), Technische Universit\"at M\"unchen, Garching b. M\"{u}nchen, Germany}
\footnotetext[2]{Laboratori de C\`alcul Num\`eric (LaC\`aN), ETS de Ingenieros de Caminos, Canales y Puertos, Universitat Polit\`ecnica de Catalunya, Barcelona, Spain
\vspace{5pt}\\
Corresponding author: Matteo Giacomini. \textit{E-mail:} \texttt{matteo.giacomini@upc.edu}
}

\begin{abstract}
A strategy to couple continuous Galerkin (CG) and hybridizable discontinuous Galerkin (HDG) discretizations based only on the HDG hybrid variable is presented for linear thermal and elastic problems.
The hybrid CG-HDG coupling exploits the definition of the numerical flux and the trace of the solution on the mesh faces to impose the transmission conditions between the CG and HDG subdomains.
The continuity of the solution is imposed in the CG problem via Nitsche's method, whereas the equilibrium of the flux at the interface is naturally enforced as a Neumann condition in the HDG global problem.
The proposed strategy does not affect the core structure of CG and HDG discretizations. In fact, the resulting formulation leads to a minimally-intrusive coupling, suitable to be integrated in existing CG and HDG libraries.
Numerical experiments in two and three dimensions show optimal global convergence of the stress and superconvergence of the displacement field, locking-free approximation, as well as the potential to treat structural problems of engineering interest featuring multiple materials with compressible and nearly incompressible behaviors.
%
%
\end{abstract}


\section{Introduction}
\label{sc:Introduction}

The finite element (FE) method has been extensively used in simulation-based engineering over the last $50$ years. Strategies to enhance FE computations by coupling them with other discretization methods have received a lot of attention both in the mathematical and engineering community. 
Examples of successful couplings include, e.g., FE and meshless methods~\cite{SFM-HF:00,SFM-FH:04,SFM-HFL:04,SFM-FBH:05}, FE and finite volume methods~\cite{Casadei-CL-11,Riviere-CMR-11,Vassilevski-COV-18}, as well as advanced discretization techniques as mixed and discontinuous Galerkin (DG) methods~\cite{Firoozabadi-MF-16,Firoozabadi-HF-18}.

The idea of coupling distinct numerical techniques in different regions of the computational domain may also be interpreted in the framework of domain decomposition (DD) methods.
In this context, and for nonverlapping strategies, a key aspect is represented by the so-called \emph{transmission conditions} to be imposed between two different regions of the domain. Such conditions enforce continuity of the solution and of the flux/stress across the interface and have been treated in the literature either via Lagrange multipliers~\cite{Bertsekas-82book} or using Nitsche's method~\cite{Nitsche1971}.
On the one hand, couplings involving FE discretizations and Lagrange multipliers has lead to the well-known Mortar element method~\cite{Bernardi-BMP-92,Bernardi-BMP-93,LeTallec-LTS-95,Wohlmuth-WW-98,Maday-AMW-99,Belgacem-99,Thomas-ALT-99,Wheeler-ACWY-00,Buffa-BMR-01}.
On the other hand, DD methods based on DG approximations have been studied using both Lagrange multipliers~\cite{Riviere-RW-02,Girault-GSWY-08,Wheeler-MW-14} and Nitsche's method~\cite{Stenberg-98,Becker-BHS-03}.

Among existing numerical methods, the coupling of CG and DG in different regions of the domain has been extensively studied in the literature to exploit the advantages of both approaches.
This is of special interest in the context of multiphysics and multimaterial problems in which different regions of the computational domain feature distinct physical properties, for which specific discretizations need to be devised. 
On the one hand, CG leads to computationally efficient discretizations with a limited number of degrees of freedom~\cite{braess2001finite}. 
On the other hand, DG provides a flexbile paradigm to handle meshes with hanging nodes and construct nonuniform polynomial degree and high-order approximations~\cite{Karniadakis-CKS-00,Riviere2008,ErnBook,Cangiani2017}. Moreover, DG methods have proved to be very efficient in stabilizing convection terms in conservation laws~\cite{Bassi-BR-97,AbgrallRicchiutoECM}.

The idea of coupling CG and DG, specifically the local discontinuous Galerkin (LDG)~\cite{Chu-CS-98}, via an appropriate definition of the numerical flux was first proposed in~\cite{Perugia-PS-01} to handle nonmatching grids.
The advantage of a CG-LDG coupling in presence of convection-dominated problems was later discussed in~\cite{Dawson-DP-02,Dawson-DP-03,Dawson-DP-04}.

More recently, hybrid discretization techniques, especially HDG~\cite{cockburn2008superconvergent,Jay-CGL:09,Nguyen-NPC:09,Nguyen-NPC:09b,Nguyen-CNP:10,Nguyen-NPC:11} and the hybrid high-order (HHO) method~\cite{DiPietro-DPEL-14,Ern-DPE-15,Ern-AEP-18,Ern-AEP-19}, have gained a lot of attention owing to their reduced computational costs with respect to classical DG approaches like LDG.
In the context of HDG, coupling with the boundary element method has been discussed in~\cite{Sayas-CGS-12,Sayas-FHS-17}, whereas a first attempt to couple HDG and CG discretizations has been proposed in~\cite{SFM-PTFM-19}.
This approach requires the introduction of an appropriate projection operator to enforce the transmission condition in the HDG local problem. This results in a coupling of local and global degrees of freedom of the HDG problem with the ones of the CG discretization, making the implementation of this strategy in existing CG and HDG libraries intrusive.

This work proposes a CG-HDG coupling in which solely the HDG hybrid variable is exploited in the transmission conditions.
This approach requires two ingredients. On the one hand, an appropriate definition of the trace of the numerical normal flux at the interface between the CG and HDG subdomains. On the other hand, the weak imposition of Dirichlet boundary conditions in the CG formulation via Nitsche's method. A formulation of the latter in an optimization framework is described in~\cite{MG-G:18}.
The resulting hybrid CG-HDG coupling does not affect the structure of the core CG and HDG matrices, thus leading to a minimally-intrusive implementation of this technique in existing CG and HDG libraries.

This approach is of special interest in the context of linear elastic problems since displacement-based formulations fail to provide locking-free approximations in nearly incompressible materials, when low-order CG discretizations are utilized~\cite{Brenner-BS-92}.
To remedy this issue, mixed~\cite{Veubeke1975,Brezzi-ABD-84,Stenberg1988,Arnold-AW-02} and equilibrium formulations~\cite{Moitinho-MM-17book}, as well as discretizations based on the nonconforming Crouzeix-Raviart element~\cite{Crouzeix-CR-73} and DG techniques~\cite{Hansbo-HL-02,Cockburn-CSW-06,Demkowicz-BDGQ-12} have been proposed.
The above mentioned HDG method also provides locking-free approximations~\cite{soon2009hybridizable,Fu-FCS-15,Cockburn-CS-13,Cockburn-CF-17,Qiu-QSS-18} while preserving the advantages of a DG discretization with hybridization. It relies on a mixed hybrid formulation~\cite{brezzi1991mixed} with polynomial approximations discontinuous element-by-element.
More precisely, the HDG formulation for linear elasticity proposed in~\cite{RS-SGKH:18,RS-SGH:19} is utilized. Exploiting the well-known Voigt notation for second-order symmetric tensors~\cite{FishBelytschko2007}, this approach provides optimal convergence of the stress tensor and superconvergence of the displacement field using equal order approximation for all the variables, even for low-order polynomial functions. 
Hence, a CG approximation in the compressible region of the domain is coupled with an HDG solver for the nearly incompressible one in order to study multimaterial problems of engineering interest.

The remaining of this paper is organized as follows.
First, a linear thermal problem is considered to introduce the CG and HDG discretizations, as well as the proposed hybrid coupling (Section~\ref{sc:CG-HDG-Poisson}).
In Section~\ref{sc:CG-HDG-Elasticity}, the coupled CG-HDG discretization is presented for a linear elastic problem exploiting the HDG-Voigt formulation.
Section~\ref{sc:NumericalStudies} is devoted to the numerical validation of the hybrid coupling in two dimensions. More precisely, optimal orders of convergence, a sensitivity analysis to the parameters of the method and robustness in the incompressible limit are verified. Special emphasis is put in the development of a strategy based on nonuniform polynomial degree approximation to obtain optimal global convergence of the stress of order $k {+} 1$ and superconvergence of the displacement field of order $k {+} 2$.
In Section~\ref{sc:multiMaterial}, two and three-dimensional elastic problems involving composite materials with space-dependent mechanical properties are analyzed and Section~\ref{sec:conclusion} summarizes the results of this work.

\section{CG-HDG coupling for thermal problems} 
\label{sc:CG-HDG-Poisson}

The hybrid coupling of CG and HDG discretizations is presented for a thermal problem described by a second-order scalar elliptic partial differential equation (PDE).
Consider an open bounded domain $\Omega \in \RR^{\nsd}$ with boundary $\partial\Omega {=} \Gamma_D \cup \Gamma_N$ such that $\Gamma_D \cap \Gamma_N {=} \emptyset$ and $\nsd$ being the number of spatial dimensions.
The strong form of the Poisson equation is
\begin{equation}\label{eq:PoissonStrongForm}
\left\{
\begin{aligned}
-\grad {\cdot} \grad u &= f           &&\text{in $\Omega$,}   \\
u &= u_D                                     &&\text{on $\Gamma_D$,} \\
\bn {\cdot} \grad u &=t               &&\text{on $\Gamma_N$,}
\end{aligned}
\right.
\end{equation}
where $u$ represents the unknown temperature, $f$ is a user-prescribed source term and $u_D$ and $t$ denote the Dirichlet and Neumann boundary data, respectively.

In the following subsections, the discrete forms of the CG and HDG approximations are first recalled separately and the coupling strategy based on Nitsche's method is thus presented.

\subsection{CG approximation}
\label{sc:PoissonCG}

Assume that $\Omega$ is partitioned in $\numel$ disjoint subdomains $\Omega_e$ such that
\begin{equation}\label{eq:brokenDomain}
\Omega = \bigcup_{e=1}^{\numel} \Omega_e,
\quad \Omega_i \cap \Omega_j = \emptyset \text{ for } i \neq j .
\end{equation}
The following discrete functional spaces are introduced
\begin{subequations}\label{eq:CGspaces}
\begin{align}
\Vh(\Omega) & {:=}
\lbrace 
v \in \Czero(\overline{\Omega}) : \, v\vert_{\Omega_e} \in \Poly[k](\Omega_e) \, \forall \Omega_e, e {=} 1,\ldots,\numel , \, v\vert_{\Gamma_D} {=} u_D
\rbrace,
\\
\Vh_0(\Omega) & {:=}
\lbrace
v \in \Czero(\overline{\Omega}) : \, v\vert_{\Omega_e} \in \Poly[k](\Omega_e) \, \forall \Omega_e, e {=} 1,\ldots,\numel , \, v\vert_{\Gamma_D} {=} 0
\rbrace ,
\end{align}
\end{subequations}
where $\Poly[k](\Omega_e)$ is the space of polynomial functions of complete degree at most $k {\geq} 1$ on each mesh element $\Omega_e, \, e {=} 1,\ldots,\numel$.
Moreover, $(\cdot,\cdot)_D$ and $\langle \cdot, \cdot \rangle_S$ are introduced to denote the usual $\eltwo$ inner products on a generic subdomain $D$ and on a generic surface $S$, namely
$$
(u, v)_D {:=} \int_D{u v \, d\Omega} , \quad \langle u, v \rangle_S {:=} \int_S{u v \, d\Gamma} .
$$

The discrete form of the CG approximation of Equation~\eqref{eq:PoissonStrongForm} is thus obtained by integrating by parts the corresponding weak form and it reads: find $u^h \in \Vh(\Omega)$ such that
\begin{equation}\label{eq:PoissonCG}
( \grad v,\grad u^h )_{\Omega} = (v, f )_{\Omega} + \langle v, t \rangle_{\Gamma_N}
\end{equation}
for all $v \in \Vh_0(\Omega)$.

\subsection{HDG approximation}
\label{sc:PoissonHDG}

Consider the partition of the domain introduced in Equation~\eqref{eq:brokenDomain} and define the internal interface as
$$
\Gamma := \left[ \bigcup_{e=1}^{\numel} \partial\Omega_e \right] \setminus \partial\Omega .
$$

The HDG formulation of the Poisson problem under analysis is obtained by rewriting Equation~\eqref{eq:PoissonStrongForm} element-by-element on $\Omega_e, \, e {=} 1,\ldots,\numel$ as a system of first-order PDEs, via the introduction of the mixed variable $\bq {=} {-} \grad u$ and the hybrid variable $\hu$ representing the trace of the solution on $\partial\Omega_e \setminus \Gamma_D$, namely
\begin{equation}\label{eq:PoissonHDGstrong}
\left\{
\begin{aligned}
\bq + \grad u &= \bm{0}                            && \text{in $\Omega_e, \, e {=} 1,\ldots,\numel$,} \\
\grad {\cdot} \bq &= f                               && \text{in $\Omega_e, \, e {=} 1,\ldots,\numel$,} \\
u &= u_D                                                   && \text{on $\partial\Omega_e \cap \Gamma_D$,} \\
u &= \hu                                                   && \text{on $\partial\Omega_e \setminus \Gamma_D$,} \\
\bn {\cdot} \bq &= -t                                && \text{on $\partial\Omega_e \cap \Gamma_N$,} \\
\jump{ u \bn } &= \bm{0}                         && \text{on $\Gamma$,}   \\
\jump{ \bn {\cdot} \bq } &= 0                  &&\text{on $\Gamma$,}
\end{aligned}
\right.
\end{equation}
where the last two equations are \emph{transmission conditions} enforcing the continuity of the solution and of the normal flux across the interface $\Gamma$. The \emph{jump} operator $\jump{\cdot}$ is defined according to~\cite{AdM-MFH:08} as the sum of the values from the elements $\Omega_l$ and $\Omega_r$, respectively on the left and on the right of a given interface, namely
$$
\jump{ \odot } = \odot_l + \odot_r .
$$

The solution of the HDG problem is performed in two stages~\cite{cockburn2008superconvergent,Jay-CGL:09,Nguyen-NPC:09,Nguyen-NPC:09b,Nguyen-CNP:10,Nguyen-NPC:11}.
First, a set of $\numel$ local problems is defined to determine $(u_e,\bq_e)$ as functions of the unknown hybrid variable $\hu$ in each element $\Omega_e, \, e {=} 1,\ldots,\numel$, that is
\begin{equation}\label{eq:PoissonHDGstrongLocal}
\left\{
\begin{aligned}
\bq_e + \grad u_e &= \bm{0}                            && \text{in $\Omega_e, \, e {=} 1,\ldots,\numel$,} \\
\grad {\cdot} \bq_e &= f                               && \text{in $\Omega_e, \, e {=} 1,\ldots,\numel$,} \\
u_e &= u_D                                                   && \text{on $\partial\Omega_e \cap \Gamma_D$,} \\
u_e &= \hu                                                   && \text{on $\partial\Omega_e \setminus \Gamma_D$.}
\end{aligned}
\right.
\end{equation}
Second, the trace of the solution on $\Gamma \cup \Gamma_N$ is computed by solving a global problem given by the above introduced transmission conditions
\begin{equation}\label{eq:PoissonHDGstrongGlobal}
\left\{
\begin{aligned}
\jump{ u \bn } &= \bm{0}                         && \text{on $\Gamma$,}   \\
\jump{ \bn {\cdot} \bq } &= 0                  &&\text{on $\Gamma$,} \\
\bn {\cdot} \bq &= -t                                && \text{on $\Gamma_N$,} 
\end{aligned}
\right.
\end{equation}
where the first condition is automatically fulfilled owing to the Dirichlet boundary condition $u_e {=} \hu$ on $\partial\Omega_e \setminus \Gamma_D$ imposed in the local problem and to the uniqueness of the hybrid variable on each mesh edge (respectively, face in 3D).

Following the rationale in Section~\ref{sc:PoissonCG}, the discrete functional spaces 
\begin{subequations}\label{eq:HDGspaces}
\begin{align}
\Wh(\Omega) & {:=}
\lbrace 
v \in \eltwo(\Omega) : \, v\vert_{\Omega_e} \in \Poly[k](\Omega_e) \forall \Omega_e, \, e {=} 1,\ldots,\numel
\rbrace,
\\
\HWh(S) & {:=}
\lbrace
\hv \in \eltwo(S) : \, \hv\vert_{\Gamma_i} \in \Poly[k](\Gamma_i) \forall \Gamma_i \subset S \subseteq \Gamma \cup \partial\Omega
\rbrace ,
\end{align}
\end{subequations}
are introduced for the HDG approximation and $\Poly[k](\Omega_e)$ and $\Poly[k](\Gamma_i)$ here denote the spaces of polynomial functions of complete degree at most $k {\geq} 1$ in $\Omega_e$ and on $\Gamma_i$, respectively.

Introduce the definition of the trace of the numerical normal flux 
\begin{equation}\label{eq:PoissonFlux}
\bn {\cdot} \reallywidehat{\bq}^h := 
\begin{cases}
\bn {\cdot} \bq_e^h + \tau (u_e^h - u_D) & \text{on $\partial\Omega_e\cap\Gamma_D$,} \\
\bn {\cdot} \bq_e^h + \tau (u_e^h - \hu^h) & \text{elsewhere,}  
\end{cases}
\end{equation}
where $\tau$ is a stabilization parameter critical to ensure stability and convergence of the HDG method, as extensively studied in a series of publications~\cite{cockburn2008superconvergent,Jay-CGL:09,Nguyen-NPC:09,Nguyen-NPC:09b,Nguyen-CNP:10,Nguyen-NPC:11} by Cockburn and coworkers.

The discrete HDG local problems are thus obtained integrating by parts the weak form of Equation~\eqref{eq:PoissonHDGstrongLocal} and exploiting the definition in Equation~\eqref{eq:PoissonFlux}. Following~\cite{RS-SH:16} and integrating the first equation by parts once and the second one twice, the resulting problem is: for $e {=} 1,\ldots,\numel$, given $u_D$ on $\Gamma_D$ and $\hu^h$ on $\Gamma\cup\Gamma_N$, find $(u_e^h,\bq_e^h) \in \Wh(\Omega_e) {\times} \left[\Wh(\Omega_e)\right]^{\nsd}$ such that
\begin{subequations}\label{eq:PoissonHDGdiscreteLocal}
\begin{align}
- (\bw, \bq_e^h)_{\Omega_e} {+} (\grad {\cdot} \bw, u_e^h)_{\Omega_e} 
&= \langle \bn {\cdot} \bw, u_D \rangle_{\partial\Omega_e\cap\Gamma_D} {+} \langle \bn {\cdot} \bw, \hu^h \rangle_{\partial\Omega_e\setminus\Gamma_D} ,
\\
(v, \grad {\cdot} \bq_e^h)_{\Omega_e} {+} \langle v, \tau u_e^h \rangle_{\partial\Omega_e} 
&= (v, f)_{\Omega_e} {+} \langle v, \tau u_D \rangle_{\partial\Omega_e\cap\Gamma_D} {+} \langle v, \tau \hu^h \rangle_{\partial\Omega_e\setminus\Gamma_D} ,
\end{align}
\end{subequations}
for all $(v,\bw) \in \Wh(\Omega_e) {\times} \left[\Wh(\Omega_e)\right]^{\nsd}$.

Similarly, the discrete form of the HDG global problem~\eqref{eq:PoissonHDGstrongGlobal} is: find $\hu^h \in \HWh(\Gamma\cup\Gamma_N)$ such that
\begin{equation}\label{eq:PoissonHDGdiscreteGlobal}
\sum_{e=1}^{\numel} \Bigl\{ 
\langle \hv, \bn {\cdot} \bq_e^h \rangle_{\partial\Omega_e\setminus\Gamma_D} {+} \langle \hv, \tau u_e^h \rangle_{\partial\Omega_e\setminus\Gamma_D} 
{-} \langle \hv, \tau \hu^h \rangle_{\partial\Omega_e\setminus\Gamma_D} \Bigr\} = {-} \sum_{e=1}^{\numel} \langle \hv, t \rangle_{\partial\Omega_e\cap\Gamma_N} ,
\end{equation}
for all $\hv \in \HWh(\Gamma\cup\Gamma_N)$.

For the sake of readability and except in case of ambiguity, the superscript $^h$ associated with the numerical counterpart of the unknowns in the continuous spaces and the subscript $_e$ related to the HDG elemental approximations will be henceforth omitted.

\subsection{Hybrid coupling based on Nitsche's method}
\label{sc:PoissonCoupling}

Consider a splitting of the domain $\Omega$ under analysis in two nonoverlapping subdomains $\CG{\Omega}$ and $\HDG{\Omega}$ such that $\Omega {=} \CG{\Omega} \cup \HDG{\Omega}$. The interface is defined as $\Gamma_I {:=} \CG{\overline{\Omega}} \cap \HDG{\overline{\Omega}}$.
In the following, the subscripts $\CG{}$ and $\HDG{}$ are employed to identify the quantities associated with the CG and HDG discretizations, respectively.

\begin{figure}
\centering
\includegraphics[width=0.6\columnwidth]{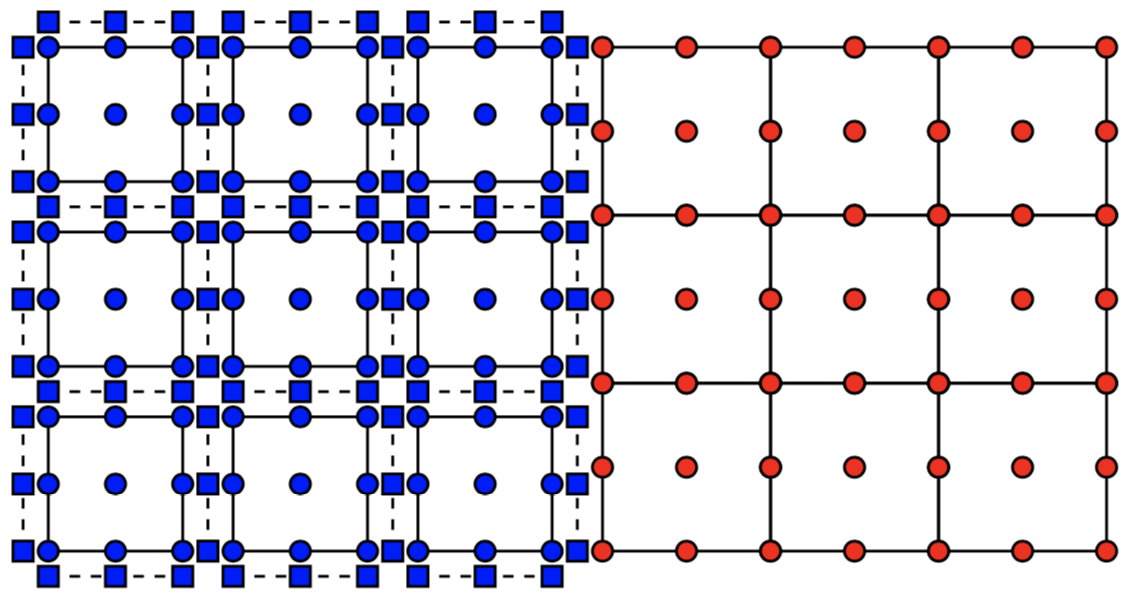}
\caption{Degrees of freedom of the coupled CG-HDG discretization using polynomial approximation of degree $k {=} 2$ in the CG subdomain $\CG{\Omega}$ (red) and in the HDG subdomain $\HDG{\Omega}$ (blue).}
\label{fig:discretization}
\end{figure}

The degrees of freedom of the  coupled problem are displayed in Figure~\ref{fig:discretization}: for the CG subdomain, the unknown $u$ is depicted in red, whereas the degrees of freedom of HDG local and global problems are represented by blue circles and blue squares, respectively.

The coupled CG-HDG strategy is obtained starting from the discretizations described in Section~\ref{sc:PoissonCG} and \ref{sc:PoissonHDG}, with two additional conditions enforcing continuity of the solution and of the normal flux across the interface $\Gamma_I$, namely
%
%
\begin{equation}\label{eq:couplingConditions}
\left\lbrace
\begin{aligned}
\CG{u} &= \HDG{\hu}                                                   &&\text{on $\Gamma_I$,} \\
-\bn_I {\cdot} \grad \CG{u} &= \bn_I {\cdot} \HDG{\bq}  &&\text{on $\Gamma_I$,}
\end{aligned}
\right.
\end{equation}
where $\CG{u}$ is solution of the CG problem in Equation~\eqref{eq:PoissonCG}, $\HDG{\hu}$ is obtained from the HDG global problem in Equation~\eqref{eq:PoissonHDGdiscreteGlobal}, the numerical flux $\bn_I {\cdot} \HDG{\bq}$ stems from the definition in Equation~\eqref{eq:PoissonFlux} and $\bn_I$ is the outer normal to the domain $\HDG{\Omega}$. 
It is worth noticing that the condition imposing the equilibrium of the flux in Equation~\eqref{eq:couplingConditions} already accounts for the uniqueness of the unit normal vector on the interface $\Gamma_I$. Henceforth, and unless in case of ambiguity, on the interface $\Gamma_I$ the normal $\bn$ is assumed to be the outer direction to the subdomain under analysis, that is $\bn {=} \bn_I$ for $\HDG{\Omega}$ and $\bn {=} {-} \bn_I$ for $\CG{\Omega}$.

From a practical point of view, on the one hand, the HDG global unknown $\HDG{\hu}$ is now defined also on the interface $\Gamma_I$ and the continuity of the solution is imposed as a Dirichlet boundary condition in the CG problem. On the other hand, the equilibrium of the flux is enforced via a Neumann boundary condition in the HDG global problem.
Hence, the proposed technique does not affect neither the HDG local problem nor the core routines of the CG solver, leading to a minimally-intrusive coupling strategy.

Recall that $\HDG{\hu}$ is solution of Equation~\eqref{eq:PoissonHDGdiscreteGlobal} on $\HDG{\Omega}$ and is thus known solely at the nodes of the HDG discretization. In this context, the Dirichlet condition $\CG{u} {=} \HDG{\hu}$ is enforced in the CG problem via the well-known Nitsche's method for the weak imposition of essential boundary conditions~\cite{Nitsche1971}.
Moreover, this choice provides a flexible framework for the coupled discretization, allowing for nonconforming meshes at the interface $\Gamma_I$, as well as nonuniform polynomial approximations in $\CG{\Omega}$ and $\HDG{\Omega}$.

Following the rationale utilized for HDG in Equation~\eqref{eq:PoissonFlux}, the trace of the CG numerical normal flux on the interface $\Gamma_I$ is defined as
\begin{equation}\label{eq:PoissonFluxCG}
- \bn_I {\cdot} \reallywidehat{\grad\CG{u}} :=
- \bn_I {\cdot} \grad\CG{u} + \dfrac{\gamma}{h}(\CG{u} - \HDG{\hu}) ,
\end{equation}
where $h$ denotes the characteristic element size of the mesh discretization on $\Gamma_I$ and $\gamma$ is a \emph{sufficiently large} positive parameter, commonly used to enforce coercivity of the discrete bilinear form in CG discretizations with Nitsche's imposition of essential boundary conditions~\cite{SFM-FH:04}.

Exploiting the definition of the CG numerical normal flux in Equation~\eqref{eq:PoissonFluxCG}, the discrete form of the coupled CG-HDG solver given by Equations~\eqref{eq:PoissonCG}-\eqref{eq:PoissonHDGdiscreteGlobal}-\eqref{eq:couplingConditions} is: find $(\CG{u}, \HDG{\hu}) \in \Vh(\CG{\Omega}) {\times} \HWh(\Gamma\cup\Gamma_N\cup\Gamma_I)$ such that
\begin{subequations}\label{eq:PoissonCoupled}
\begin{align}
&\begin{aligned}
(\grad v, \grad \CG{u})_{\CG{\Omega}} {-}\langle  v, \bn {\cdot} \grad & \CG{u} \rangle_{\Gamma_I} {-} \langle \bn {\cdot} \grad v, \CG{u} \rangle_{\Gamma_I} {+} \langle  v, \tfrac{\gamma}{h}\CG{u} \rangle_{\Gamma_I} 
\\
&{=} (v, f)_{\CG{\Omega}} {+} \langle  v,t \rangle_{\Gamma_N} {-} \langle \bn {\cdot} \grad v, \HDG{\hu} \rangle_{\Gamma_I} {+} \langle  v, \tfrac{\gamma}{h} \HDG{\hu} \rangle_{\Gamma_I} ,
\end{aligned} 
\label{eq:PoissonCoupledCG} \\
&\begin{aligned}
\sum_{e=1}^{\numel} \Bigl\{ 
\langle \hv, \bn {\cdot} \HDG{\bq} \rangle_{\partial\Omega_e\setminus\Gamma_D} {+}& \langle \hv, \tau \HDG{u} \rangle_{\partial\Omega_e\setminus\Gamma_D} {-} \langle \hv, \tau \HDG{\hu} \rangle_{\partial\Omega_e\setminus\Gamma_D} {+} \langle \hv, \tfrac{\gamma}{h} \HDG{\hu} \rangle_{\partial\Omega_e\cap\Gamma_I} \Bigr\} 
\\
&{=} \sum_{e=1}^{\numel} \Bigl\{ {-} \langle \hv, t \rangle_{\partial\Omega_e\cap\Gamma_N} {-} \langle  \hv, \bn {\cdot} \grad \CG{u} \rangle_{\partial\Omega_e\cap\Gamma_I} {+} \langle \hv , \tfrac{\gamma}{h} \CG{u} \rangle_{\partial\Omega_e\cap\Gamma_I} \Bigr\} ,
\end{aligned}
\label{eq:PoissonCoupledHDG}
\end{align}
\end{subequations}
for all $(v, \hv) \in \Vh_0(\CG{\Omega}) {\times} \HWh(\Gamma\cup\Gamma_N\cup\Gamma_I)$.
Note that the third and fourth terms on the left-hand side of Equation~\eqref{eq:PoissonCoupledCG} respectively enforce the symmetry and the coercivity of the resulting CG bilinear form, being $\gamma$ the above mentioned Nitsche's parameter~\cite{SFM-FH:04}. An alternative approach which does not require introducing the penalty term to guarantee the stability of the resulting numerical method is presented in~\cite{Baumann-OBB-98} and relies on a nonsymmetric formulation of the discrete equation.
It is worth noticing that in Equation~\eqref{eq:PoissonCoupledCG} $\bn$ refers to the outer normal to the domain $\CG{\Omega}$, whereas in Equation~\eqref{eq:PoissonCoupledHDG} $\bn$ is the outer normal to the element $\Omega_e \subset \HDG{\Omega}$.

The linear system arising from the discretization of the coupled problem in Equation~\eqref{eq:PoissonCoupled} is
\begin{equation}\label{eq:PoissonLinearSystem}
\begin{bmatrix}
\CG{\mat{K}}          & \KInt \\
\KInt^T & \HDG{\mat{K}} \\
\end{bmatrix}
\begin{bmatrix}
\CG{\vect{u}} \\
\HDG{\vect{\hu}}
\end{bmatrix}
=
\begin{bmatrix}
\CG{\vect{f}} \\
\HDG{\vect{f}}
\end{bmatrix} ,
\end{equation}
where the matrices $\CG{\mat{K}}$ and $\HDG{\mat{K}}$ are symmetric and feature the usual structure of the matrices of the CG and HDG global problem, respectively, whereas $\KInt$ is responsible for the coupling at the interface and stems from the last two terms on the right-hand side of Equations~\eqref{eq:PoissonCoupled}.

It is worth noticing that despite the structure of the block matrix in Equation~\eqref{eq:PoissonLinearSystem}, which is similar to the one presented in~\cite{SFM-PTFM-19}, its construction is extremely different.
The strategies described in~\cite{SFM-PTFM-19} couple CG and HDG discretizations at the elemental level. More precisely, they rely on the introduction of an appropriate projection operator to define the trace of the numerical flux and to impose the Dirichlet boundary condition on $\Gamma_I$ in the HDG local problems. This leads to the intrusive modification of either the block matrix of the elemental problem or the one of the global problem.
On the contrary, the proposed coupling strategy solely relies on the hybrid variable $\HDG{\hu}$ both to impose the Dirichlet boundary condition in the CG problem, see Equation~\eqref{eq:PoissonCoupledCG}, and the Neumann one involving the numerical flux in the HDG global problem, see Equation~\eqref{eq:PoissonCoupledHDG}. 
Moreover, the HDG local problems are not affected by the coupling and the implementation of the resulting strategy is thus minimally-intrusive with respect to existing CG and HDG solvers.

\section{CG-HDG coupling for linear elastic problems}
\label{sc:CG-HDG-Elasticity}

In this section, the coupling strategy is presented for the linear elasticity equation. 
%
%
%
The goal is to exploit in a minimally-intrusive way both the computational efficiency of CG approximations in presence of compressible materials and the accuracy and robustness of HDG discretizations in the nearly incompressible limit.

\subsection{Governing equations}
\label{sc:ElasticityPbStatement}

First, recall the governing equations of a linear elastic material
\begin{equation}\label{eq:Elasticity}
\left\lbrace
\begin{aligned}
-\grad {\cdot} \stress &= \bm{f}          &&\text{in $\Omega$,}   \\
\stress &= \stress^T,          &&\text{in $\Omega$,}   \\
\bu &= \bu_D                                        &&\text{on $\Gamma_D$,} \\
\bn {\cdot} \stress &= \bm{t}                    &&\text{on $\Gamma_N$,}
\end{aligned}
\right.
\end{equation}
where $\bu$ is the unknown displacement field and $\stress$ is the Cauchy stress tensor. Equation~\eqref{eq:Elasticity} enforces the equilibrium of linear and angular momentum of an elastic structure $\Omega$ subject to a volume force $\bm{f}$, a tension $\bm{t}$ on the boundary $\Gamma_N$ and an imposed displacement $\bu_D$ on $\Gamma_D$.
For a linear elastic material, Hooke's law describes the relationship between the Cauchy stress tensor $\stress$ and the linearized strain rate tensor $\defo(\bu) {:=} (\grad \bu {+} \grad\bu^T)/2$, namely
\begin{equation}\label{eq:elasticityConstitutiveLaw}
\stress = \frac{E}{1 + \nu} \defo(\bu) + \frac{E \nu}{(1 + \nu)(1 - 2 \nu)} \tr(\defo(\bu)) \Insd, 
\end{equation}
where $\Insd$ is the $\nsd {\times} \nsd$ identity matrix, $\tr(\cdot)$ is the trace operator and $(E,\nu)$ are the Young's modulus and the Poisson's ratio describing the mechanical properties of the material under analysis.
Henceforth, the material is assumed to be homogeneous and isotropic inside the subdomains $\CG{\Omega}$ and $\HDG{\Omega}$. Thus, the above mentioned material coefficients depend neither on the spatial coordinate nor on the direction of the main strains. 

Following~\cite{FishBelytschko2007}, Equation~\eqref{eq:Elasticity} is rewritten exploiting the Voigt notation for second-order symmetric tensors. The rationale of this approach is to enforce the symmetry of $\stress$ and $\defo(\bu)$ pointwise by storing solely the $\msd {:=} \nsd(\nsd {+} 1)/2$ nonredundant components of the tensor, namely
\begin{equation*}
\stressV := \begin{cases}
\bigl[\sigma_{xx} ,\; \sigma_{yy} ,\; \tau_{xy} \bigr]^T
&\text{in 2D,} \\
\bigl[\sigma_{xx} ,\; \sigma_{yy} ,\; \sigma_{zz} ,\; \tau_{xy} ,\; \tau_{xz} ,\; \tau_{yz} \bigr]^T
&\text{in 3D,} 
\end{cases}
\end{equation*}
and
\begin{equation*}
\strainV := \begin{cases}
\bigl[\varepsilon_{xx} ,\; \varepsilon_{yy} ,\; \gamma_{xy} \bigr]^T
&\text{in 2D,} \\
\bigl[\varepsilon_{xx} ,\; \varepsilon_{yy} ,\; \varepsilon_{zz} ,\; \gamma_{xy} ,\; \gamma_{xz} ,\; \gamma_{yz} \bigr]^T
&\text{in 3D.} 
\end{cases}
\end{equation*}

The strain rate tensor is thus rewritten as $\strainV {=} \gradS \bu$, where the symmetric gradient operator is expressed in matrix form by defining $\gradS \in \RR^{\msd {\times} \nsd}$ such that
\begin{equation} \label{eq:symmGrad}
%
\gradS :=\begin{cases}
\begin{bmatrix}
\displaystyle\frac{\partial}{\partial x} & 0 & \displaystyle\frac{\partial}{\partial y} \\
0 & \displaystyle\frac{\partial}{\partial y} & \displaystyle\frac{\partial}{\partial x}
\end{bmatrix}^T
&\text{in 2D,} \\
\begin{bmatrix}
\displaystyle\frac{\partial}{\partial x} & 0 & 0 & \displaystyle\frac{\partial}{\partial y} & \displaystyle\frac{\partial}{\partial z} & 0 \\
0 & \displaystyle\frac{\partial}{\partial y} & 0 & \displaystyle\frac{\partial}{\partial x} & 0 & \displaystyle\frac{\partial}{\partial z} \\
0 & 0 & \displaystyle\frac{\partial}{\partial z} & 0 & \displaystyle\frac{\partial}{\partial x} & \displaystyle\frac{\partial}{\partial y}
\end{bmatrix}^T
&\text{in 3D,} 
\end{cases}
\end{equation}

Moreover, the constitutive law is expressed in matrix form as $\stressV {=} \bD \strainV$, where the fourth-order elasticity tensor linking $\stress$ and $\defo(\bu)$, see Equation~\eqref{eq:elasticityConstitutiveLaw}, is rewritten by means of the $\msd {\times} \msd$ matrix
\begin{equation} \label{eq:LawVoigt}
\bD :=\begin{cases}
\lambda
\begin{bmatrix}
1{+}(1{-}\vartheta)\nu & \nu & 0 \\
\nu & 1{+}(1{-}\vartheta)\nu & 0 \\
0 & 0 & \displaystyle\frac{1{-}\vartheta\nu}{2}
\end{bmatrix}
&\text{in 2D,} \\[2em]
\lambda
\begin{bmatrix}
1{-}\nu & \nu & \nu & \\
\nu & 1{-}\nu & \nu & \bm{0}_{\nsd} \\
\nu & \nu & 1{-}\nu & \\
& \bm{0}_{\nsd} & & \displaystyle\frac{1{-}2\nu}{2}\Insd
\end{bmatrix}
&\text{in 3D,} 
\end{cases}
\end{equation}
where 
\begin{equation} \label{eq:lambda}
\lambda :=\begin{cases}
\displaystyle\frac{E}{(1+\nu)(1-\vartheta\nu)}
&\text{in 2D,} \\[1em]
\displaystyle\frac{E}{(1+\nu)(1-2\nu)}
&\text{in 3D,} 
\end{cases}
\end{equation}
and the parameter $\vartheta$ denotes either a plane stress model ($\vartheta {=} 1$) or a plane strain model ($\vartheta=2$) in 2D.

In a similar fashion, the Neumann boundary condition is formulated as $\bN^T \stressV {=} \bm{t}$ by introducing the $\msd {\times} \nsd$ matrix describing the outer unit normal vector to the boundary
\begin{equation} \label{eq:normalVoigt}
\bN :=\begin{cases}
\begin{bmatrix}
n_x & 0 & n_y \\
0 & n_y & n_x
\end{bmatrix}^T
&\text{in 2D,} \\
\begin{bmatrix}
n_x & 0 & 0 & n_y & n_z & 0\\
0 & n_y & 0 & n_x & 0 & n_z \\
0 & 0 & n_z & 0 & n_x & n_y
\end{bmatrix}^T
&\text{in 3D.} 
\end{cases}
\end{equation}

Hence, the linear elastic problem in Equation~\eqref{eq:Elasticity}-\eqref{eq:elasticityConstitutiveLaw} is rewritten using Voigt notation as 
\begin{equation} \label{eq:elasticitySystemVoigt}
\left\{\begin{aligned}
-\gradS^T \stressV  &= \bm{f}       &&\text{in $\Omega$,}\\
\stressV &= \bD \strainV      &&\text{in $\Omega$,}\\
\bu &= \bu_D  &&\text{on $\Gamma_D$,}\\
\bN^T \stressV &= \bm{t}        &&\text{on $\Gamma_N$.}\\
\end{aligned}\right.
\end{equation}

\subsection{Hybrid coupling based on Nitsche's method}
\label{sc:ElasticityCoupling}

Following the rationale discussed in Section~\ref{sc:CG-HDG-Poisson}, the strong form of the coupled CG-HDG problem is introduced. More precisely, in the CG subdomain the classical formulation discussed in~\cite{FishBelytschko2007} holds
\begin{equation}\label{eq:ElasticityCG}
\left\lbrace
\begin{aligned}
- \gradS^T \bD \gradS \CG{\bu} &= \bm{f} &&\text{in $\CG{\Omega}$,}  \\
\CG{\bu} &= \bu_D                                                        &&\text{on $\Gamma_D$,}\\
\bN^T \bD \gradS \CG{\bu} &= \bm{t}     &&\text{on $\Gamma_N$,}
\end{aligned}
\right.
\end{equation}
whereas in each element $\Omega_e, \, e {=} 1,\ldots,\numel$ of the HDG subdomain, the recently proposed HDG-Voigt formulation~\cite{RS-SGKH:18,MG-GKSH:18} is considered
\begin{equation}\label{eq:ElasticityHDG}
\left\lbrace
\begin{aligned}
\HDG{\bL} + \bDHalf \gradS \HDG{\bu} &= \bm{0}                    &&\text{in $\Omega_e \subset \HDG{\Omega}$,} \\
\gradS^T \bDHalf \HDG{\bL} &= \bm{f}                             &&\text{in $\Omega_e \subset \HDG{\Omega}$,} \\
\HDG{\bu} &= \bu_D                                                            &&\text{on $\partial\Omega_e \cap \Gamma_D$,} \\
\HDG{\bu} &= \HDG{\bhu}                                                            &&\text{on $\partial\Omega_e \setminus \Gamma_D$,} \\
\bN^T \bDHalf \HDG{\bL} &= - \bm{t}                           &&\text{on $\partial\Omega_e \cap \Gamma_N$,} \\
\jump{ \HDG{\bu} {\otimes} \bn}  &= \bm{0}              &&\text{on $\Gamma$,}   \\
\jump{ \bN^T \bDHalf \HDG{\bL} } &= \bm{0} &&\text{on $\Gamma$.}
\end{aligned}
\right.
\end{equation}
The two methods are thus coupled by a set of conditions equivalent to the ones introduced in Equation~\eqref{eq:couplingConditions}, namely
\begin{equation}\label{eq:ElasticityCoupling}
\left\lbrace
\begin{aligned}
\CG{u} &= \HDG{\hu}                                                   &&\text{on $\Gamma_I$,} \\
-\bN_I^T \bD \gradS \CG{\bu} &= \bN_I^T \bDHalf\HDG{\bL}  &&\text{on $\Gamma_I$,}
\end{aligned}
\right. 
\end{equation}
where the first condition enforces the continuity of the displacement field and the second one the equilibrium of the normal traction at the interface, $\bN_I$ being the outer normal to the domain $\HDG{\Omega}$.
Henceforth, and unless in case of ambiguity, on the interface $\Gamma_I$ the normal $\bN$ is assumed to be the outer direction to the subdomain under analysis, that is $\bN {=} \bN_I$ for $\HDG{\Omega}$ and $\bN {=} {-} \bN_I$ for $\CG{\Omega}$.

As previously discussed, the discrete form of the coupled problem is obtained by the weak form of the CG problem with Dirichlet boundary conditions imposed via Nitsche's method on $\Gamma_I$ and by the HDG global problem, that is: find $(\CG{\bu},\HDG{\bhu}) \in [\Vh(\CG{\Omega})]^{\nsd} {\times}$ $[\HWh(\Gamma \cup \Gamma_N \cup \Gamma_I)]^{\nsd}$ such that
\begin{subequations}\label{eq:ElasticityCoupled}
\begin{align}
&\begin{aligned}
(\gradS \bv, \bD \gradS \CG{\bu})_{\CG{\Omega}} {-}& \langle  \bv, \bN^T \bD \gradS \CG{\bu} \rangle_{\Gamma_I} {-} \langle \bN^T \bD \gradS \bv, \CG{\bu} \rangle_{\Gamma_I} {+} \langle  \bv, \tfrac{\gamma}{h}\CG{\bu} \rangle_{\Gamma_I} 
\\
&{=} (\bv, \bm{f})_{\CG{\Omega}} {+} \langle  \bv, \bm{t} \rangle_{\Gamma_N} {-} \langle \bN^T \bD \gradS \bv, \HDG{\bhu} \rangle_{\Gamma_I} {+} \langle  \bv, \tfrac{\gamma}{h} \HDG{\bhu} \rangle_{\Gamma_I} ,
\end{aligned} 
\label{eq:ElasticityCoupledCG} \\
&\begin{aligned}
\sum_{e=1}^{\numel} \Bigl\{ 
\langle \bhv, \bN^T \bDHalf \HDG{\bL} & \rangle_{\partial\Omega_e\setminus\Gamma_D} {+} \langle \bhv, \tau \HDG{\bu} \rangle_{\partial\Omega_e\setminus\Gamma_D} {-} \langle \bhv, \tau \HDG{\bhu} \rangle_{\partial\Omega_e\setminus\Gamma_D} {+} \langle \bhv, \tfrac{\gamma}{h} \HDG{\bhu} \rangle_{\partial\Omega_e\cap\Gamma_I} \Bigr\} 
\\
&{=} \sum_{e=1}^{\numel} \Bigl\{ {-} \langle \bhv, \bm{t} \rangle_{\partial\Omega_e\cap\Gamma_N} {-} \langle  \hv, \bN^T \bD \gradS \CG{\bu} \rangle_{\partial\Omega_e\cap\Gamma_I} {+} \langle \bhv , \tfrac{\gamma}{h} \CG{\bu} \rangle_{\partial\Omega_e\cap\Gamma_I} \Bigr\} ,
\end{aligned}
\label{eq:ElasticityCoupledHDG}
\end{align}
\end{subequations}
for all $(\bv, \bhv) \in [\Vh_0(\CG{\Omega})]^{\nsd} {\times} [\HWh(\Gamma\cup\Gamma_N\cup\Gamma_I)]^{\nsd}$, where the definition of the trace of the numerical normal flux used for CG on the interface $\Gamma_I$ is
\begin{equation}\label{eq:ElasticityFluxCG}
- \bN_I^T \reallywidehat{\bD \gradS \CG{\bu}} :=
- \bN_I^T \bD \gradS \CG{\bu} + \dfrac{\gamma}{h}(\CG{\bu} - \HDG{\bhu}) ,
\end{equation}
and for HDG on all internal and boundary faces is
\begin{equation} \label{eq:ElasticityFlux}
\bN^T \reallywidehat{\bDHalf \HDG{\bL}} := 
\begin{cases}
\bN^T \bDHalf \HDG{\bL} + \tau (\HDG{\bu} - \bu_D) & \text{on $\partial\Omega_e\cap\Gamma_D$}, \\
\bN^T \bDHalf \HDG{\bL} + \tau (\HDG{\bu} - \HDG{\bhu}) & \text{elsewhere}.
\end{cases}
\end{equation}

As previously remarked, the proposed coupling only affects the HDG discretization in the global problem, whereas the local element-by-element problems present the usual structure of an HDG approximation.
More precisely, the choice of the HDG formulation under analysis relies on its capability to achieve optimal convergence of the stress tensor and superconvergence of the displacement field using equal order polynomial approximation for all the variables.
This is due to the definition of a pointwise symmetric mixed variable, namely the strain rate tensor, via Voigt notation~\cite{RS-SGKH:18,MG-GKSH:18}.
The resulting HDG local problems are: for $e {=} 1,\dots,\numel$ find $(\HDG{\bu},\HDG{\bL}) \in [\Wh(\Omega_e)]^{\nsd} {\times} [\Wh(\Omega_e)]^{\msd}$ such that
\begin{subequations}\label{eq:ElasticityHDGdiscreteLocal}
\begin{align}
- (\bw, \HDG{\bL})_{\Omega_e} {+} (\gradS \bDHalf \bw, \HDG{\bu})_{\Omega_e} 
&= \langle \bN^T \bDHalf \bw, \bu_D \rangle_{\partial\Omega_e\cap\Gamma_D} {+} \langle \bN^T \bDHalf \bw, \bhu^h \rangle_{\partial\Omega_e\setminus\Gamma_D} ,
\\
(\bv , \gradS \bDHalf \HDG{\bL})_{\Omega_e} {+} \langle \bv, \tau \HDG{\bu} \rangle_{\partial\Omega_e}
&= (\bv,\bm{f})_{\Omega_e} {+} \langle \bv, \tau \bu_D \rangle_{\partial\Omega_e\cap\Gamma_D} {+} \langle \bv, \tau \bhu^h \rangle_{\partial\Omega_e\setminus\Gamma_D} ,
\end{align}
\end{subequations}
for all $(\bv,\bw) \in [\Wh(\Omega_e)]^{\nsd} {\times} [\Wh(\Omega_e)]^{\msd}$.

Moreover, given the discrete functional space 
\begin{equation*}
\Wh_{\star}(\Omega) {:=}
\lbrace 
v \in \eltwo(\Omega) : \, v\vert_{\Omega_e} \in \Poly[k+1](\Omega_e) \forall \Omega_e, \, e {=} 1,\ldots,\numel
\rbrace,
\end{equation*}
where $\Poly[k+1](\Omega_e)$ is the space of polynomial functions of complete degree at most $k {+} 1$, a superconvergent approximation $\HDG{\bu}^\star \in [\Wh_{\star}(\Omega)]^{\nsd}$ of the displacement field is computed for each element $\Omega_e, \, e {=} 1,\dots,\numel$ by solving the postprocessed problem
\begin{equation}\label{eq:ElasticityPostProcess}
\left\lbrace
\begin{aligned}
- \gradS^T \bDHalf \gradS \HDG{\bu}^{\star} &= \gradS^T \HDG{\bL} &&\text{in $\Omega_e$,}  \\
\bN^T \bDHalf \gradS \HDG{\bu}^{\star} &= - \bN^T \HDG{\bL}     &&\text{on $\partial\Omega_e$,}
\end{aligned}
\right.
\end{equation}
with the constraint 
\begin{equation}\label{eq:ElasticityConstraintTranslation}
(\HDG{\bu}^\star, 1)_{\Omega_e} = (\HDG{\bu}, 1)_{\Omega_e} ,
\end{equation}
to remove the underdetermination associated with rigid body translations and 
\begin{equation}\label{eq:ElasticityConstraintRotation}
(\gradW \HDG{\bu}^\star, 1)_{\Omega_e} = \ \langle \bT \bu_D, 1 \rangle_{\partial\Omega_e\cap\Gamma_D} {+} \langle \bT \HDG{\bhu}, 1 \rangle_{\partial\Omega_e\setminus\Gamma_D} ,
\end{equation}
to treat rigid body rotations, where the \emph{curl} operator $\gradW \in \RR^{\nrr {\times} \nsd}$ and the tangential direction to the boundary $\bT \in \RR^{\nrr {\times} \nsd}$, with $\nrr {:=} \nsd(\nsd {-} 1)/2$ being the number of rigid rotational body modes, are written in matrix form as
\begin{equation}\label{eq:VoigtCurl}
\gradW {:=}
\left\lbrace
\begin{aligned}
&\begin{bmatrix}
-\partial/\partial y &  \partial/\partial x
\end{bmatrix}
&\text{in 2D,}
\\
&\begin{bmatrix}
                   0 & -\partial/\partial z &  \partial/\partial y \\
 \partial/\partial z &                    0 & -\partial/\partial x \\
-\partial/\partial y &  \partial/\partial x &                    0
\end{bmatrix}
&\text{in 3D,}
\end{aligned}
\right.
\end{equation}
and
\begin{equation}\label{eq:VoigtTangent}
\bT {:=}
\left\lbrace
\begin{aligned}
&\begin{bmatrix}
-n_y &  n_x
\end{bmatrix}
&\text{in 2D,}
\\
&\begin{bmatrix}
                  0 & -n_z &  n_y \\
 n_z &                   0 & -n_x \\
-n_y &  n_x &                   0
\end{bmatrix}
&\text{in 3D.}
\end{aligned}
\right.
\end{equation}

\section{Numerical studies}
\label{sc:NumericalStudies}

\begin{figure}
\centering
\includegraphics[width=0.3\textwidth]{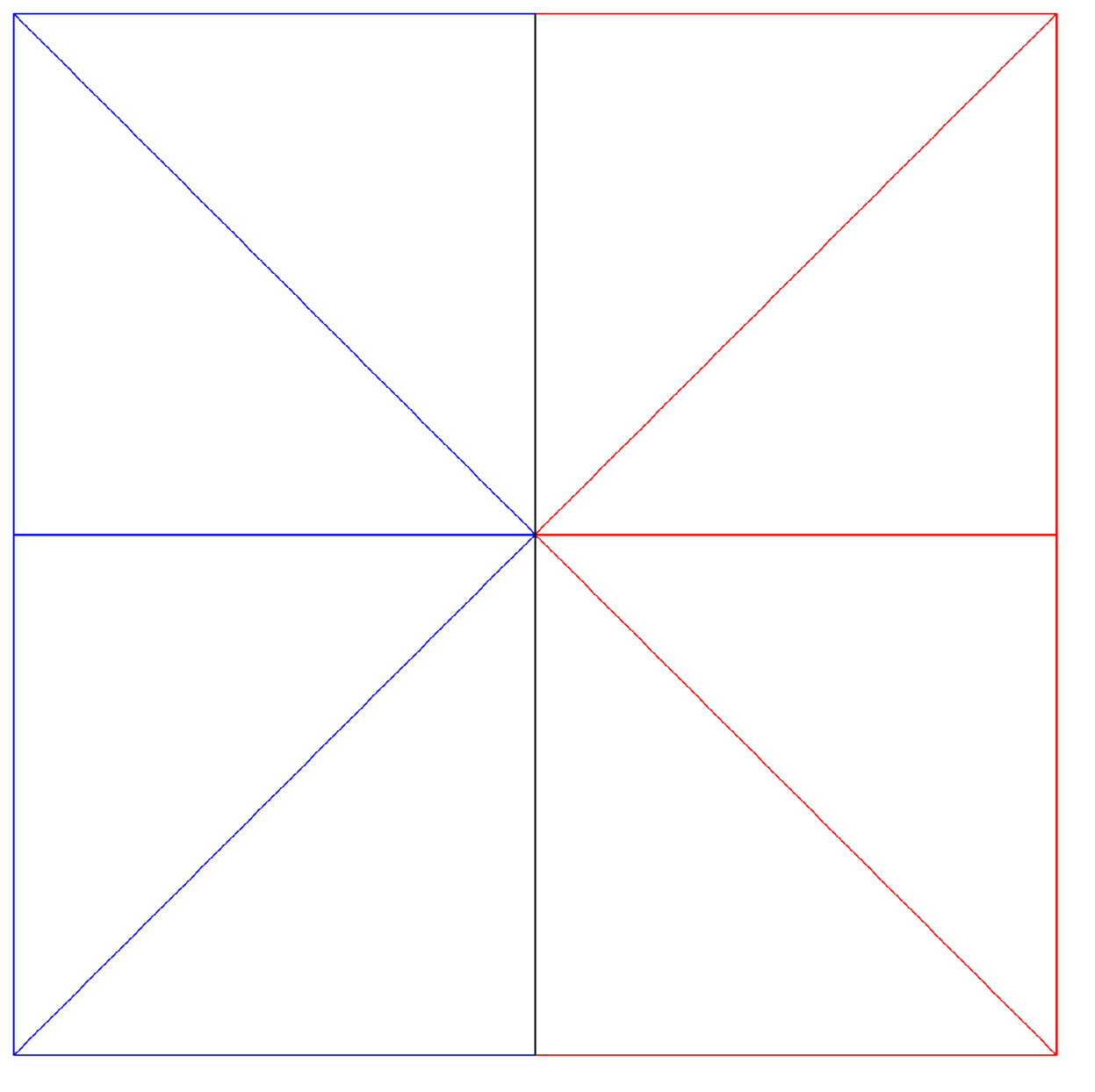}
\hfill
\includegraphics[width=0.3\textwidth]{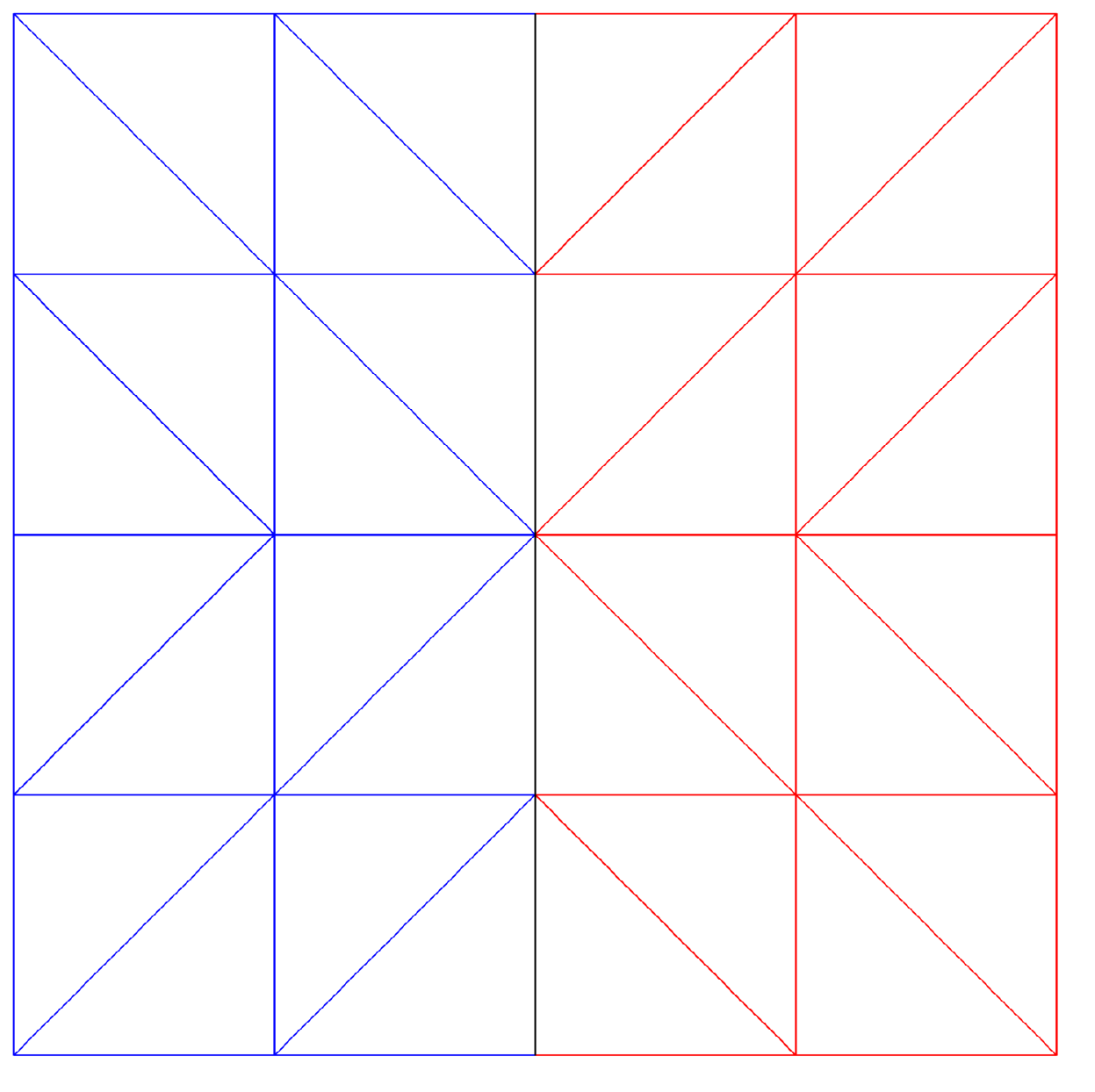}
\hfill
\includegraphics[width=0.3\textwidth]{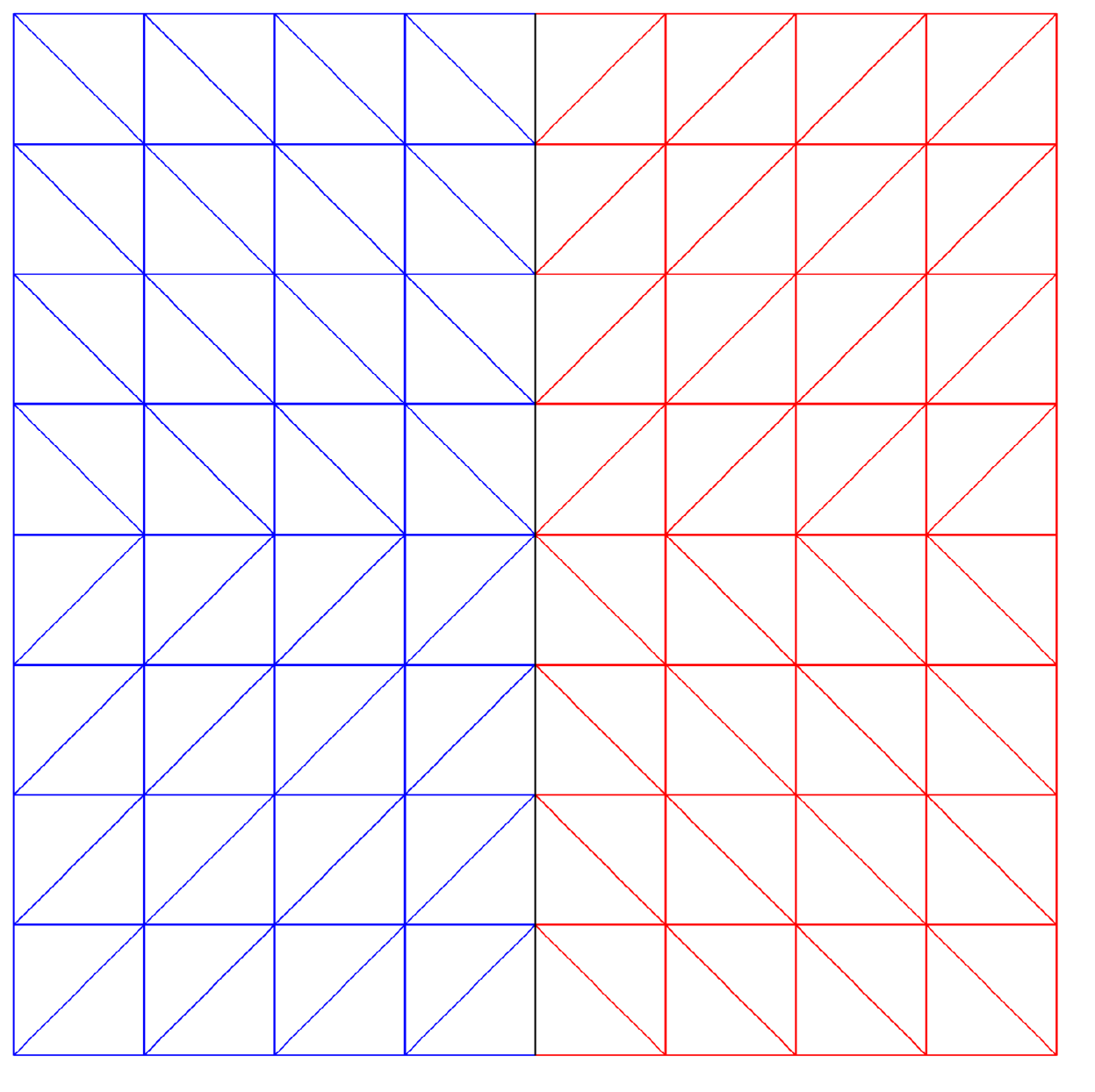}
\caption{First three levels of refinement of the mesh used for the convergence study of the two-dimensional thermal problem. Red: CG subdomain $\CG{\Omega}$. Blue: HDG subdomain $\HDG{\Omega}$. Black: interface $\Gamma_I$.}
\label{fig:poissonMesh}
\end{figure}

In this section, several numerical examples are presented to test the optimal convergence properties of the proposed hybrid CG-HDG coupling.
Special emphasis is devoted to highlight the advantages of this methodology in terms of accuracy, robustness and minimal intrusiveness of its implementation in existing CG and HDG libraries.

\subsection{Two-dimensional thermal problem}
\label{sc:2Dthermal}

The first example considers the problem in Equation~\eqref{eq:PoissonStrongForm}, in two dimensions, to assess optimal convergence and robustness of the coupling approach to the involved parameters.
The computational domain $\Omega {=} [-1,1] {\times} [-1,1]$ is decomposed in two nonoverlapping subdomains, namely $\CG{\Omega} {=} [0,1] {\times} [-1,1]$ and 
$\HDG{\Omega} {=} \Omega \setminus \CG{\Omega}$. The interface is thus identified by $\Gamma_I {=} \left\{ (x,y) \in \RR^2 : \right.$ $\left. x=0 \right\}$.
The domain $\Omega$ is discretized using uniform meshes of triangular elements constructed by means of the mesh generator EZ4U~\cite{2000-IJNME-SH,2001-CNME-SH}.
The first three levels of mesh refinement are presented in Figure~\ref{fig:poissonMesh}. Elements in red (resp., blue) belong to the CG (resp., HDG) subdomain, whereas the interface $\Gamma_I$ is drawn in black.
%
%

The source term is selected so that the analytical solution is
\begin{equation*}
u(x,y) = \cos \left( \frac{\pi}{2} \sqrt{x^2+y^2} \right) ,
\end{equation*}
and Dirichlet boundary conditions, corresponding to the analytical solution, are imposed on $\Gamma_D {=} \partial\Omega$.

\begin{figure}
\centering
\includegraphics[width=0.45\columnwidth]{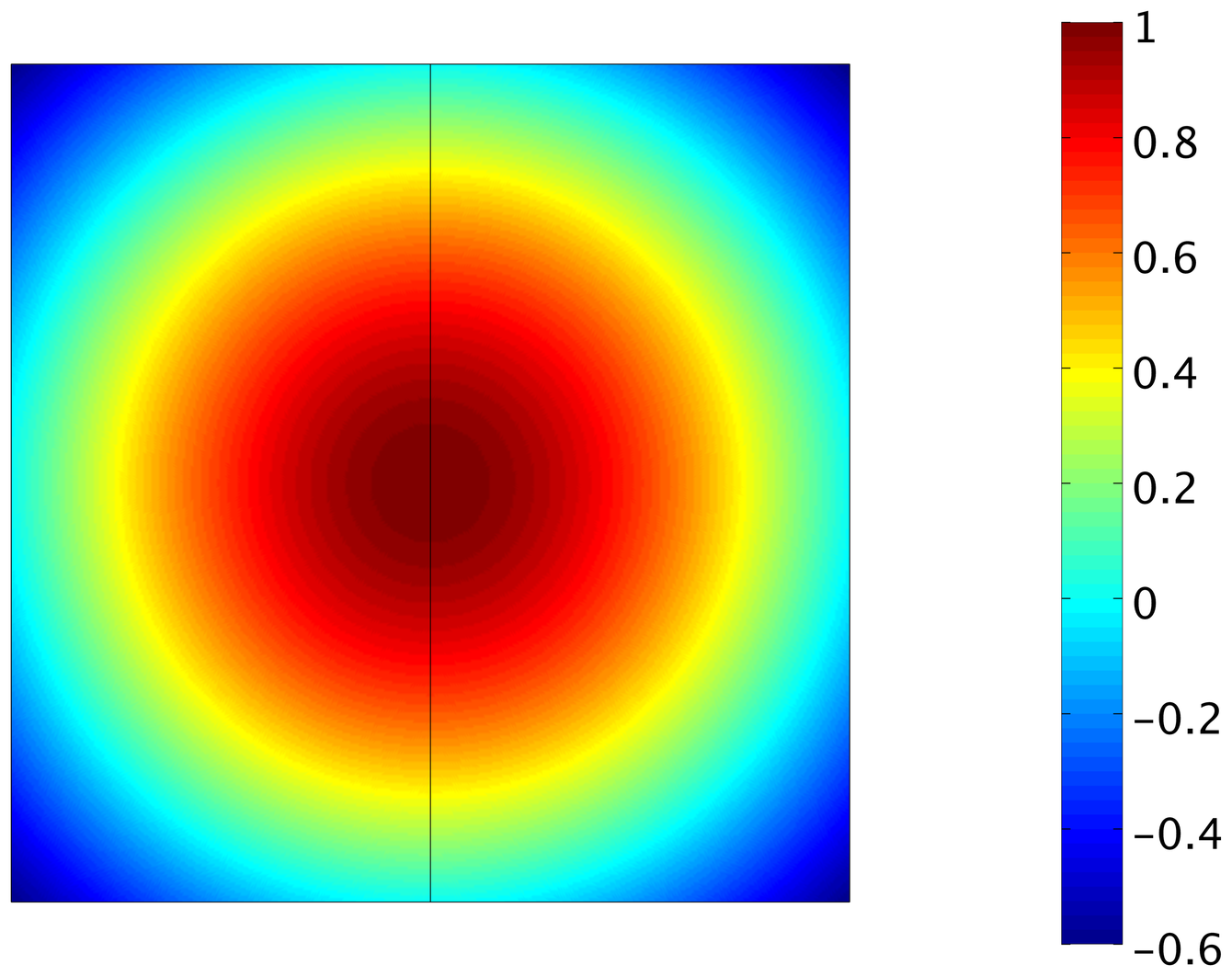}
\caption{Approximate temperature field computed using the coupled CG-HDG strategy on the fifth level of mesh refinement with polynomial of degree $k {=} 5$. The solid black line represents the interface $\Gamma_I$ identifying the HDG solution in the left subdomain and the CG one on the right.}
\label{fig:poissonSolution}
\end{figure}

The solution of the thermal problem computed on the fifth level of refinement of the mesh with polynomial approximation of degree $k {=} 5$ in both the CG and HDG subdomains is displayed in Figure~\ref{fig:poissonSolution}.

\subsubsection{Sensitivity to Nitsche's parameter}
\label{sc:2DthermalSensitivity}

The influence of Nitsche's parameter on the accuracy of the proposed hybrid CG-HDG coupling is investigated.
The HDG stabilization parameter $\tau$ is considered constant in the subdomain $\HDG{\Omega}$ and equal to $10$. It is worth recalling that a stabilization parameter $\tau {=} C \kappa/\ell$, where $\kappa$ is the thermal conductivity, equal to $1$ for the problem under analysis, $\ell$ the characteristic length of the problem and $C$ a positive constant scaling factor, guarantees stability and optimal convergence of the HDG discretization for the Poisson equation~\cite{Jay-CGL:09}.

\begin{figure}
\centering
\includegraphics[width=0.45\columnwidth]{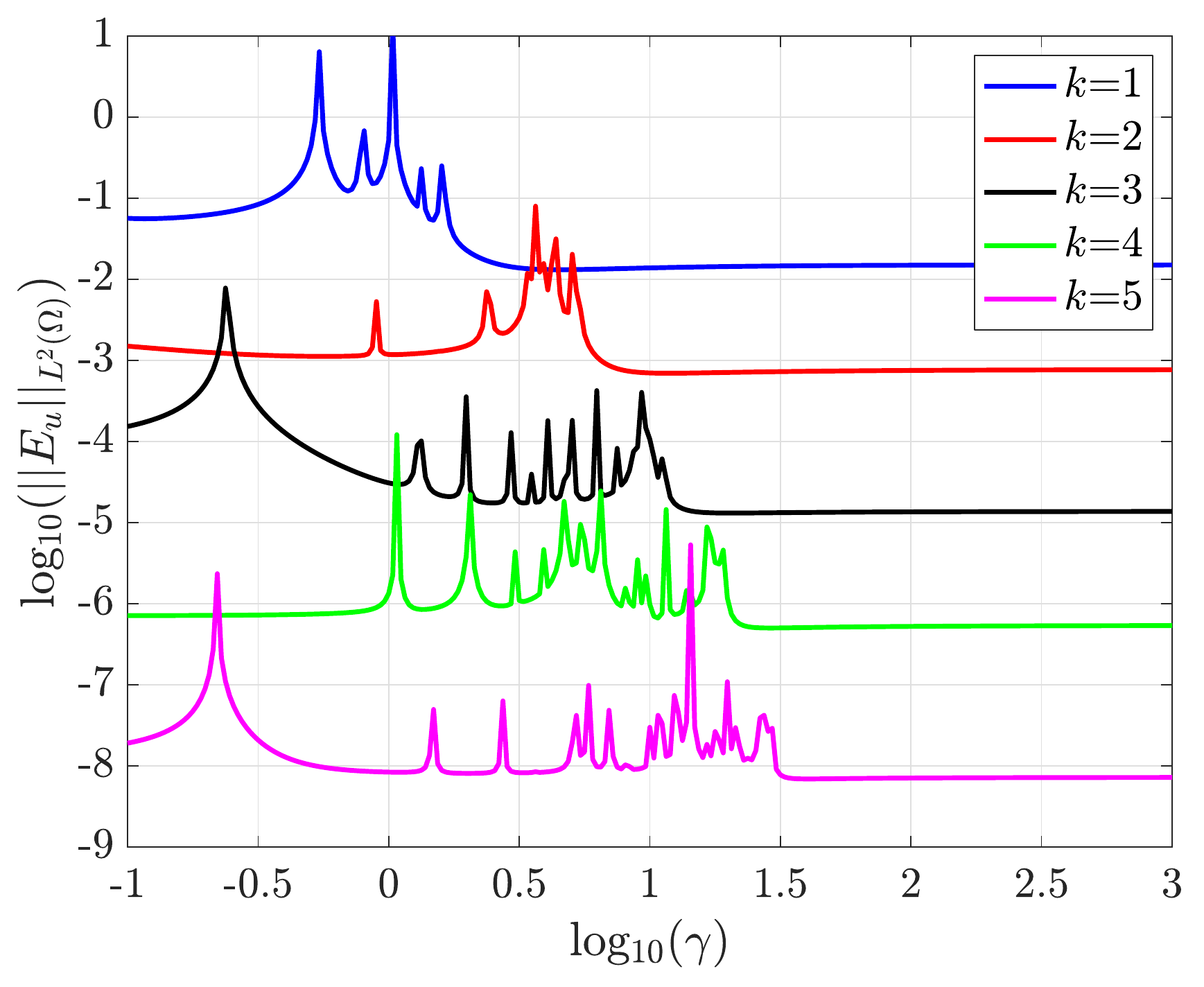}
\caption{Error of the temperature computed using the hybrid CG-HDG coupling in the $\eltwo$ norm on the domain $\Omega$ as a function of Nitsche's parameter $\gamma$.}
\label{fig:poissonNitscheParam}
\end{figure}

Figure~\ref{fig:poissonNitscheParam} displays the evolution of the error of the primal variable $u$ on the whole domain $\Omega$ as a function of Nitsche's parameter $\gamma$ using the third level of mesh refinement and for different degrees of the polynomial approximation.
For low values of Nitsche's parameter, oscillations appear in the error, whereas stability is achieved choosing a \emph{sufficiently large} $\gamma$.
It is worth noticing that the minimum value of $\gamma$ guaranteeing stability of the numerical method varies with the degree of the polynomial approximation.
Henceforth, $\gamma {=} 10^2$ is considered for the following numerical experiments.

\begin{figure}
\centering
\includegraphics[width=0.45\columnwidth]{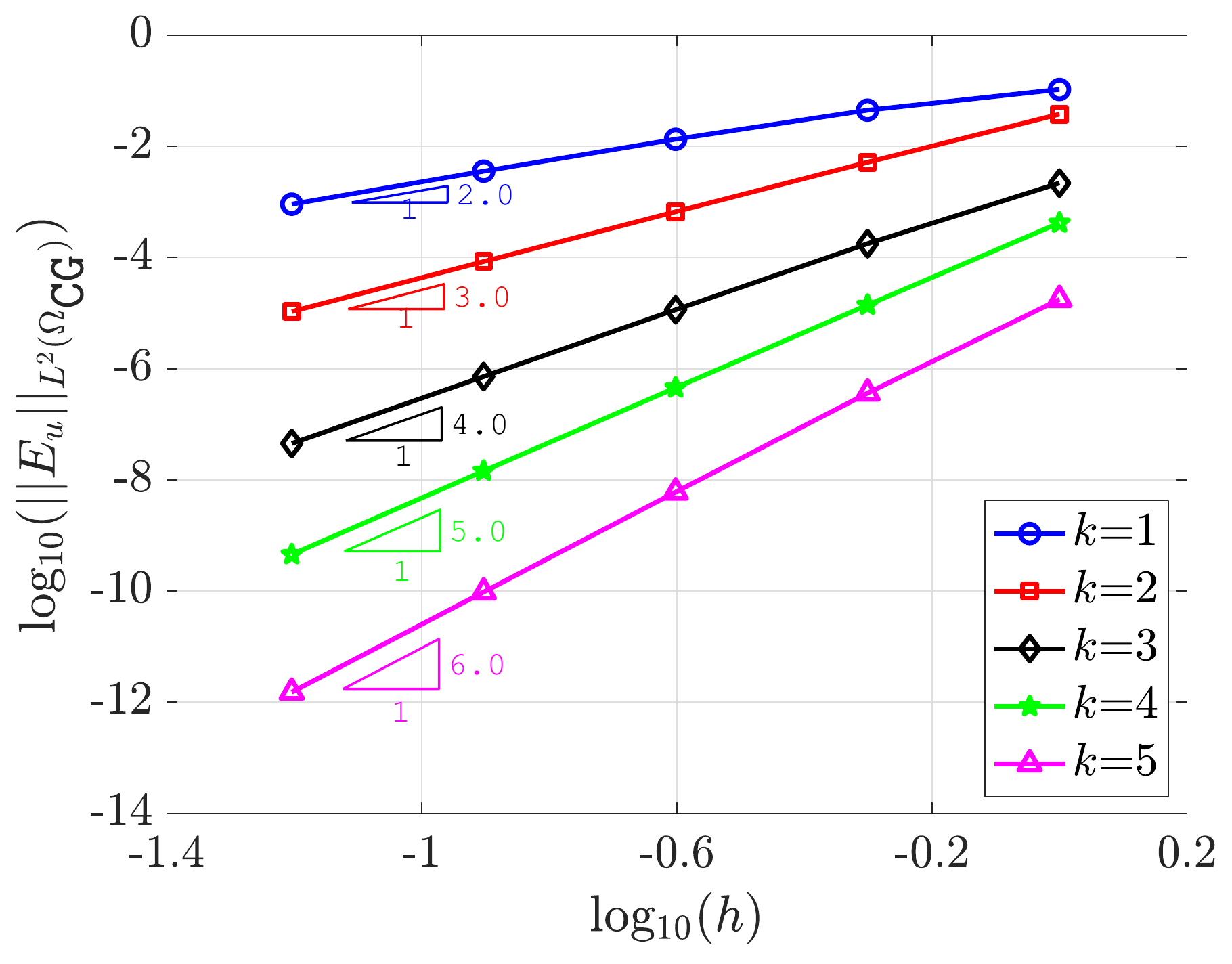}
\hfill
\includegraphics[width=0.45\columnwidth]{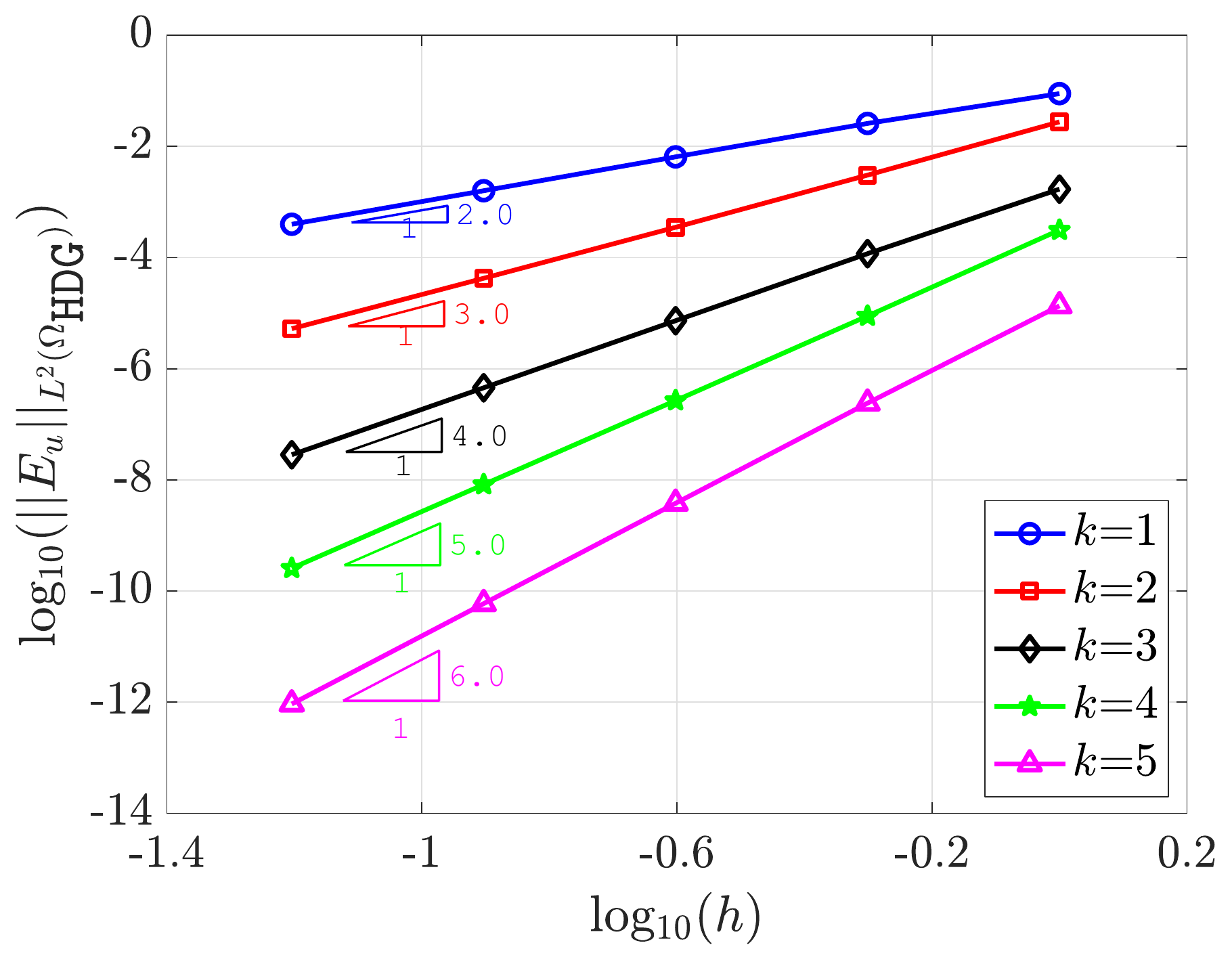}
\caption{$h$-convergence of the error of the temperature computed using the hybrid CG-HDG coupling in the $\eltwo$ norm on the subdomain $\CG{\Omega}$ (left) and $\HDG{\Omega}$ (right).}
\label{fig:poissonConvergenceCoupling}
\end{figure}

\subsubsection{Optimal convergence of the CG-HDG coupling}
\label{sc:2DthermalConvergence}

The convergence of the error of the temperature measured in the $\eltwo$ norm as a function of the characteristic element size $h$ is presented in Figure~\ref{fig:poissonConvergenceCoupling} for polynomial degree of approximation $k {=} 1,\ldots,5$.
Optimal convergence of order $k {+} 1$ is achieved both in the CG subdomain (Fig.~\ref{fig:poissonConvergenceCoupling}, left) and in the HDG one (Fig.~\ref{fig:poissonConvergenceCoupling}, right) by the proposed hybrid CG-HDG coupling.
%
%

Comparable results in terms of accuracy are obtained using the coupling strategy discussed in~\cite{SFM-PTFM-19} and the convergence studies are omitted for the sake of brevity.
It is worth recalling that the main advantage of the proposed hybrid coupling is represented by its minimally-intrusive nature which makes this approach extremely easy to implement in existing CG and HDG libraries.

\subsection{Two-dimensional elastic problem}
\label{sc:2Delastic}

\begin{figure}
\centering
\includegraphics[width=0.3\textwidth]{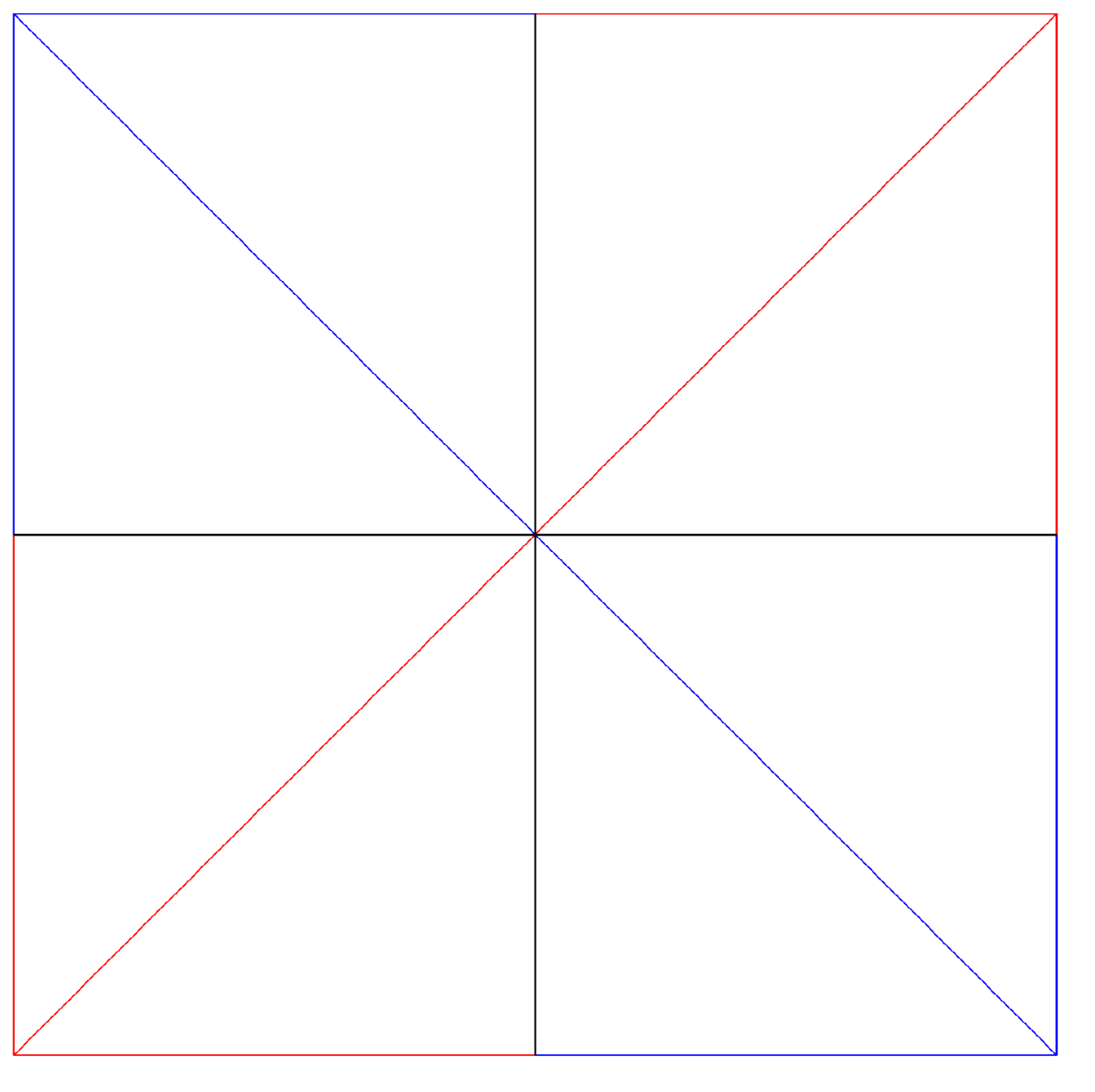}
\hfill
\includegraphics[width=0.3\textwidth]{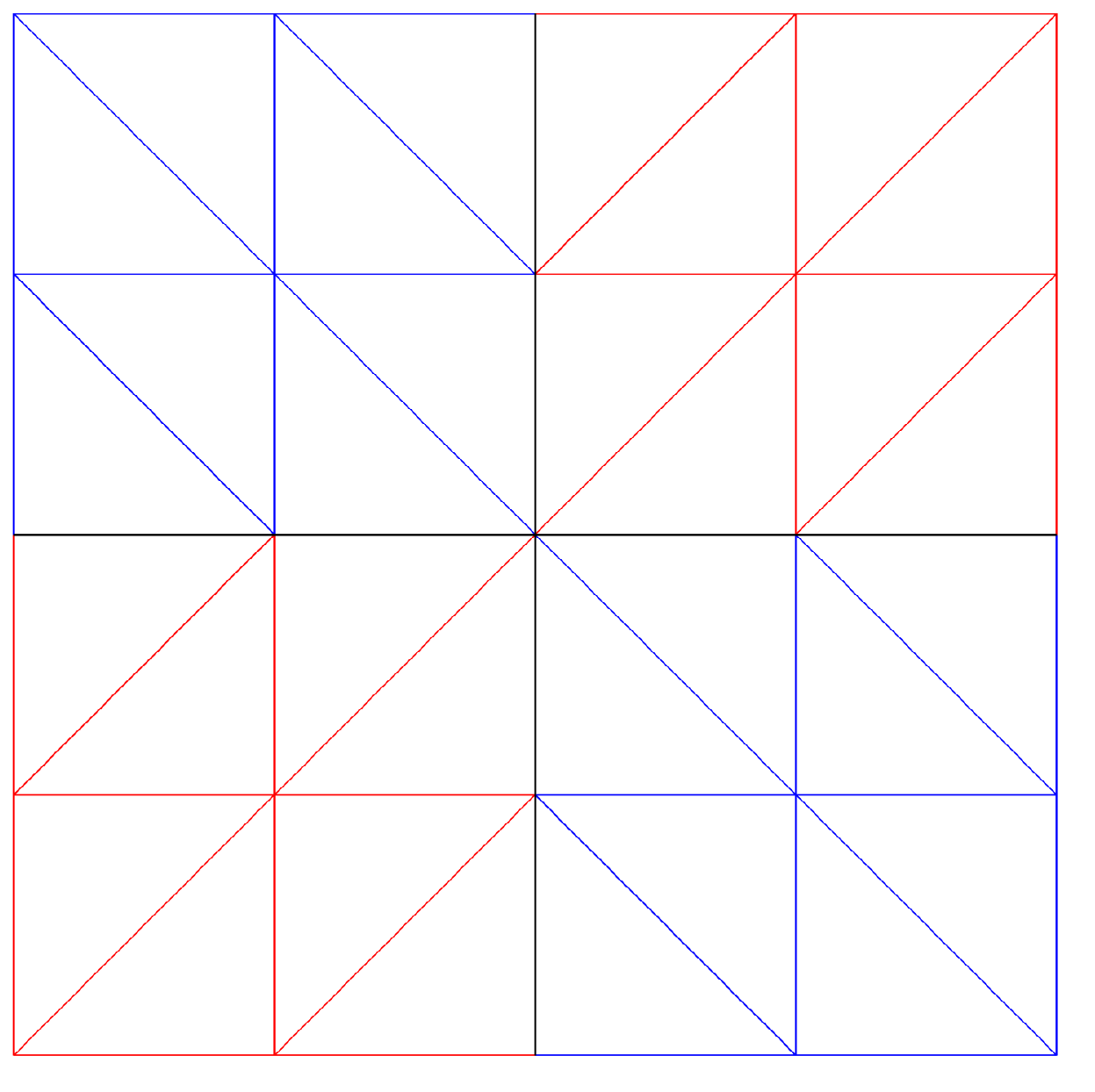}
\hfill
\includegraphics[width=0.3\textwidth]{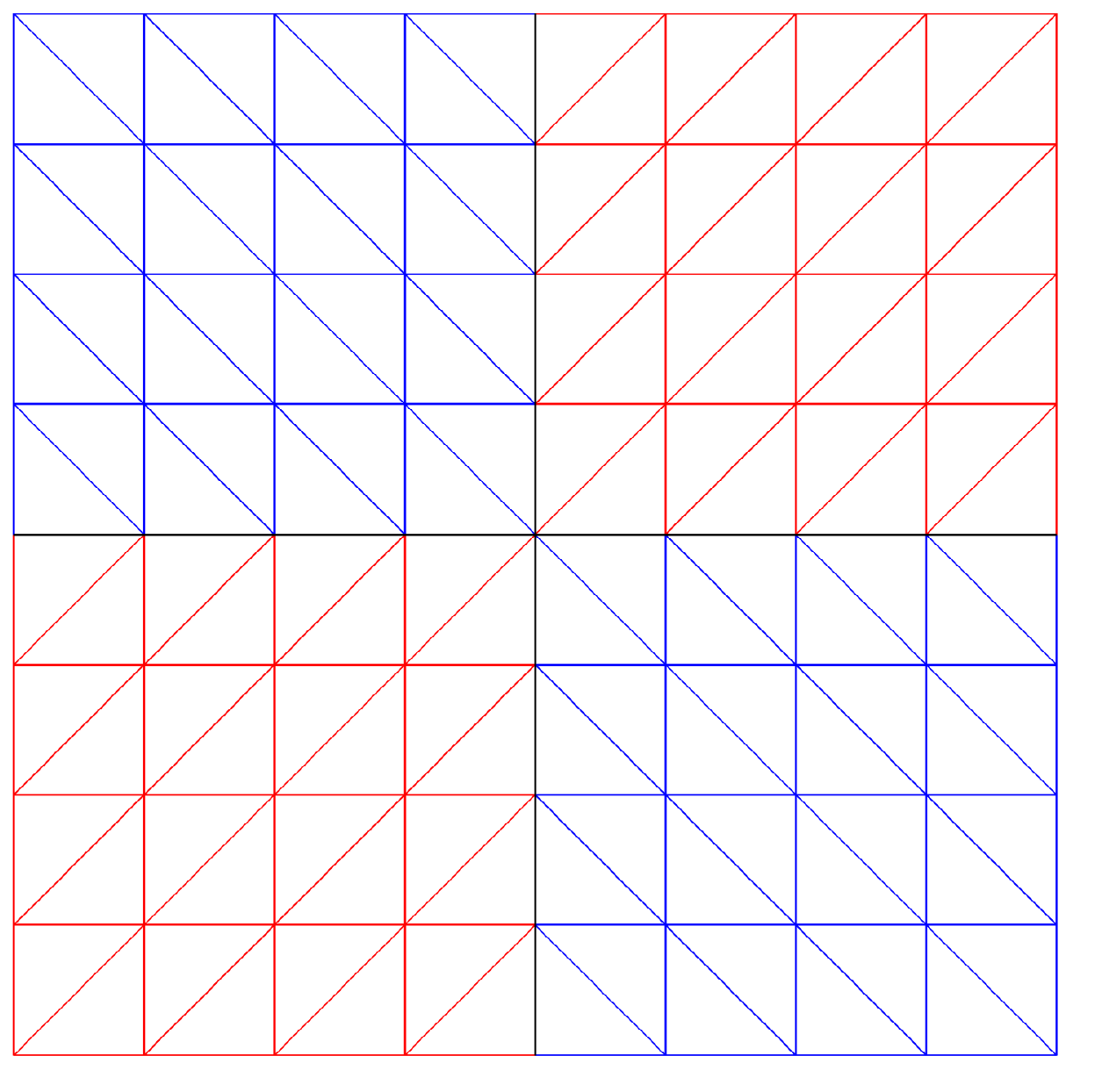}
\caption{First three levels of refinement of the mesh used for the convergence study of the two-dimensional elastic problem. Red: CG subdomain $\CG{\Omega}$. Blue: HDG subdomain $\HDG{\Omega}$. Black: interface $\Gamma_I$.}
\label{fig:elasticityMesh}
\end{figure}

The second example considers the linear elastic problem in Equation~\eqref{eq:elasticitySystemVoigt} with compressible and nearly incompressible materials.
The domain $\Omega {=} [-1,1] {\times} [-1,1]$ is decomposed in two nonoverlapping subdomains, namely $\CG{\Omega} {=} [-1,0] {\times} [-1,0] \cup [0,1] {\times} [0,1]$ and 
$\HDG{\Omega} {=} \Omega \setminus \CG{\Omega}$. The interface is thus identified by $\Gamma_I {=} \left\{ (x,y) \in \RR^2 : \right.$ $\left. x=0 \text{ or } y=0 \right\}$.
The computational meshes are constructed as for the above described thermal problem. The first three levels of mesh refinement are presented in Figure~\ref{fig:elasticityMesh}: the subdomain $\CG{\Omega}$ in displayed in red, $\HDG{\Omega}$ in blue, whereas the interface $\Gamma_I$ is drawn in black.
It is worth noticing that for the problem under analysis both the CG and the HDG subdomains are not simply connected
%

The source term is selected so that the analytical solution of the problem is $\bu {=} (u_x,u_y)$ such that
\begin{align*}\label{eq:elasticitySolution}
u_x(x,y) =&
\frac{2(1 {+}\nu)}{E} \sin(2\pi y) [{-}1 {+} \cos(2\pi x) ] 
+ \frac{(1 {+} \nu)(1 {-} 2\nu)}{(1 {+} \nu)(1 {-} 2\nu) {+} \nu E} xy \sin(\pi x) \sin(\pi y), 
\\
u_y(x,y) =&
\frac{2(1 {+}\nu)}{E} \sin(2\pi x) [1 {-} \cos(2\pi y) ] 
+ \frac{(1 {+} \nu)(1 {-} 2\nu)}{(1 {+} \nu)(1 {-} 2\nu) {+} \nu E} xy \sin(\pi x) \sin(\pi y),
\end{align*}
and Dirichlet boundary conditions, corresponding to the analytical solution, are imposed on $\Gamma_D {=} \partial\Omega$.
This numerical experiment is inspired by the work in~\cite{Lamichhane2009} and considers inhomogeneous mechanical properties in the two subdomains, namely a compressible and stiff material with Young's modulus $E {=} 250$ and Poisson's ratio $\nu {=} 0.3$ in the red domain $\CG{\Omega}$ and a nearly incompressible and soft one such that $E {=} 25$ and $\nu {=} 0.49999$ in the blue domain $\HDG{\Omega}$.

\begin{figure}
\centering
\subfigure[CG]{\includegraphics[width=0.45\textwidth]{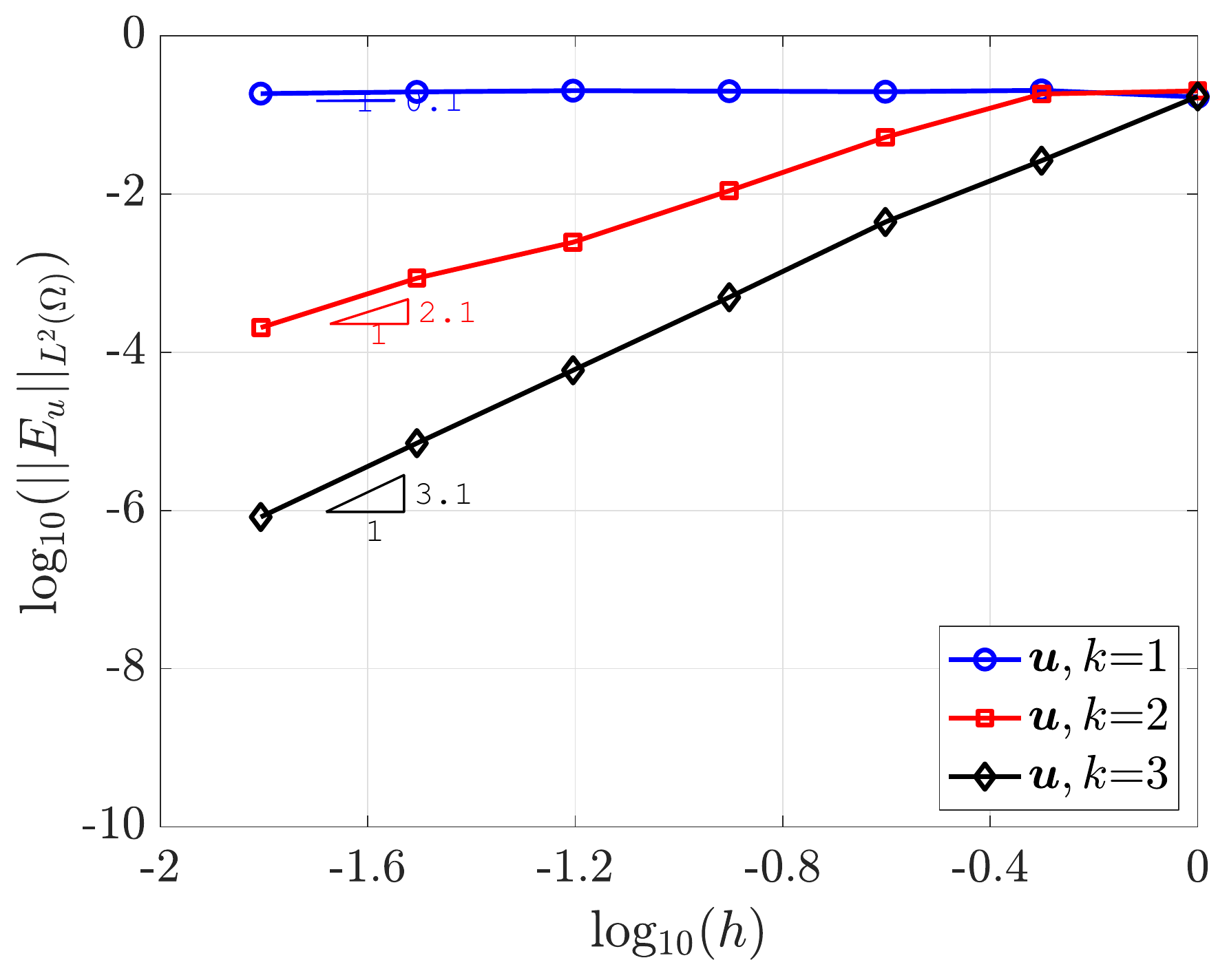}}
\hfill
\subfigure[HDG]{\includegraphics[width=0.45\textwidth]{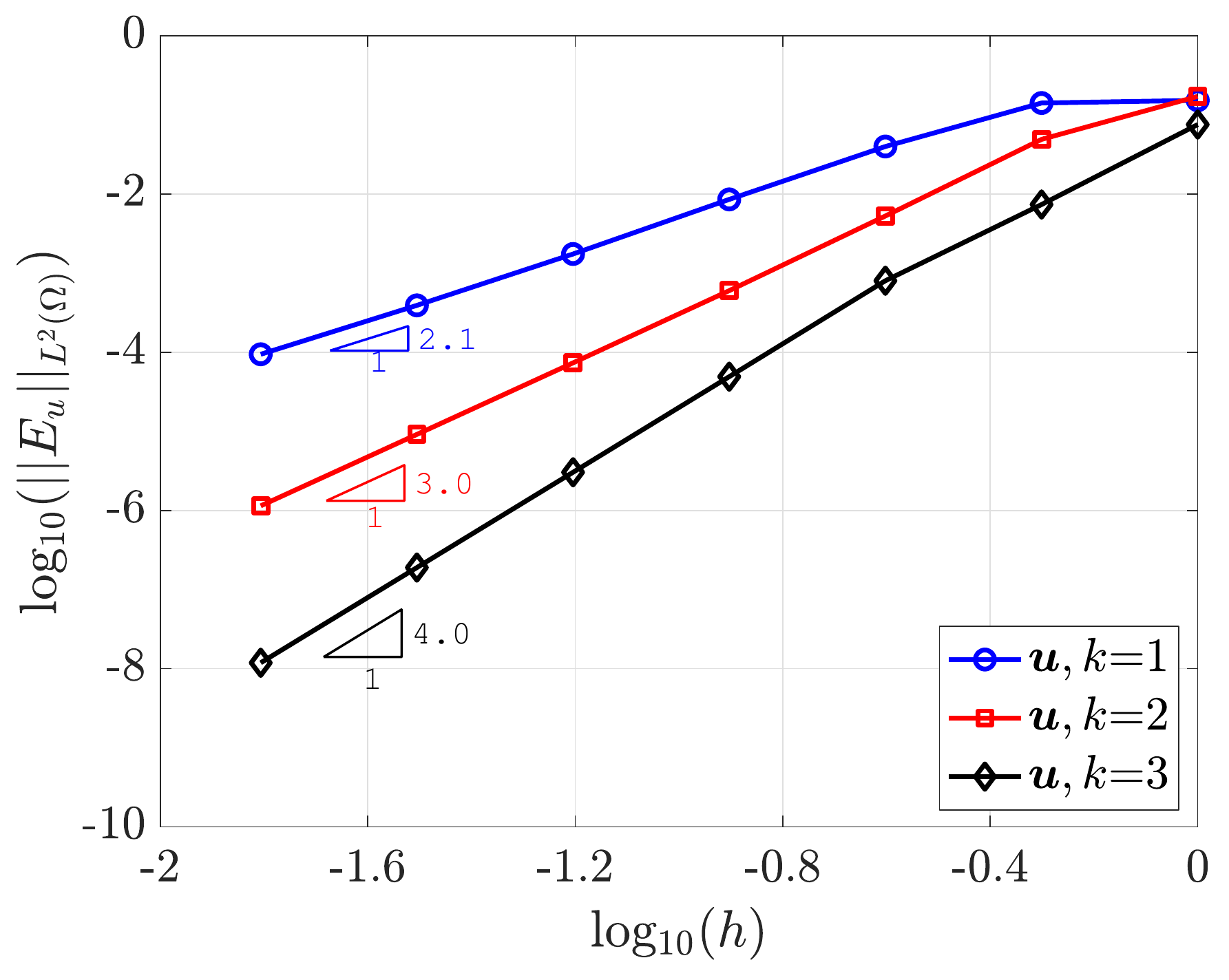}}

\subfigure[$\Kdeg{CG}$-$\Kdeg{HDG}$]{\includegraphics[width=0.45\textwidth]{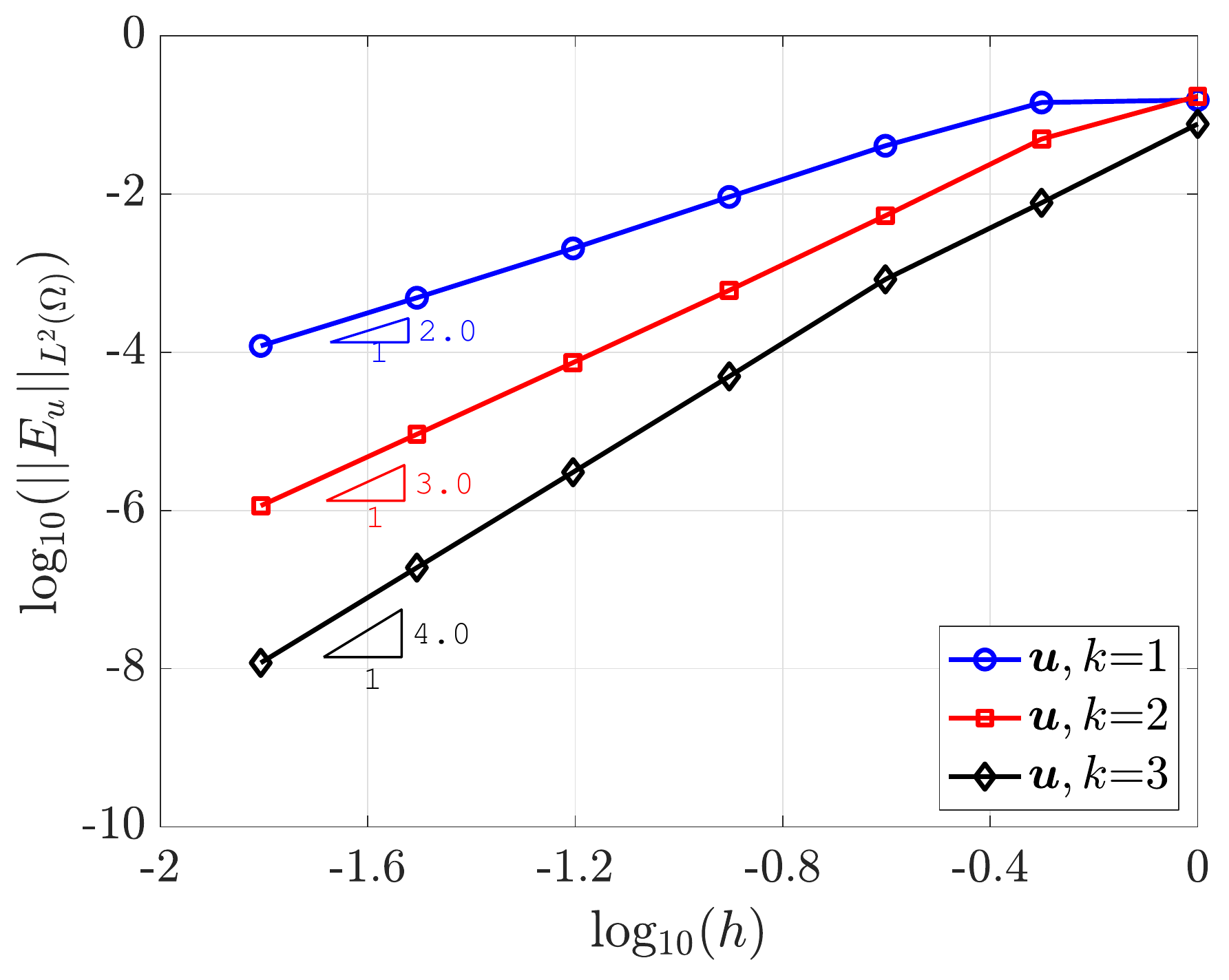}}
\hfill
\subfigure[$\KPdeg{CG}$-$\Kdeg{HDG}$ with postprocess]{\includegraphics[width=0.45\textwidth]{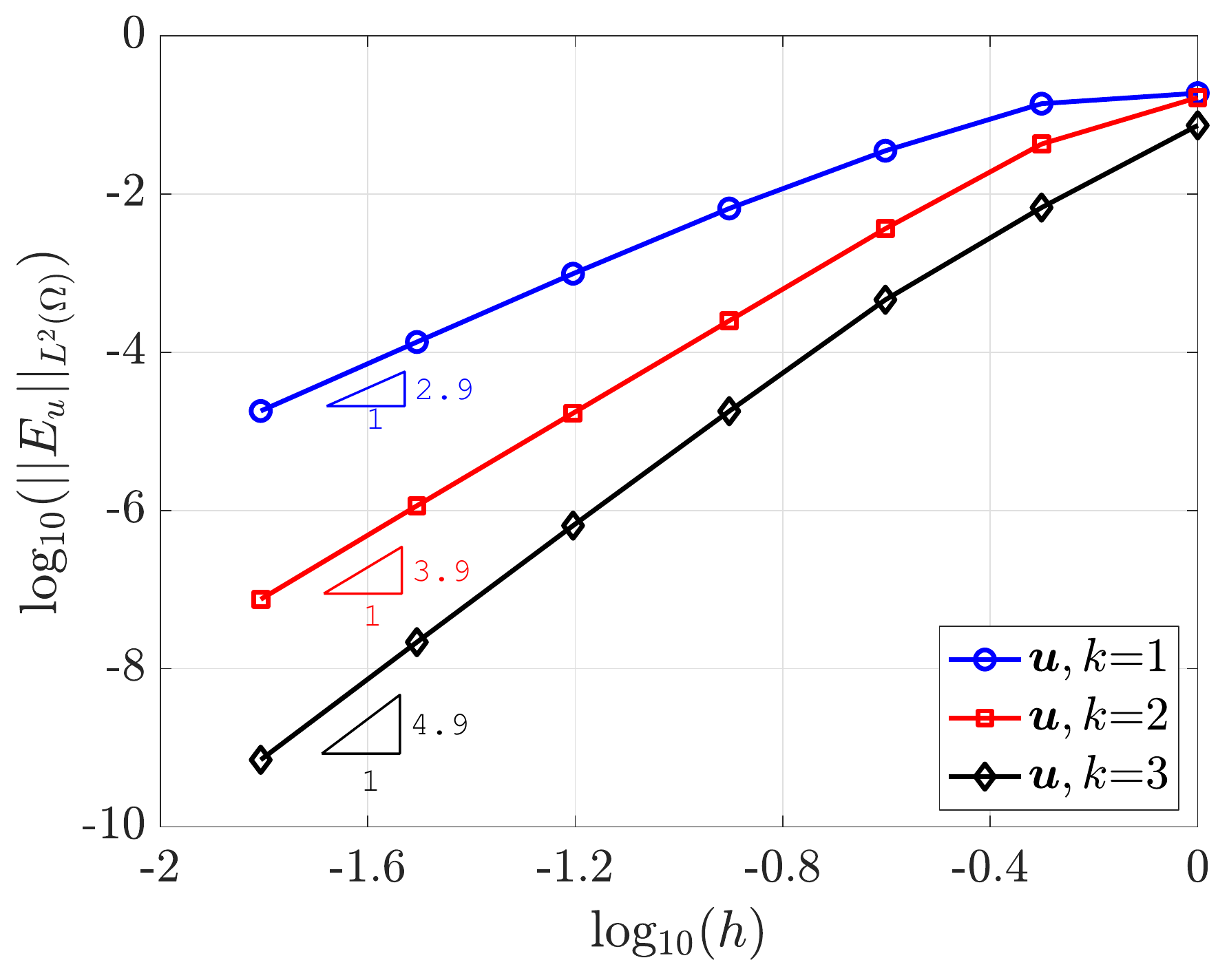}}
\caption{$h$-convergence of the error of the displacement field in the $\eltwo$ norm computed in the whole domain $\Omega$ using (a) CG, (b) HDG, (c) hybrid $\Kdeg{CG}$-$\Kdeg{HDG}$ coupling with polynomial approximation of degree $k$ in both subdomains and (d) hybrid $\KPdeg{CG}$-$\Kdeg{HDG}$ coupling with polynomial approximations of degree $k{+}1$ for CG and $k$ for HDG, with local postprocess of the primal variable.}
\label{fig:elasticityConvergenceU}
\end{figure}

\subsubsection{Locking-free approximation of the CG-HDG coupling}
\label{sc:2DelasticLocking}

Four strategies are considered for the discretization of the previously introduced elastic problem:
\begin{enumerate}
\item CG discretization in the whole domain $\Omega$;
\item HDG discretization in the whole domain $\Omega$;
\item hybrid CG-HDG discretization with polynomial approximation of degree $k$ both in $\CG{\Omega}$ and $\HDG{\Omega}$;
\item hybrid CG-HDG discretization with polynomial approximation of degree $k{+}1$ in $\CG{\Omega}$ and $k$ with local HDG postprocess in $\HDG{\Omega}$.
\end{enumerate}
It is worth noticing that the presence of nearly incompressible regions makes a CG discretization unfeasible in the whole domain, whereas employing HDG everywhere leads to a higher computational cost than CG in the subdomain $\CG{\Omega}$. Coupling the two approaches allows to devise a robust and flexible numerical scheme in $\Omega$, especially benefitting from nonuniform polynomial degrees of approximation.

For the following simulations, the HDG stabilization parameter $\tau$ is considered constant and equal to $2.5 {\times} 10^2$. It is worth recalling that a stabilization parameter $\tau {=} C E/\ell$, where $E$ is the Young's modulus, spanning from $25$ to $250$ for the problem under analysis, $\ell$ the characteristic length of the problem and $C$ a positive constant scaling factor guarantees stability and optimal convergence of the HDG discretization for the linear elasticity equation~\cite{RS-SGH:19}.
The Nitsche's parameter is set to $\gamma {=} 2.5 {\times} 10^3$ for the coupling based on uniform degree of approximation, whereas $\gamma {=} 2.5 {\times} 10^4$ is considered for the case of nonuniform approximations.

The convergence of the error of the displacement field measured in the $\eltwo$ norm as a function of the characteristic element size $h$ is presented in Figure~\ref{fig:elasticityConvergenceU}, for polynomial degree of approximation $k {=} 1,\ldots,3$.
%
%

Due to the presence of a nearly incompressible material, the CG approximation suffers from classical locking phenomena~\cite{Brenner-BS-92} which prevent convergence for $k {=} 1$, whereas suboptimal rates of order $k$ are obtained for higher-degree of polynomial approximation (Fig. \ref{fig:elasticityConvergenceU}a).
On the contrary, both HDG (Fig. \ref{fig:elasticityConvergenceU}b) and the hybrid $\Kdeg{CG}$-$\Kdeg{HDG}$ coupling (Fig. \ref{fig:elasticityConvergenceU}c) present optimal convergence of order $k {+} 1$ for the displacement field without locking effects~\cite{soon2009hybridizable} when polynomial of degree $k$ are utilized in both $\CG{\Omega}$ and $\HDG{\Omega}$.
Moreover, considering a nonuniform degree of approximation with polynomial functions of degree $k{+}1$ in $\CG{\Omega}$ and $k$ in $\HDG{\Omega}$ and exploiting the HDG postprocess strategy in Equation~\eqref{eq:ElasticityPostProcess}, a displacement field superconverging with order $k {+} 2$ in the whole domain $\Omega$ is obtained (Fig. \ref{fig:elasticityConvergenceU}d).

\begin{figure}
\centering
\subfigure[CG]{\includegraphics[width=0.45\textwidth]{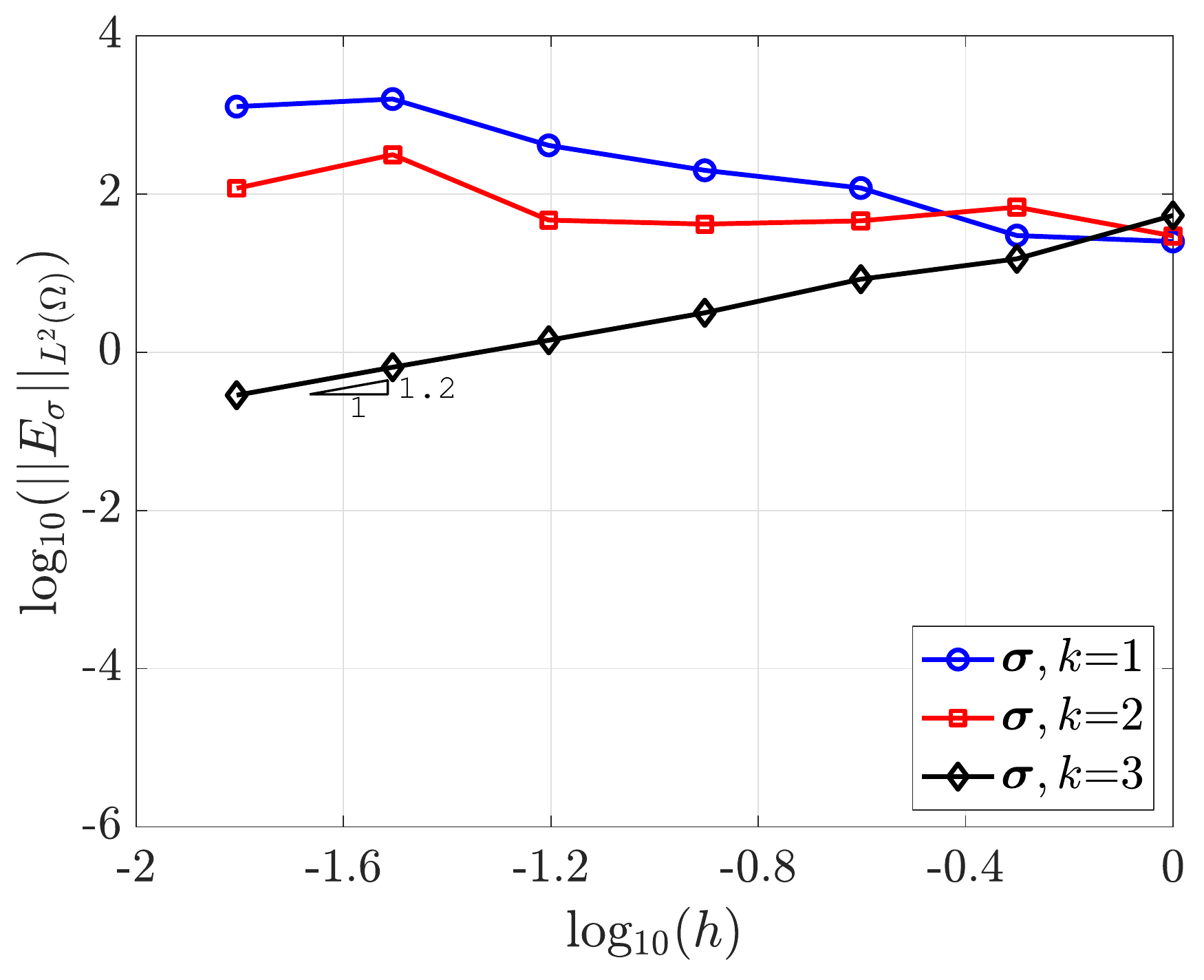}}
\hfill
\subfigure[HDG]{\includegraphics[width=0.45\textwidth]{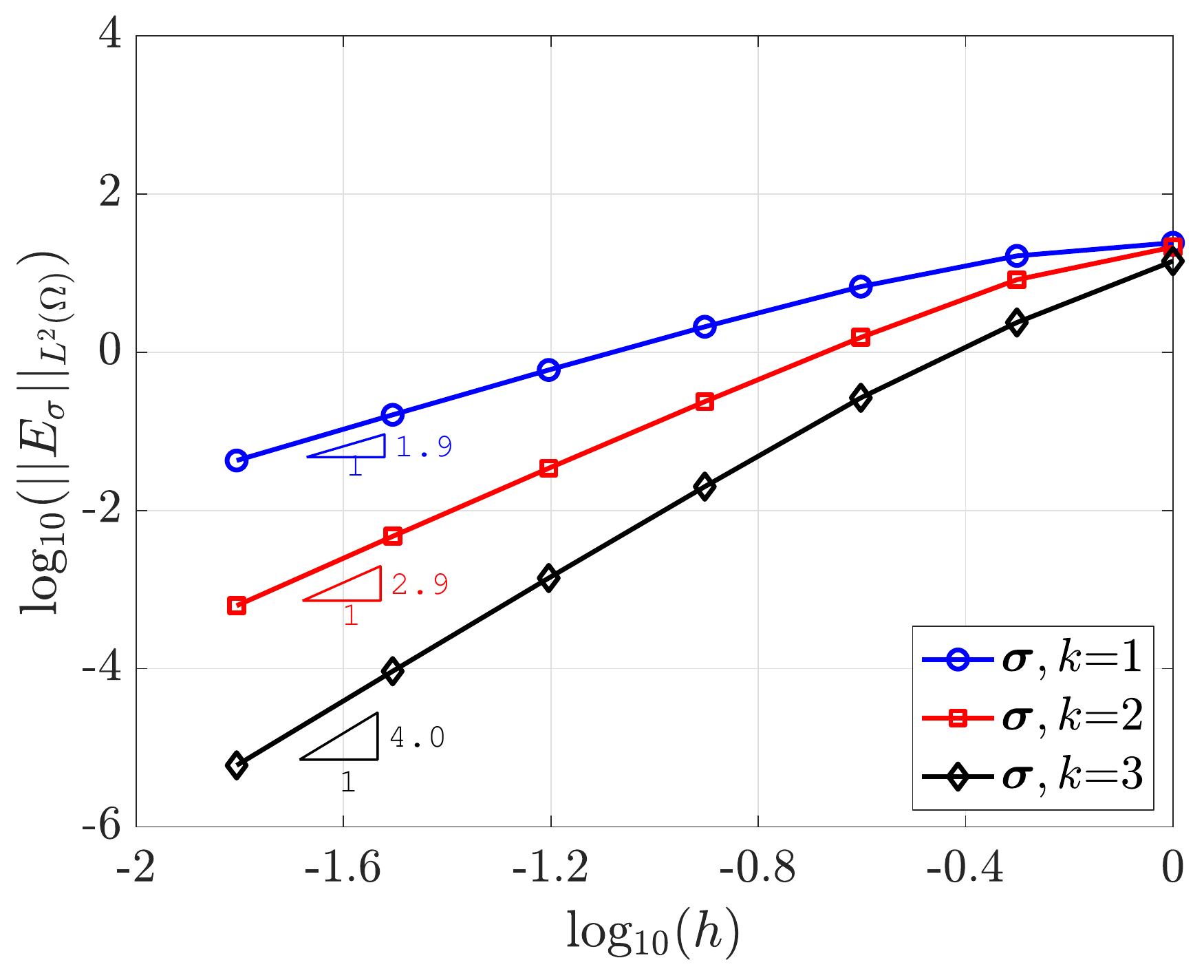}}

\subfigure[$\Kdeg{CG}$-$\Kdeg{HDG}$]{\includegraphics[width=0.45\textwidth]{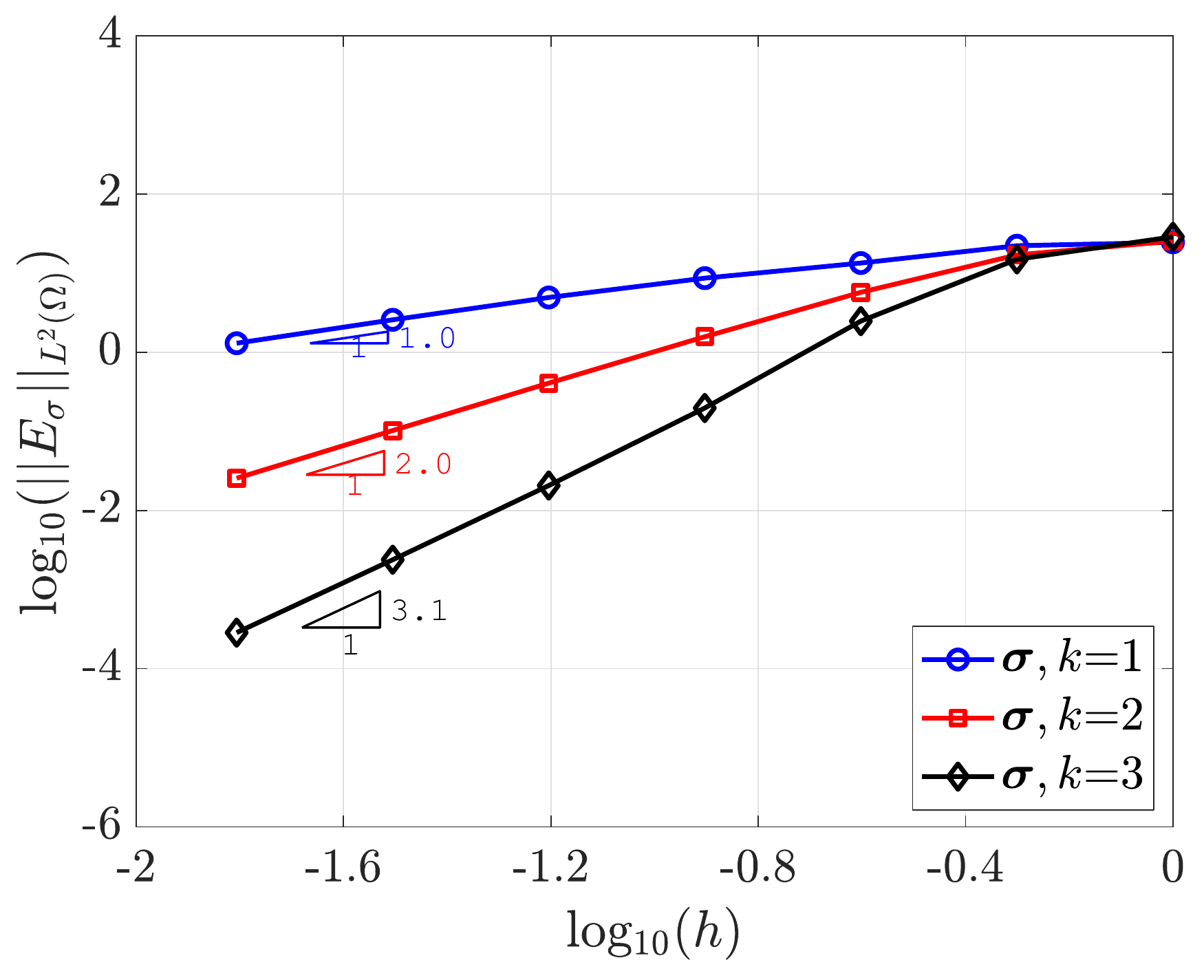}}
\hfill
\subfigure[$\KPdeg{CG}$-$\Kdeg{HDG}$]{\includegraphics[width=0.45\textwidth]{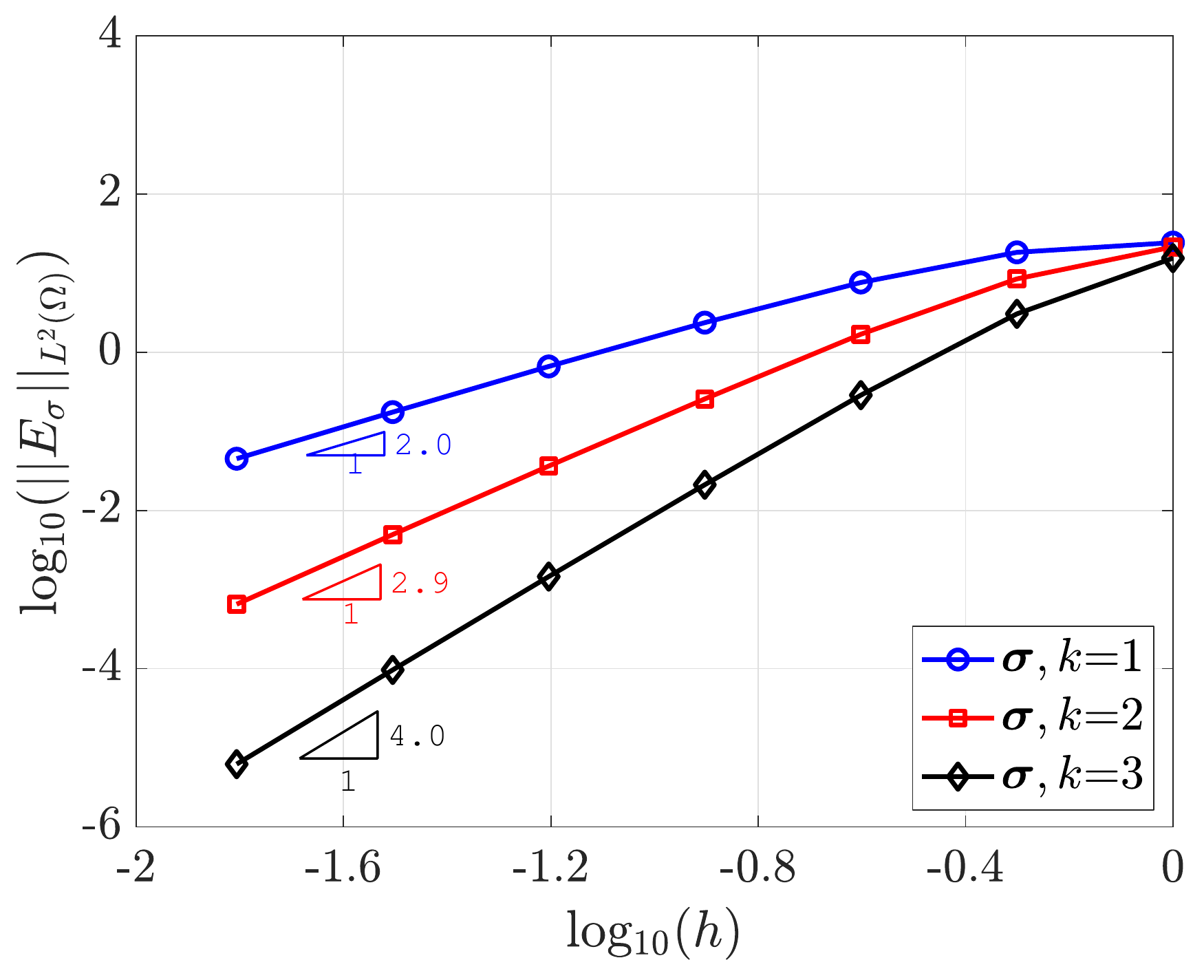}}
\caption{$h$-convergence of the error of the stress field in the $\eltwo$ norm computed in the whole domain $\Omega$ using (a) CG, (b) HDG, c) hybrid $\Kdeg{CG}$-$\Kdeg{HDG}$ coupling with polynomial approximation of degree $k$ in both subdomains and (d) hybrid $\KPdeg{CG}$-$\Kdeg{HDG}$ coupling with polynomial approximations of degree $k{+}1$ for CG and $k$ for HDG.}
\label{fig:elasticityConvergenceS}
\end{figure}

Figure~\ref{fig:elasticityConvergenceS} displays the convergence of the error of the stress field measured in the $\eltwo$ norm as a function of the characteristic element size $h$ for polynomial degree of approximation $k {=} 1,\ldots,3$.

The CG approximation of the stress field being computed as a postprocess of the displacement field leads to unreliable results due to the locking effects (Fig. \ref{fig:elasticityConvergenceS}a).
On the contrary, mixed formulations provide a direct approximation of the stress tensor~\cite{brezzi1991mixed}.
More precisely, HDG exploits the definition of the pointwise symmetric mixed variable $\HDG{\bL}$ in~\cite{RS-SGKH:18,MG-GKSH:18} to obtain optimal convergence of the stress with order $k {+} 1$ (Fig.~\ref{fig:elasticityConvergenceS}b).

To reduce the computational cost of the approximation of the linear elastic problem in the whole domain $\Omega$, a CG discretization is considered in the compressible subdomain $\CG{\Omega}$ and an HDG one in the nearly incompressible region $\HDG{\Omega}$.
The approximation of the stress computed using the $\Kdeg{CG}$-$\Kdeg{HDG}$ coupling in Figure~\ref{fig:elasticityConvergenceS}c presents suboptimal convergence of order $k$ in $\Omega$ since $\stress$ is computed as a postprocess of the primal solution $\CG{\bu}$ in the subdomain $\CG{\Omega}$.
Optimal convergence of the stress field in the whole domain $\Omega$ is recovered by constructing the hybrid $\KPdeg{CG}$-$\Kdeg{HDG}$ coupling based on a nonuniform degree of approximation with polynomial functions of degree $k{+}1$ in $\CG{\Omega}$ and $k$ in $\HDG{\Omega}$ (Fig.~\ref{fig:elasticityConvergenceS}d).

\section{Application to multimaterial elastostatics problems}
\label{sc:multiMaterial}

In this section, the proposed hybrid CG-HDG coupling is applied to linear elastostastics problems of engineering interest. More precisely, problems with domains featuring multiple materials and inhomogeneous mechanical properties are investigated.

\begin{figure}
\centering
\includegraphics[width=0.3\textwidth]{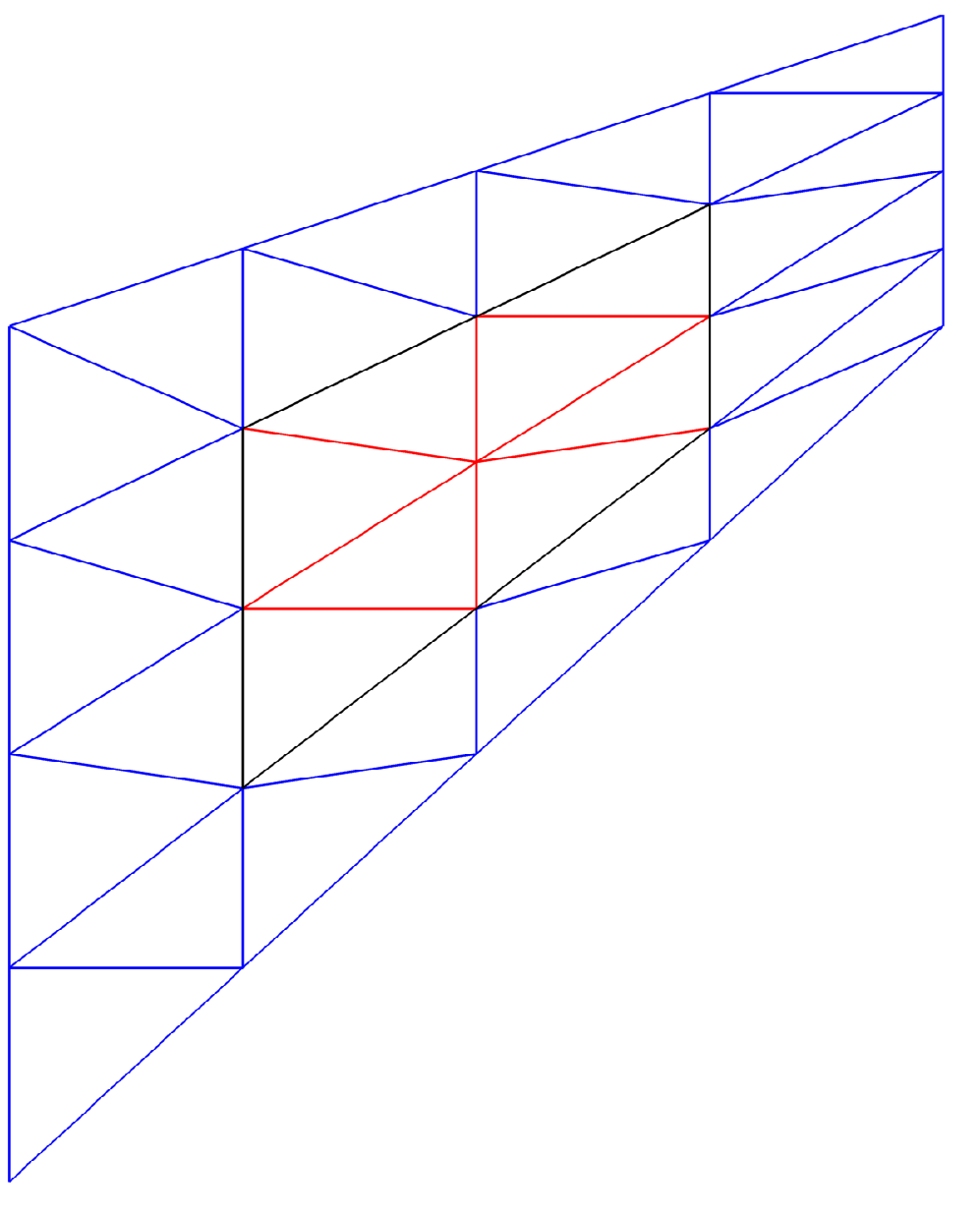}
\hfill
\includegraphics[width=0.3\textwidth]{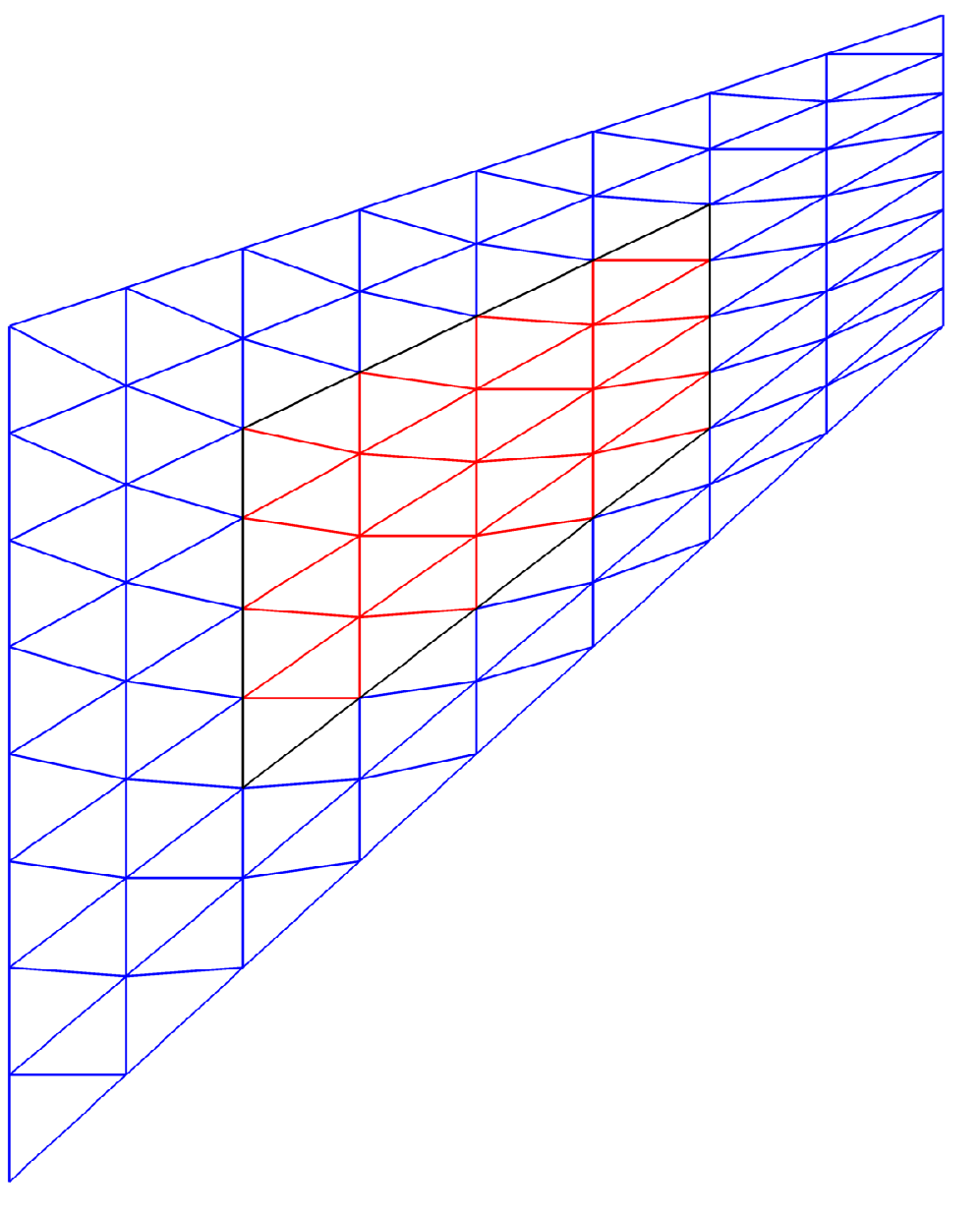}
\hfill
\includegraphics[width=0.3\textwidth]{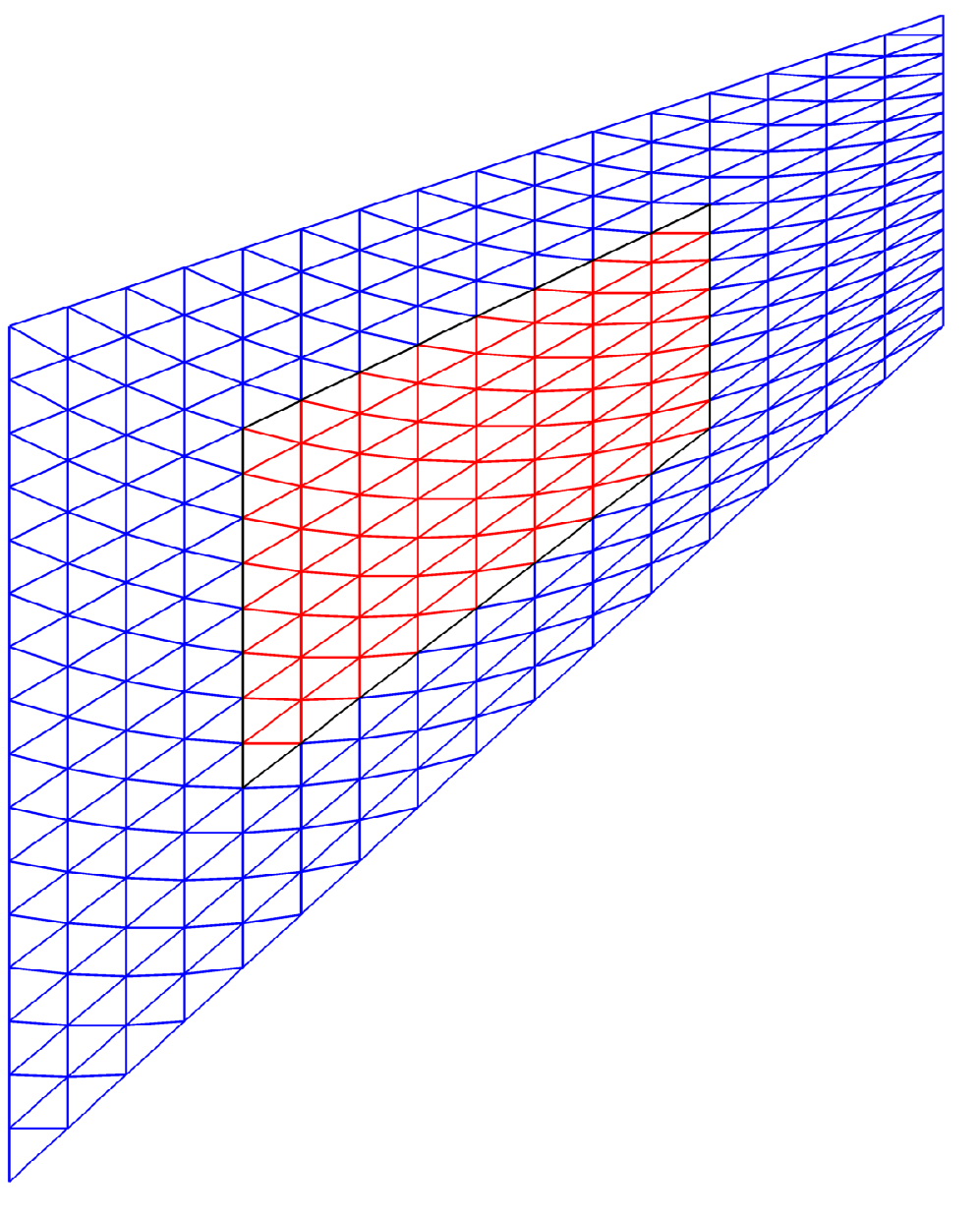}
\caption{First three levels of refinement of the mesh used for the Cook's membrane problem. Red: CG subdomain $\CG{\Omega}$. Blue: HDG subdomain $\HDG{\Omega}$. Black: interface $\Gamma_I$.}
\label{fig:cookMesh}
\end{figure}
\begin{figure}
	\centering
	\begin{tabular}[c]{@{}l@{}c@{}c@{}c@{}c@{}}
		$u_x$ & 
		\parbox[c]{0.24\textwidth}{\includegraphics[width=0.24\textwidth]{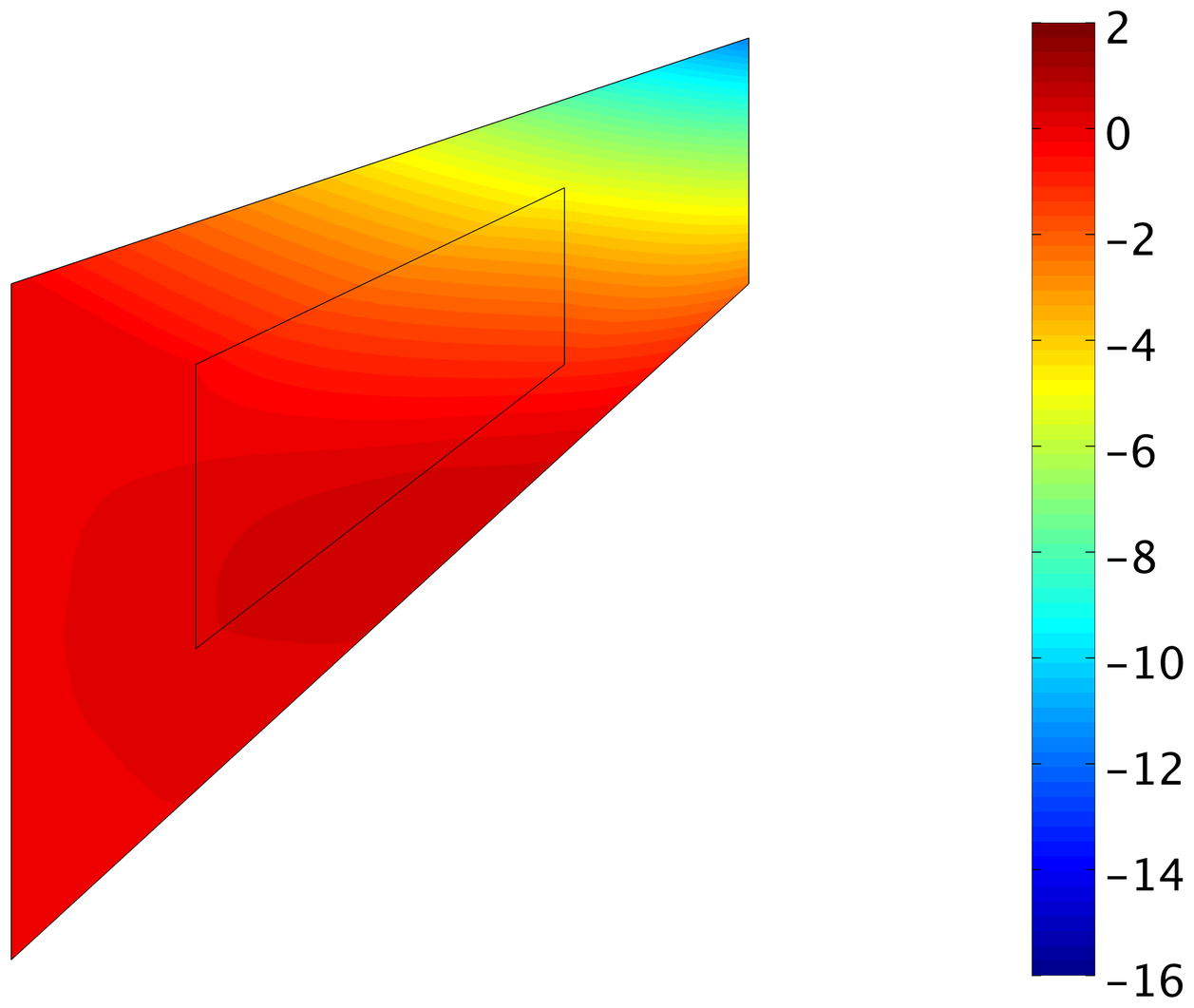}} & 
	 	\parbox[c]{0.24\textwidth}{\includegraphics[width=0.24\textwidth]{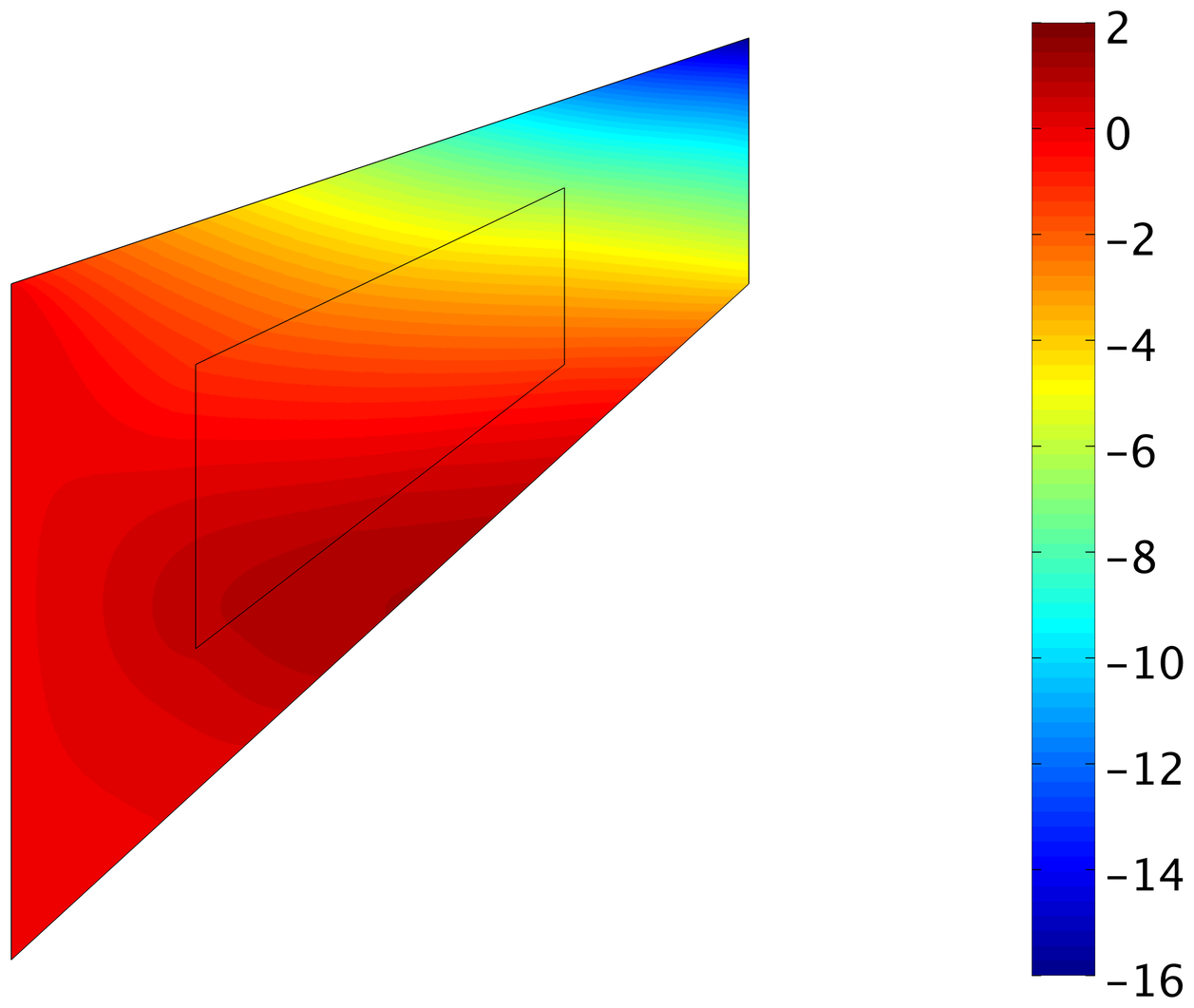}} &
		\parbox[c]{0.24\textwidth}{\includegraphics[width=0.24\textwidth]{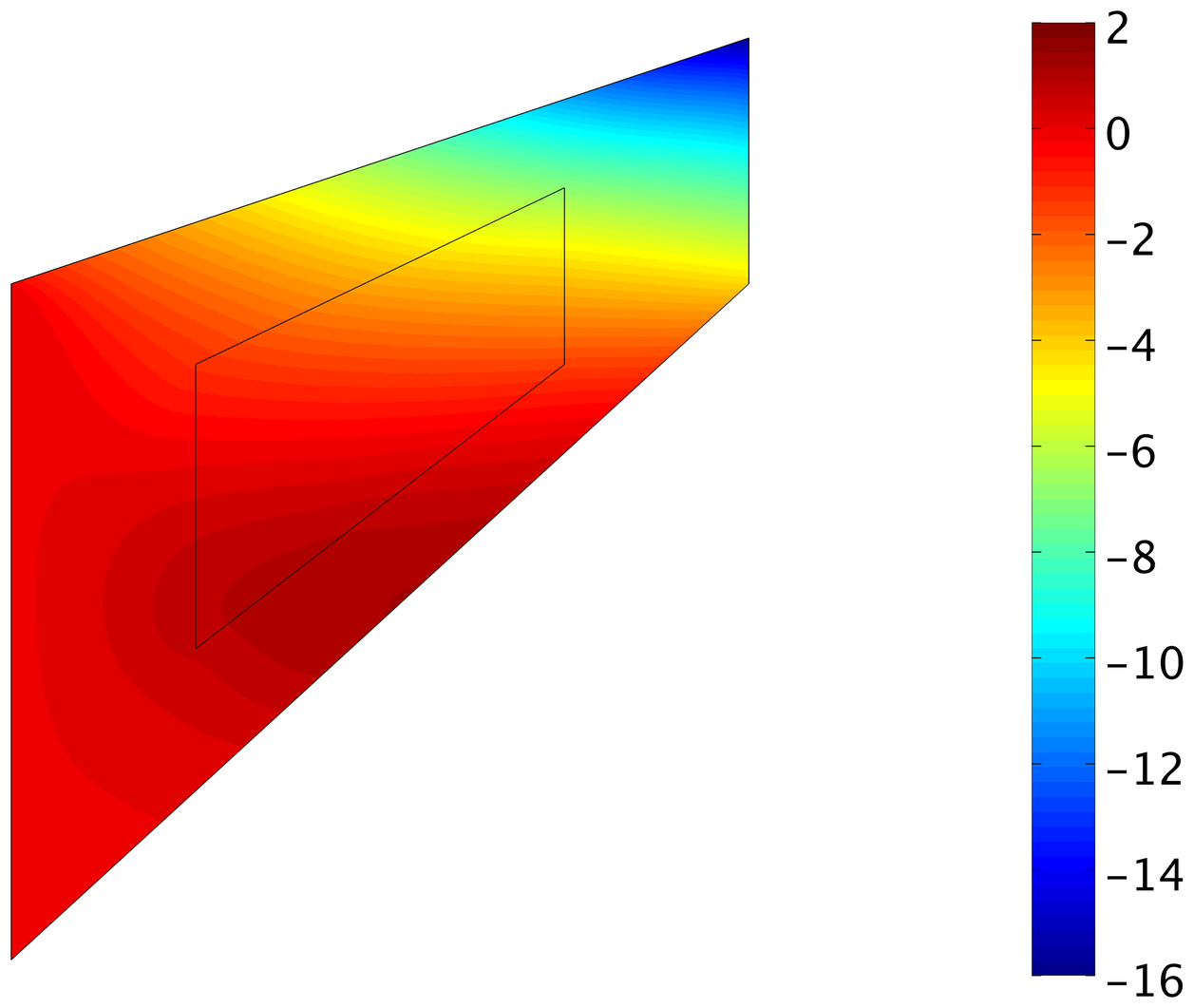}} &
		\parbox[c]{0.24\textwidth}{\includegraphics[width=0.24\textwidth]{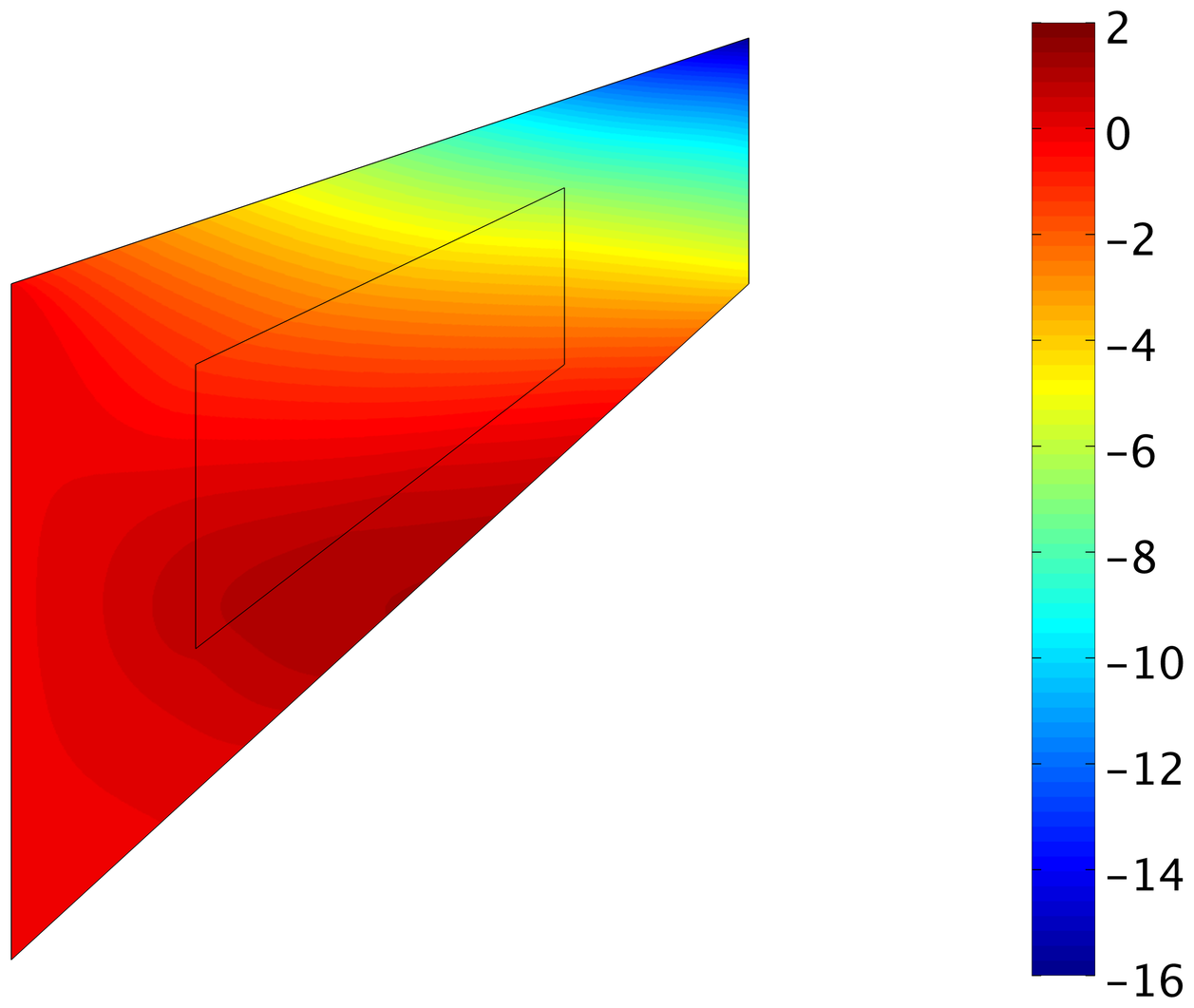}} \\[5em]
		$u_y$ & 
		\parbox[c]{0.24\textwidth}{\includegraphics[width=0.24\textwidth]{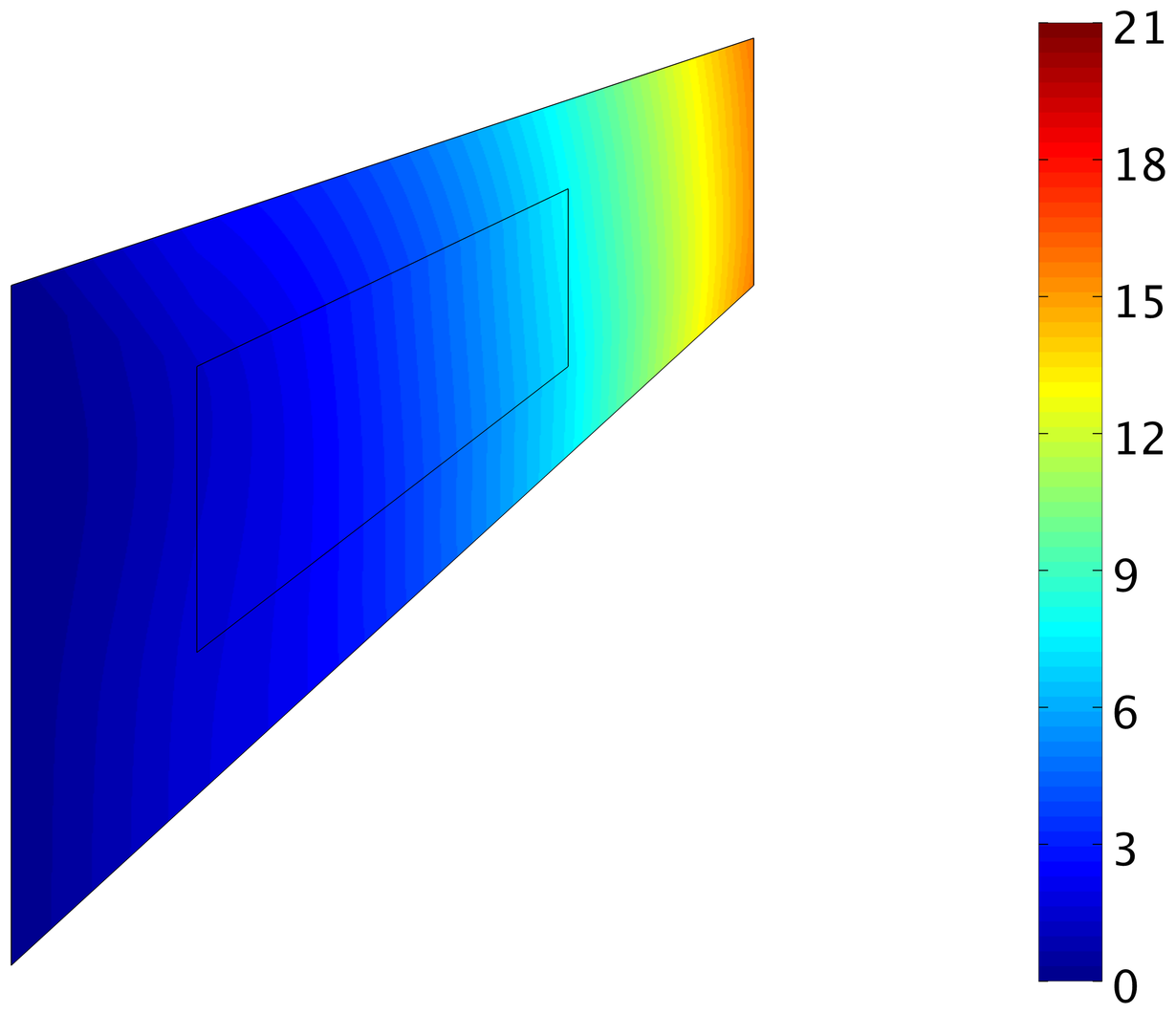}} & 
	 	\parbox[c]{0.24\textwidth}{\includegraphics[width=0.24\textwidth]{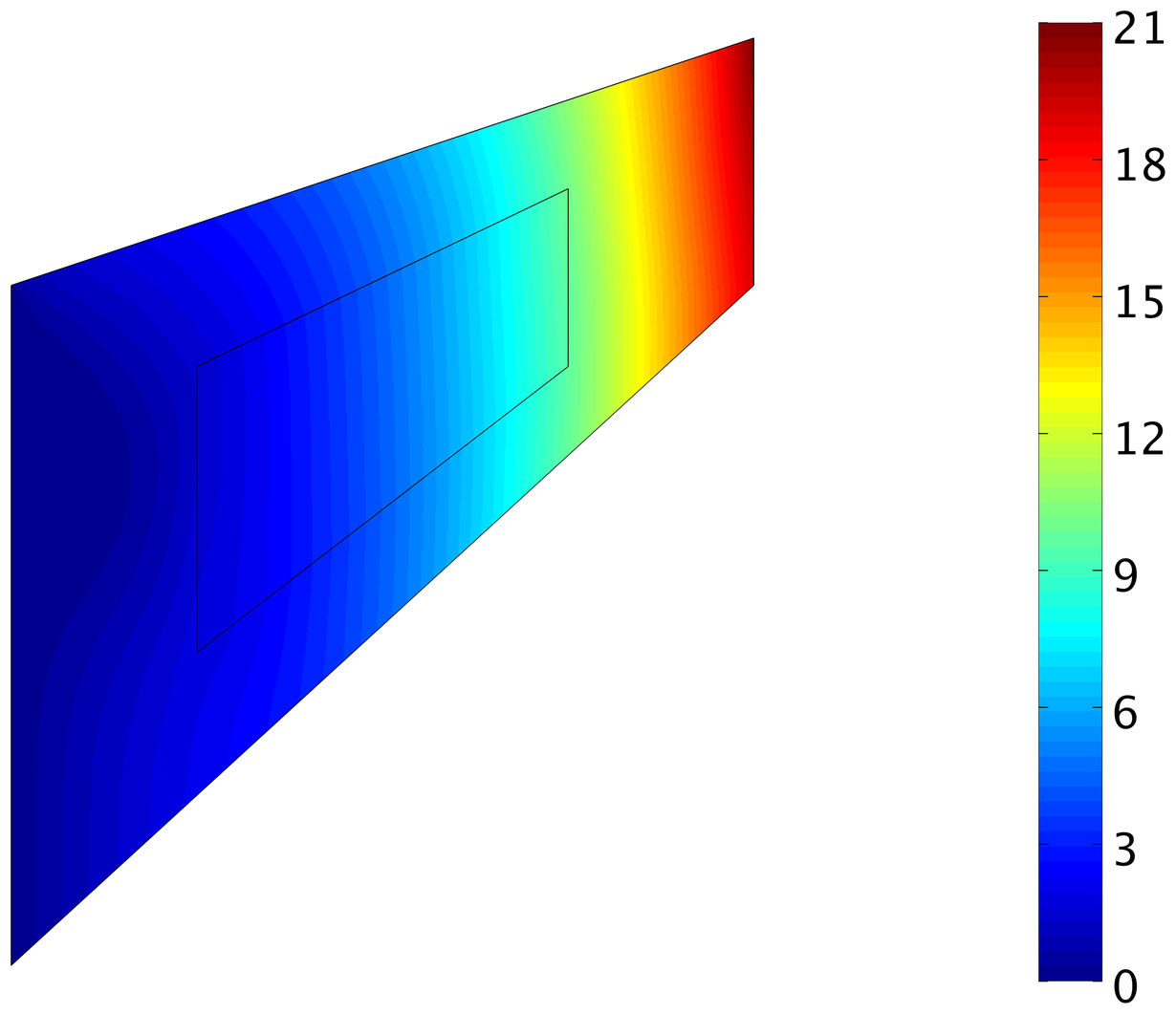}} &
		\parbox[c]{0.24\textwidth}{\includegraphics[width=0.24\textwidth]{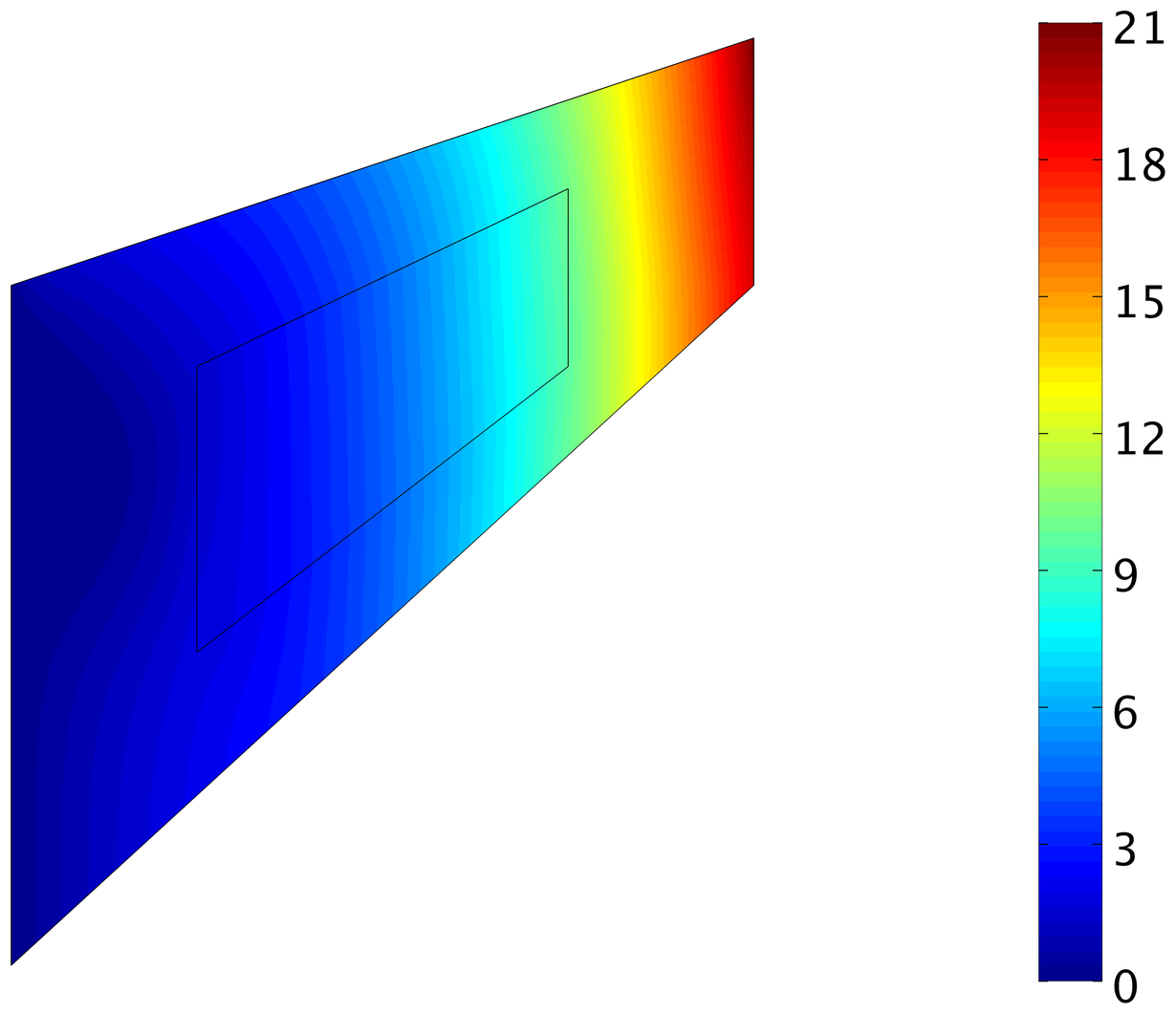}} &
		\parbox[c]{0.24\textwidth}{\includegraphics[width=0.24\textwidth]{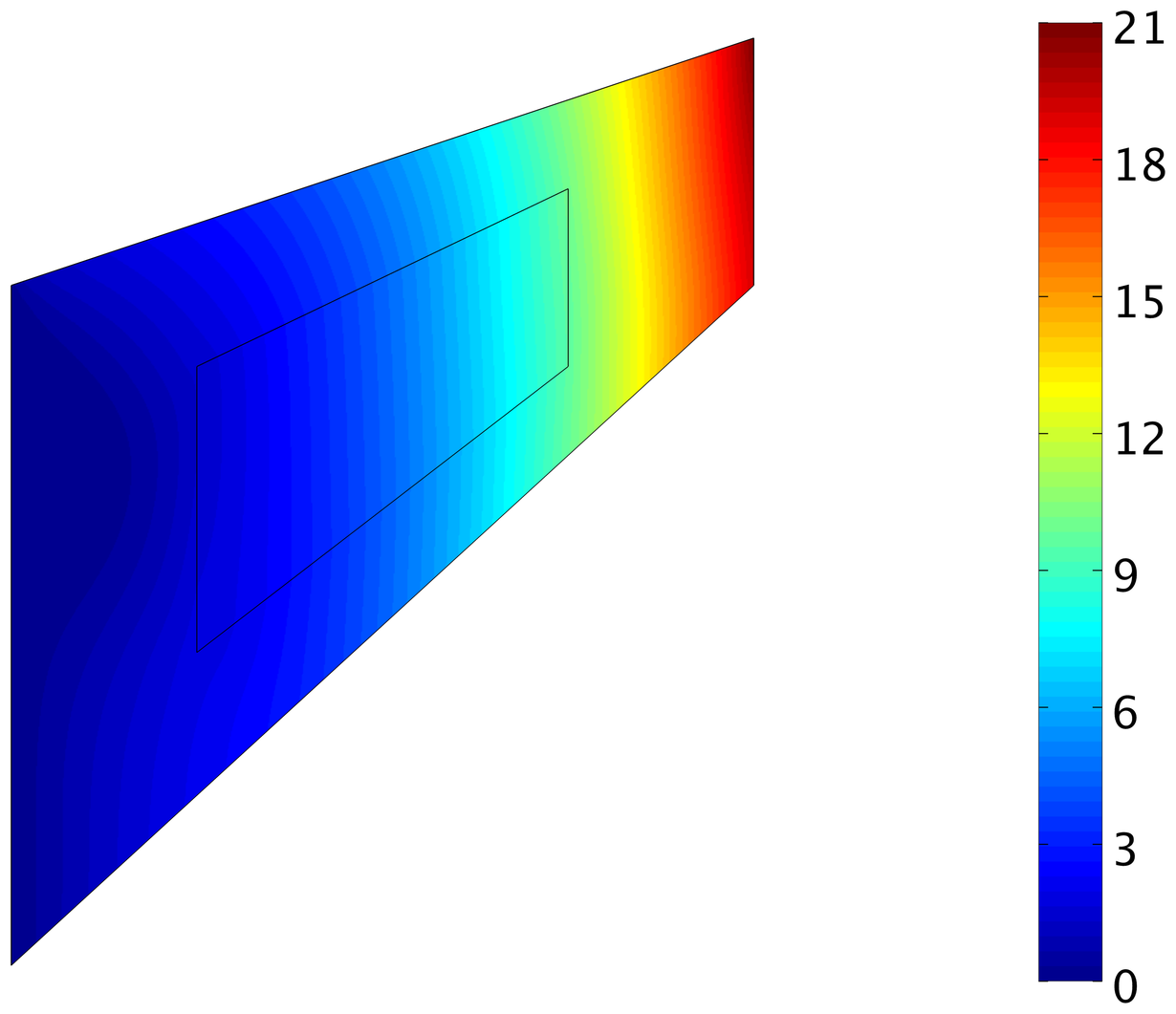}} \\[5em]
		$\sigma_{xx}$ & 
		\parbox[c]{0.24\textwidth}{\includegraphics[width=0.24\textwidth]{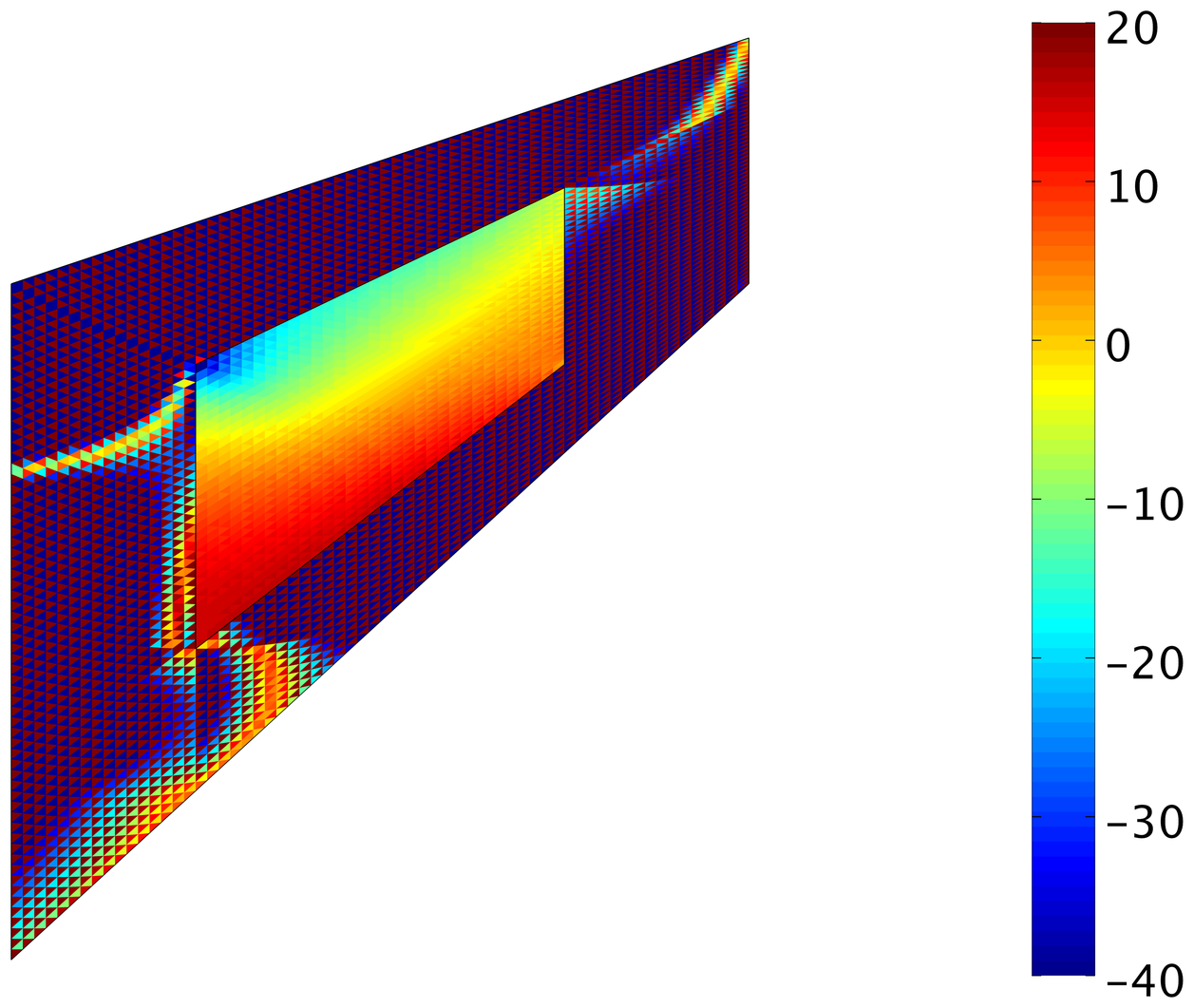}} & 
	 	\parbox[c]{0.24\textwidth}{\includegraphics[width=0.24\textwidth]{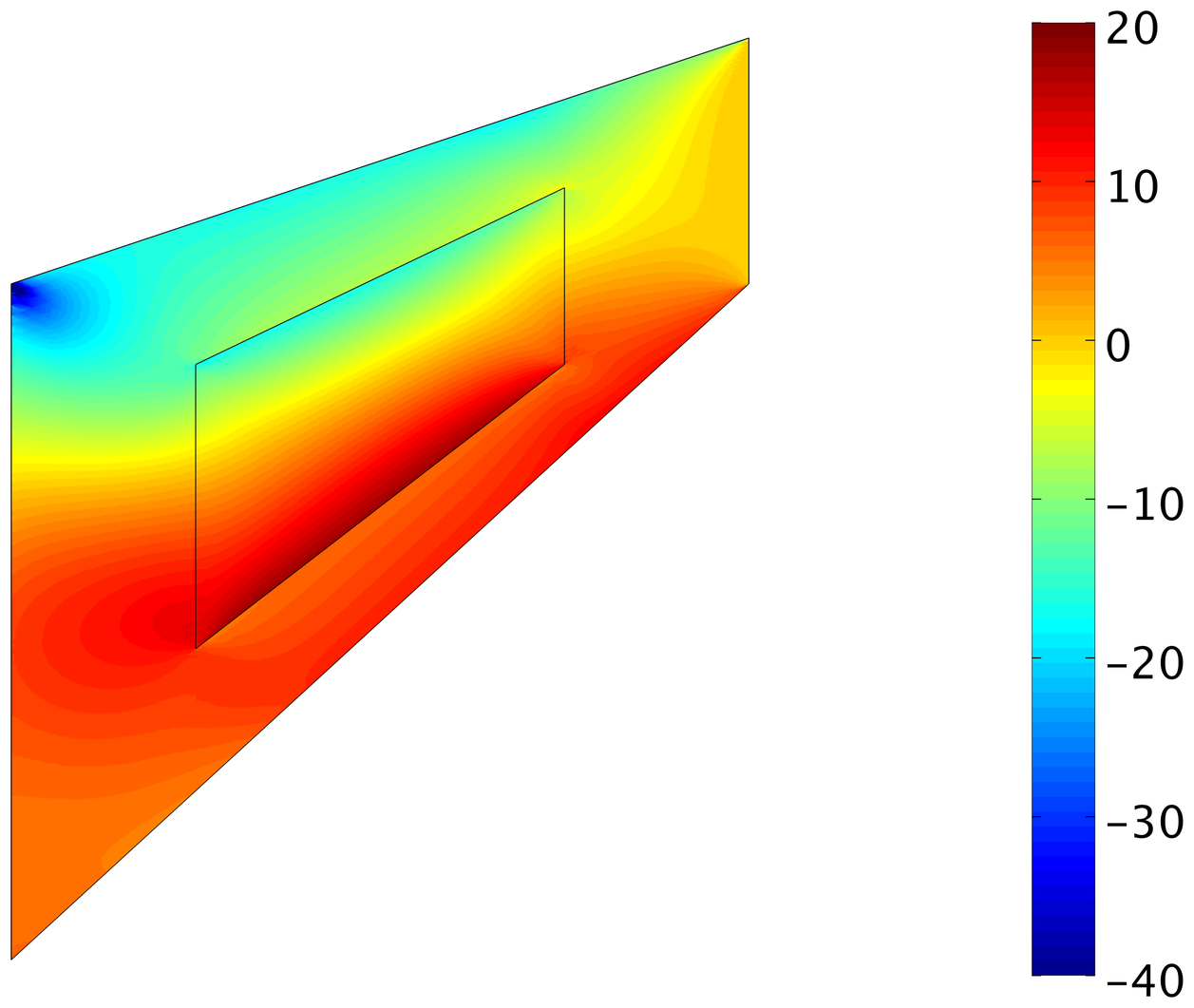}} &
		\parbox[c]{0.24\textwidth}{\includegraphics[width=0.24\textwidth]{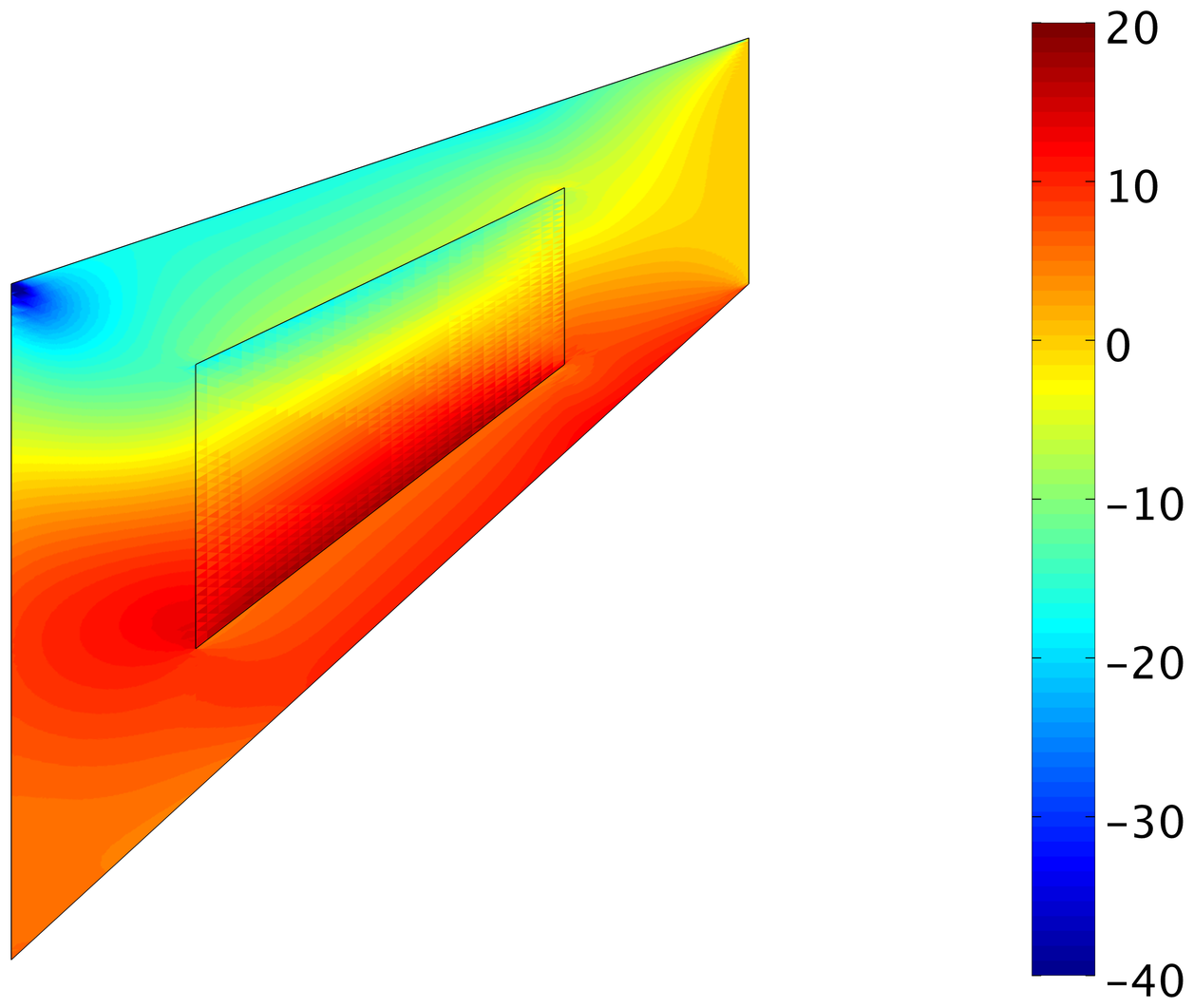}} &
		\parbox[c]{0.24\textwidth}{\includegraphics[width=0.24\textwidth]{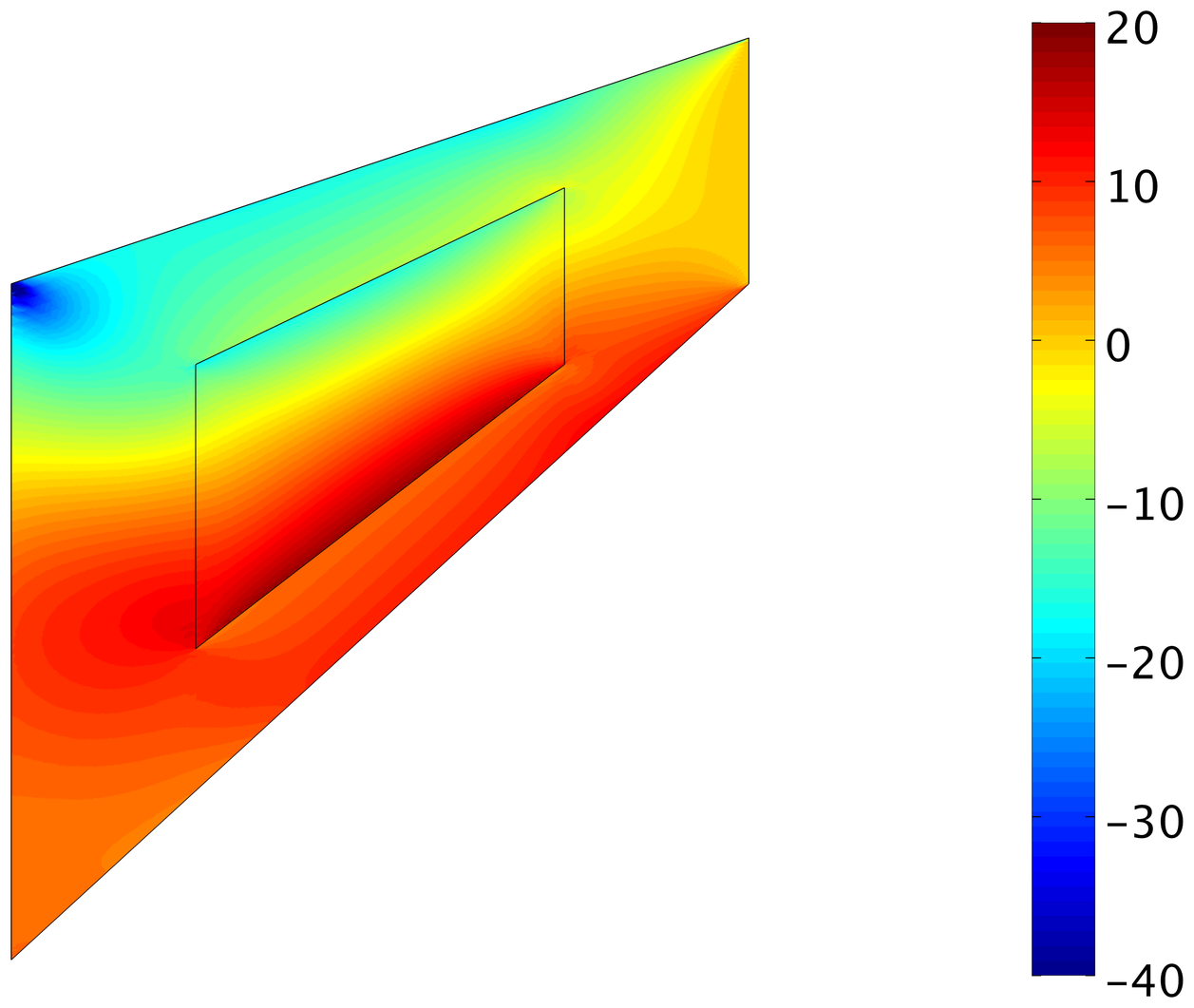}} \\[5em]
		$\sigma_{yy}$ & 
		\parbox[c]{0.24\textwidth}{\includegraphics[width=0.24\textwidth]{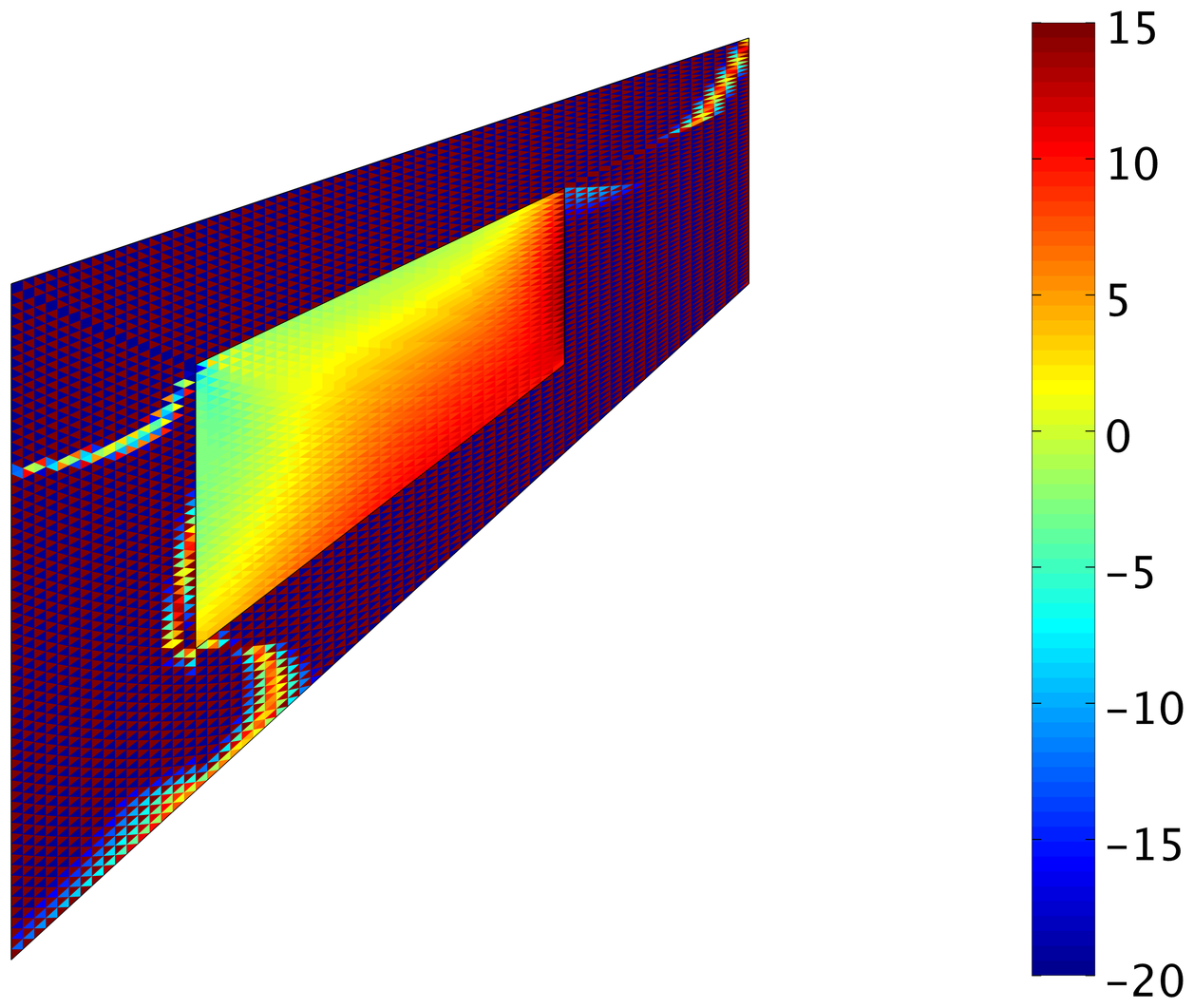}} & 
	 	\parbox[c]{0.24\textwidth}{\includegraphics[width=0.24\textwidth]{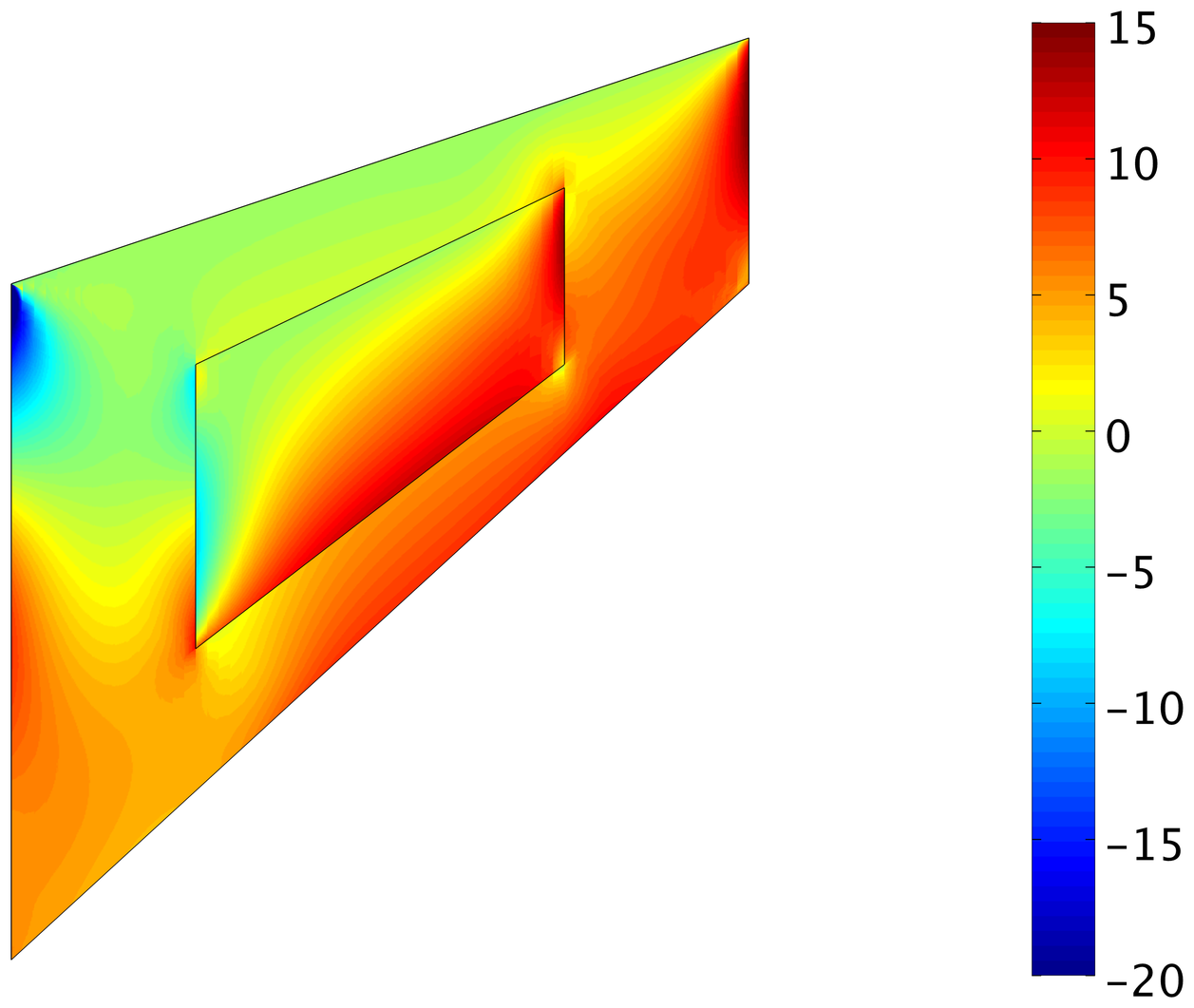}} &
		\parbox[c]{0.24\textwidth}{\includegraphics[width=0.24\textwidth]{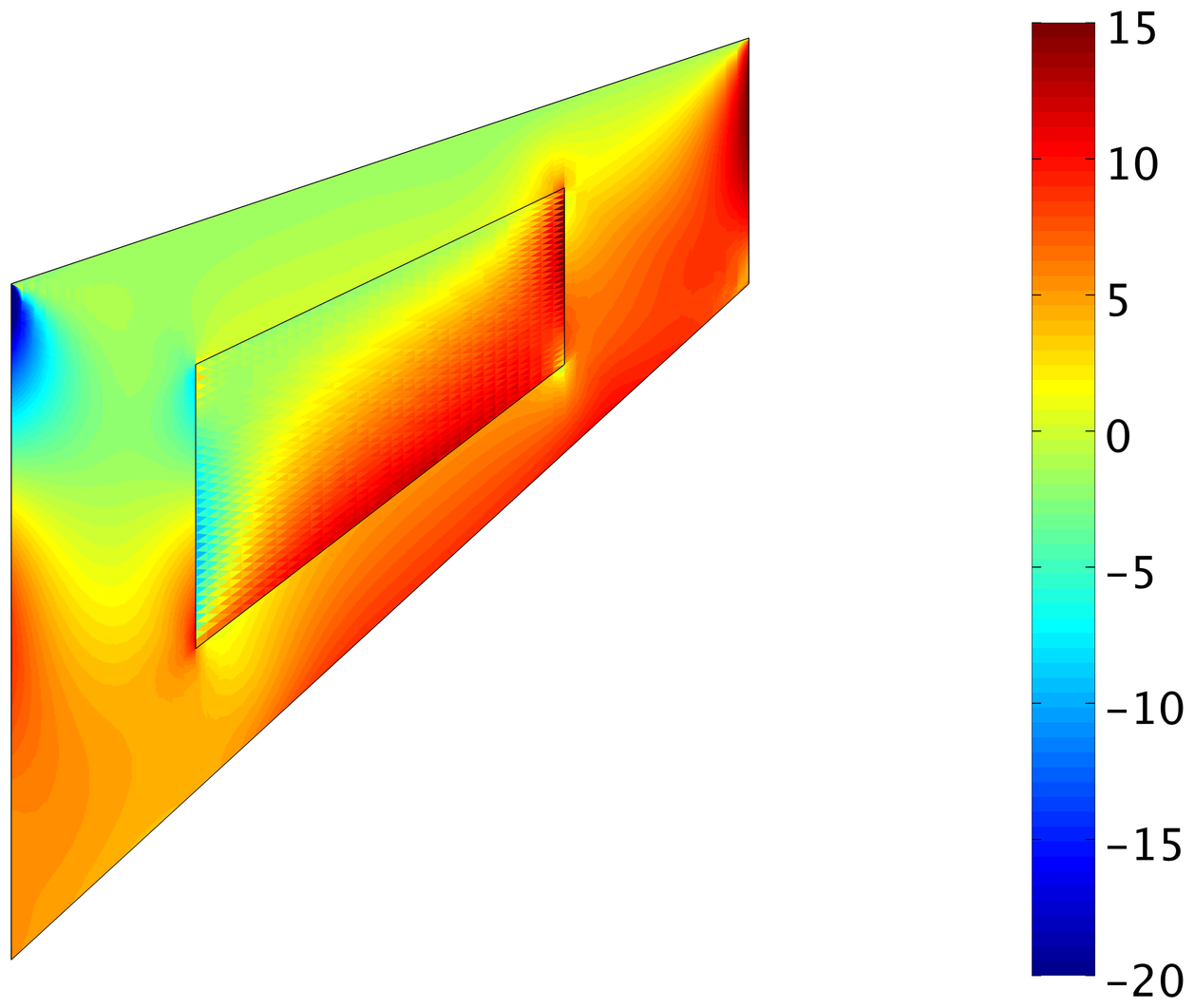}} &
		\parbox[c]{0.24\textwidth}{\includegraphics[width=0.24\textwidth]{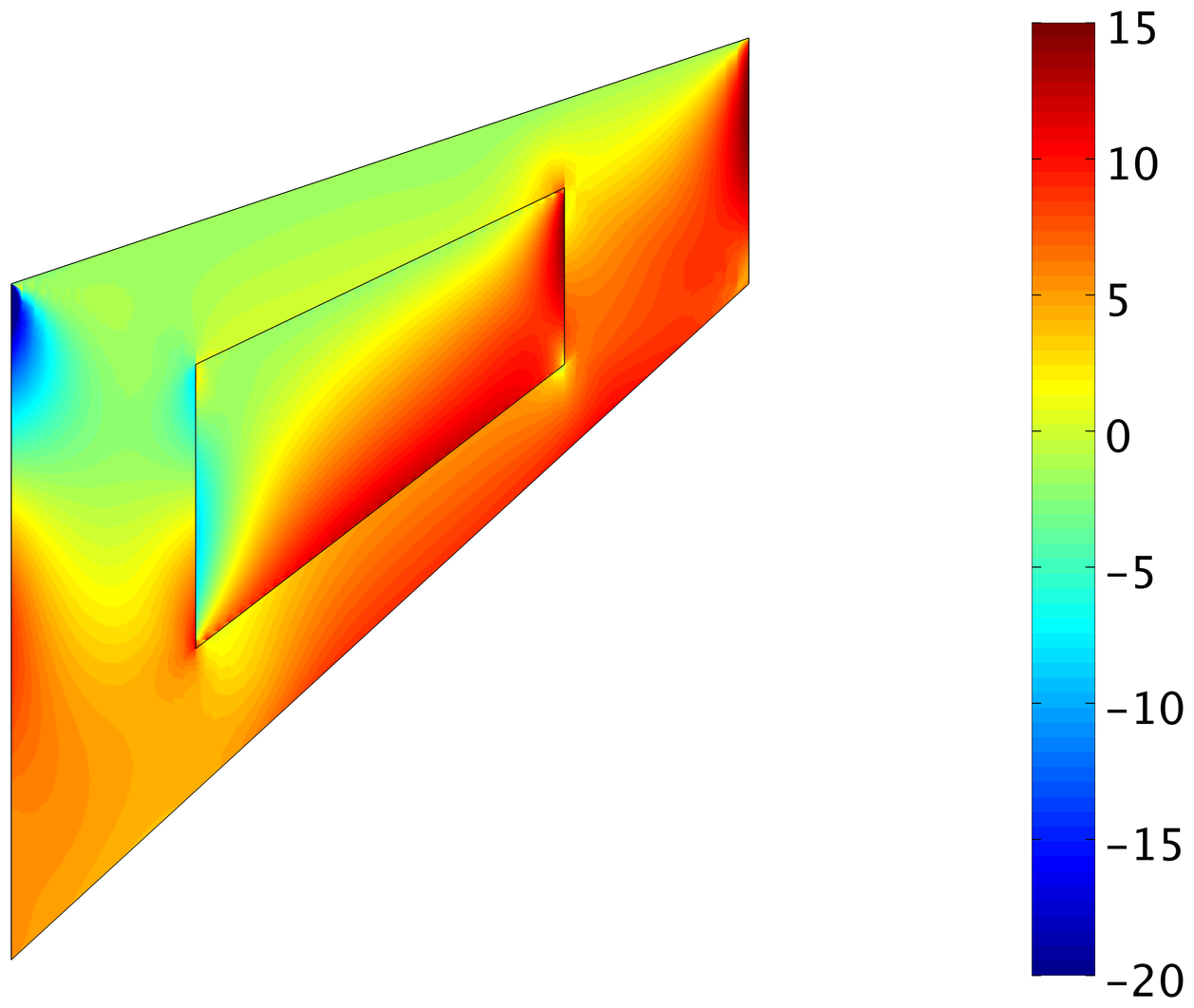}} \\[5em]
		$\sigma_{xy}$ &
		\parbox[c]{0.24\textwidth}{\includegraphics[width=0.24\textwidth]{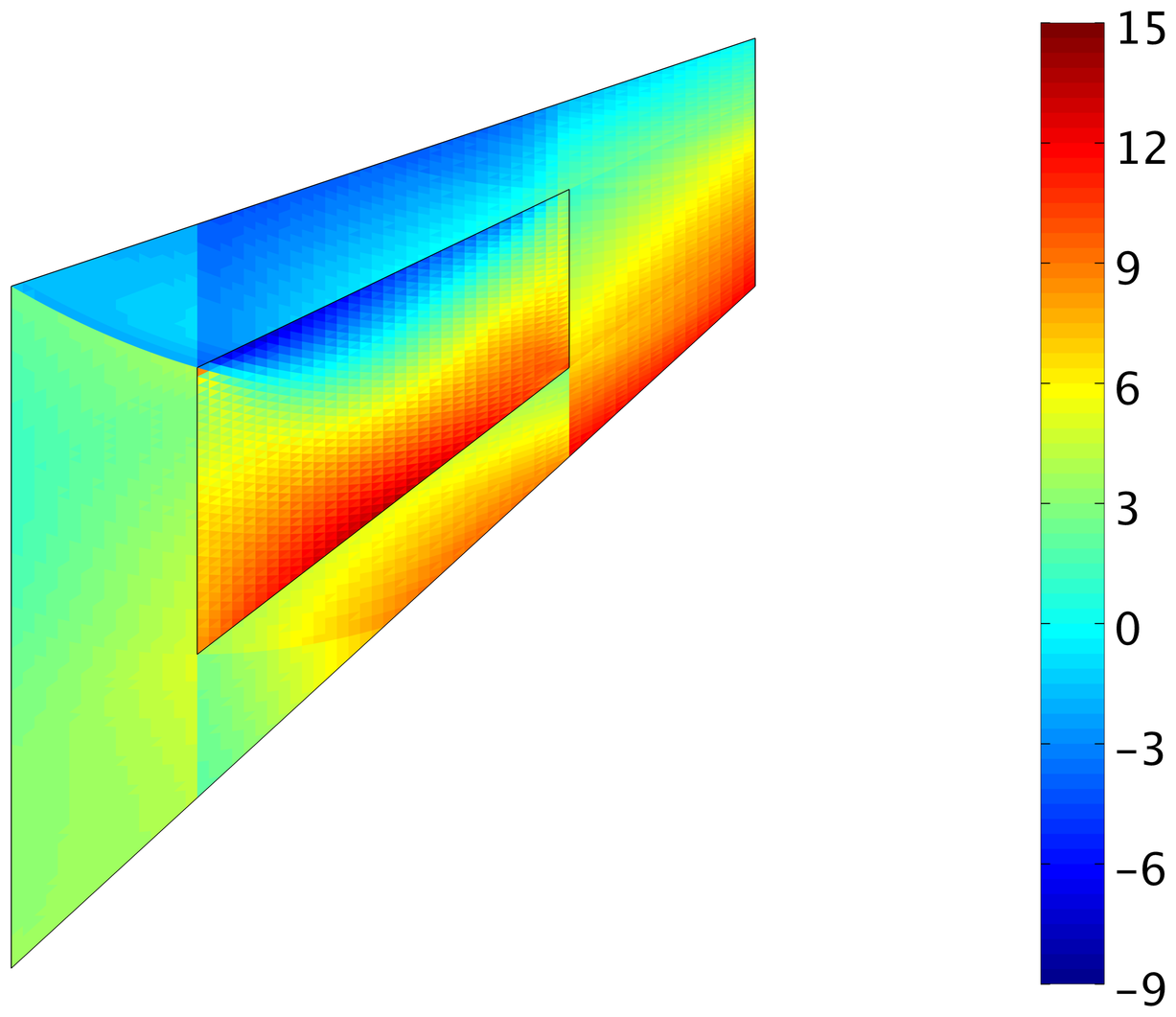}} &
		\parbox[c]{0.24\textwidth}{\includegraphics[width=0.24\textwidth]{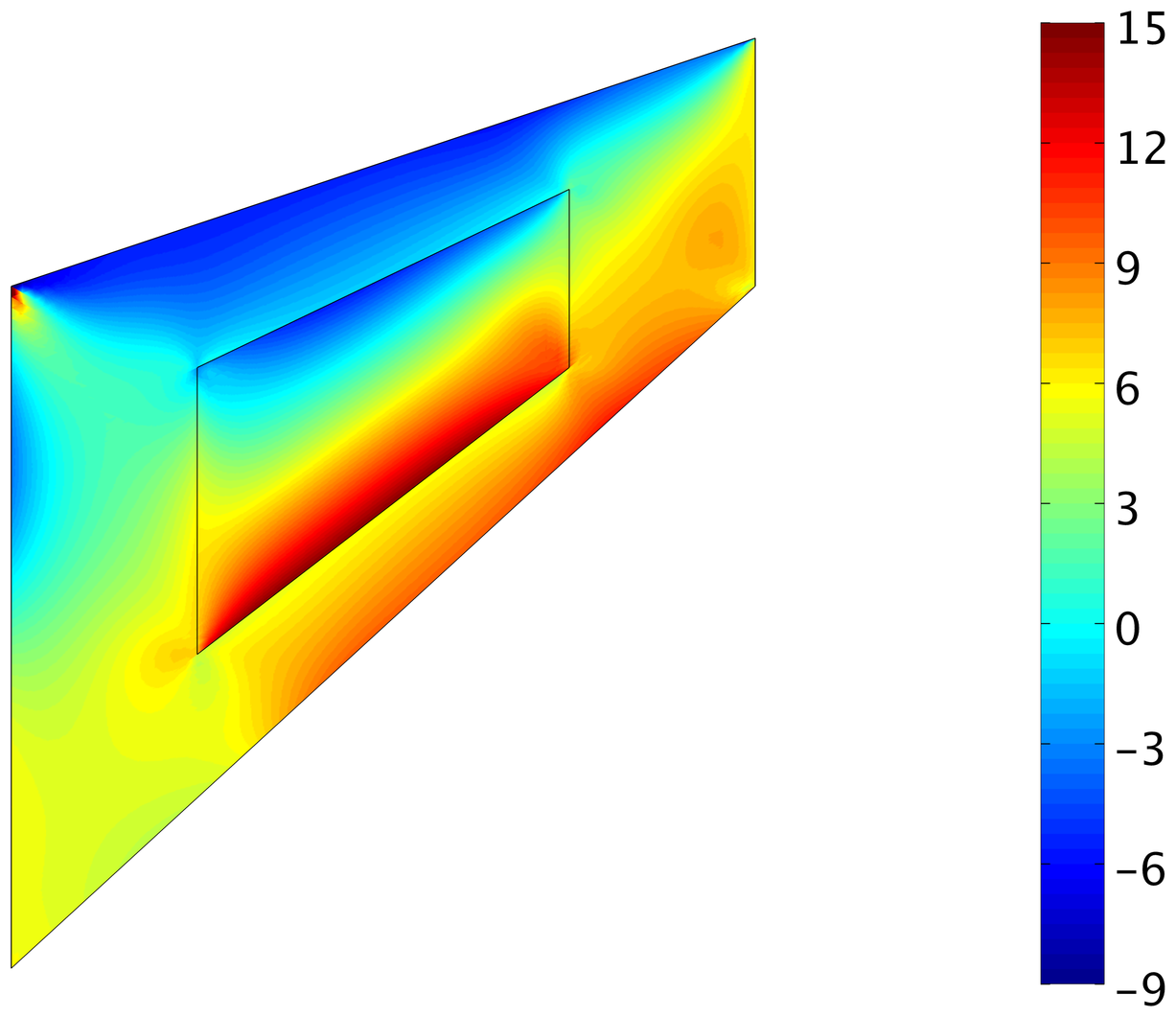}} &
		\parbox[c]{0.24\textwidth}{\includegraphics[width=0.24\textwidth]{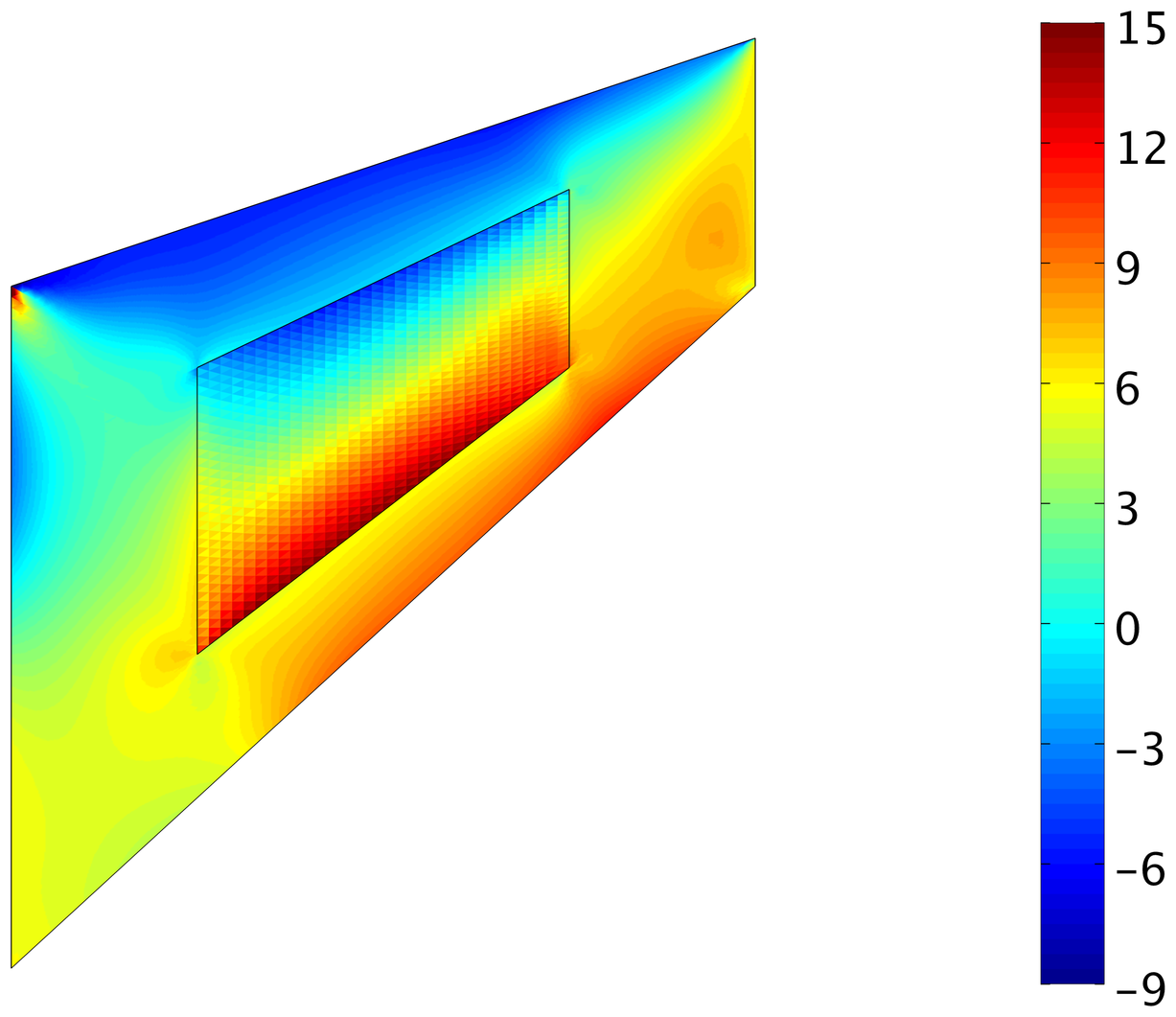}} &
		\parbox[c]{0.24\textwidth}{\includegraphics[width=0.24\textwidth]{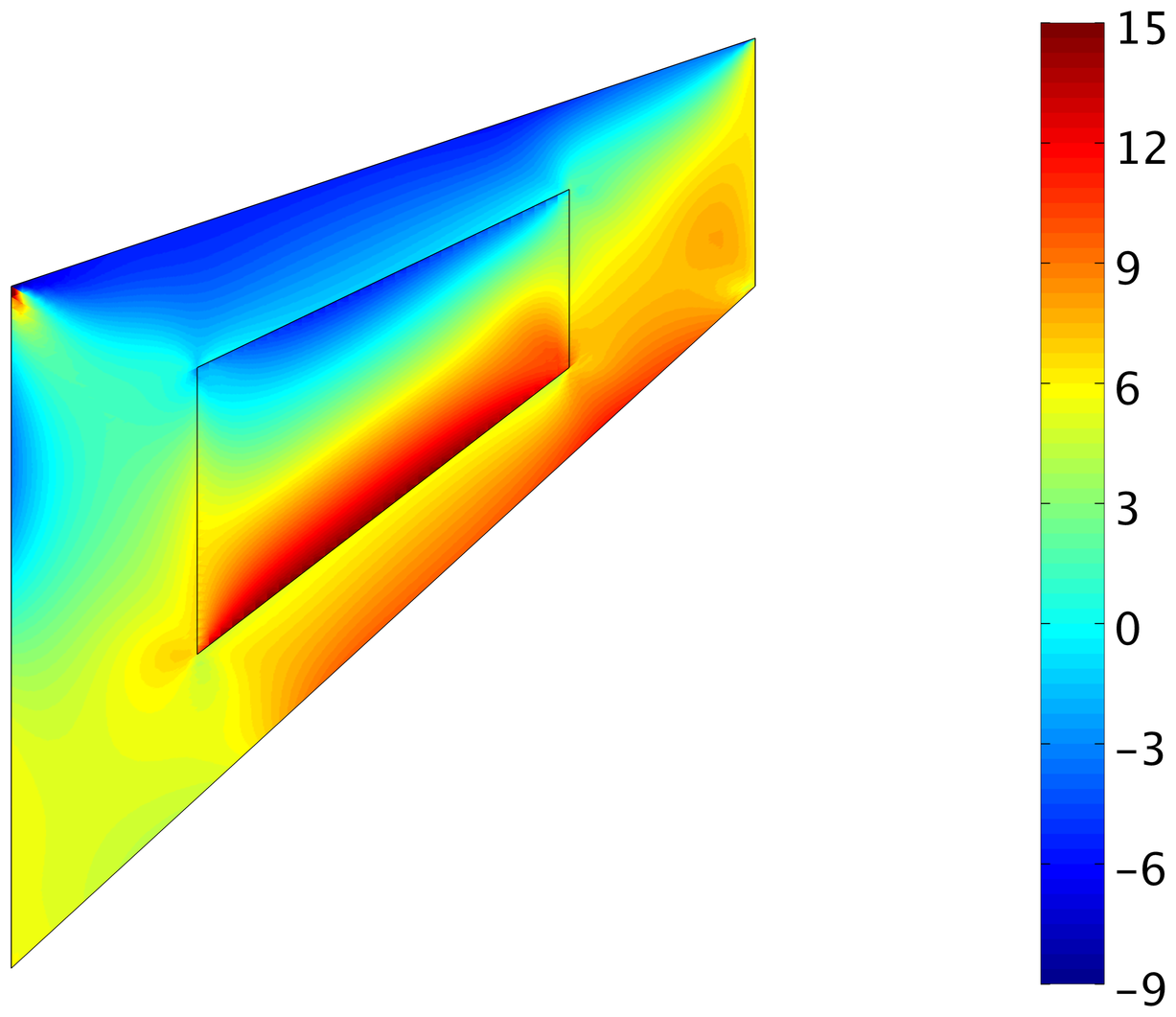}} \\[5em]
		  & CG & HDG & $\Kdeg{CG}$-$\Kdeg{HDG}$ & $\KPdeg{CG}$-$\Kdeg{HDG}$ \\
	\end{tabular}
\caption{Approximation of the displacement and the stress fields in the bimaterial Cook's membrane problem using the fifth level of mesh refinement for a polynomial approximation of degree $k {=} 1$. The interface $\Gamma_I$ is displayed in black, the compressible subdomain $\CG{\Omega}$ being on the inside and the nearly incompressible region $\HDG{\Omega}$ on the outside.}
\label{fig:cooksMembraneSolution}
\end{figure}

\subsection{Bimaterial Cook's membrane problem}
\label{sc:cooksMembrane}

The classical bending-dominated problem of Cook's membrane~\cite{cook2001concepts} is revisited to take into account two materials with different mechanical properties as proposed in~\cite{Lamichhane2009}.
The tapered plate $\Omega$ defined by the points $(0,0)$, $(48,44)$, $(48,60)$, $(0,44)$ is decomposed in two nonoverlapping subdomains such that $\CG{\Omega}$ is the convex region identified by the points $(12,20.25)$, $(36,38.75)$, $(36,50.25)$, $(12,38.75)$ and $\HDG{\Omega} {=} \Omega \setminus \CG{\Omega}$.
The first three levels of mesh refinement are displayed in Figure~\ref{fig:cookMesh}, where the subdomain $\CG{\Omega}$ is represented in red, $\HDG{\Omega}$ in blue and the interface $\Gamma_I$ in black.
%

The subdomain $\CG{\Omega}$ features a compressible and stiff material with Young's mudulus $E {=} 250$ and Poisson's ratio $\nu {=} 0.35$, whereas the properties of the nearly incompressible region $\HDG{\Omega}$ are $E {=} 80$ and $\nu$ spanning among the values $\left\{ 0.49, 0.4999, 0.499999 \right\}$.
The plate is clamped on the left and is subject to a shear load uniformly distributed along the positive $y$-direction with total load $F {=} 100$ on the right.
The remaining boundaries are free surfaces on which a homogenous Neumann condition is imposed.

The bimaterial Cook's membrane problem is solved using the four strategies discussed in Section~\ref{sc:2DelasticLocking}.
The HDG stabilization parameter $\tau$ is considered constant an equal to $10$, whereas the Nitsche's parameter $\gamma$ is set to $10^4$.
The HDG approximation being locking-free~\cite{soon2009hybridizable} is considered as reference solution for the problem under analysis.

Figure~\ref{fig:cooksMembraneSolution} compares the approximation of the displacement and the stress fields in the nearly incompressible case of  $\nu {=} 0.4999$, for the different discretization strategies described above, using polynomial of degree $k {=} 1$ on the fifth level of mesh refinement.
It is straightforward to observe that CG suffers from locking phenomena providing an unreliable approximation of both displacement and stress fields.
The two coupling strategies provide comparable results for the displacement field, whereas a slightly more accurate approximation of the stress field is achieved by the $\KPdeg{CG}$-$\Kdeg{HDG}$ coupling in the compressible region. This is due to the increased accuracy of the CG approximation in the subdomain $\CG{\Omega}$, already observed for this approach in Section~\ref{sc:2DelasticLocking}.

\begin{figure}
\centering
\subfigure[$\nu {=} 0.49$]{\includegraphics[width=0.32\columnwidth]{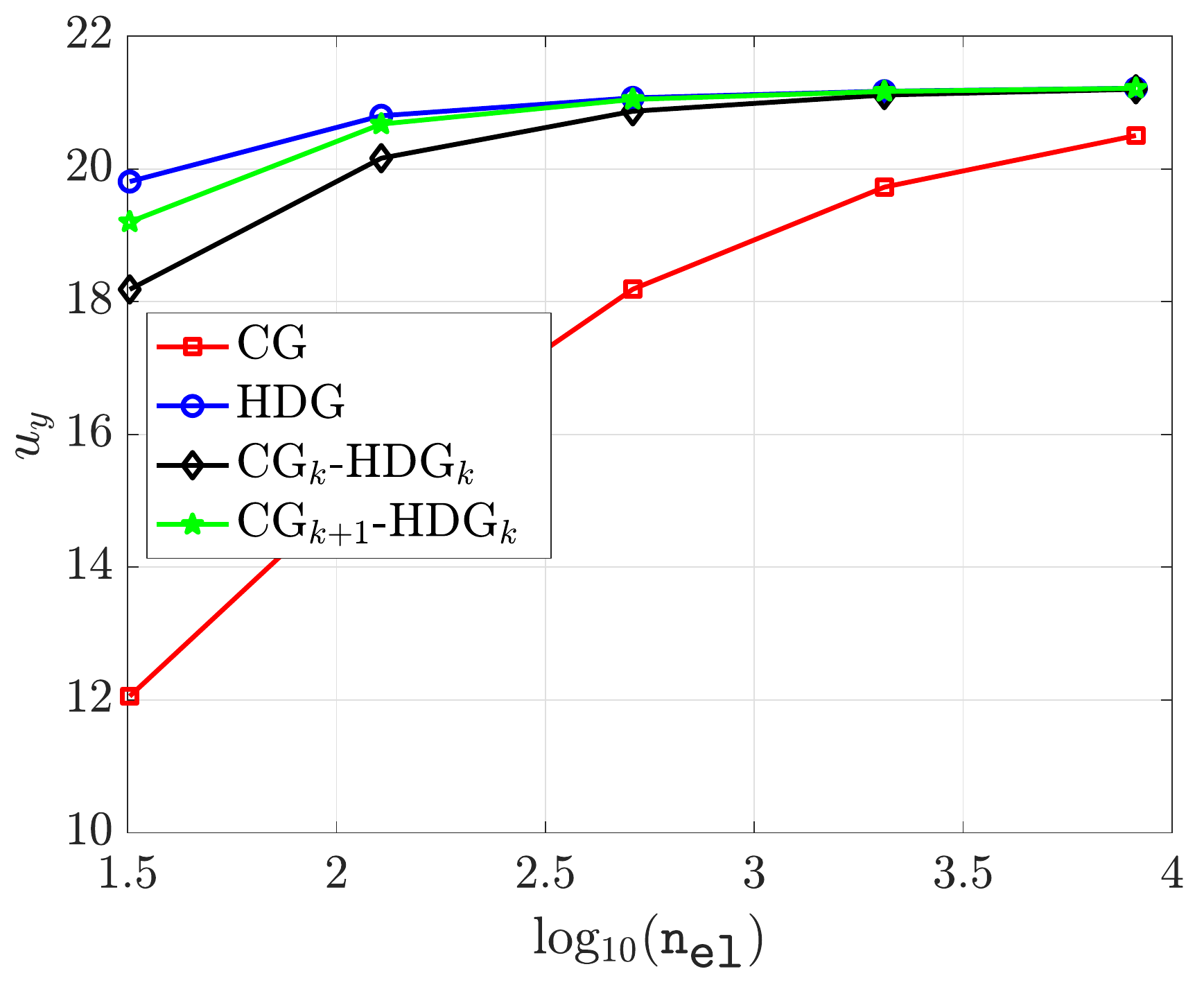}}
\hfill
\subfigure[$\nu {=} 0.4999$]{\includegraphics[width=0.32\columnwidth]{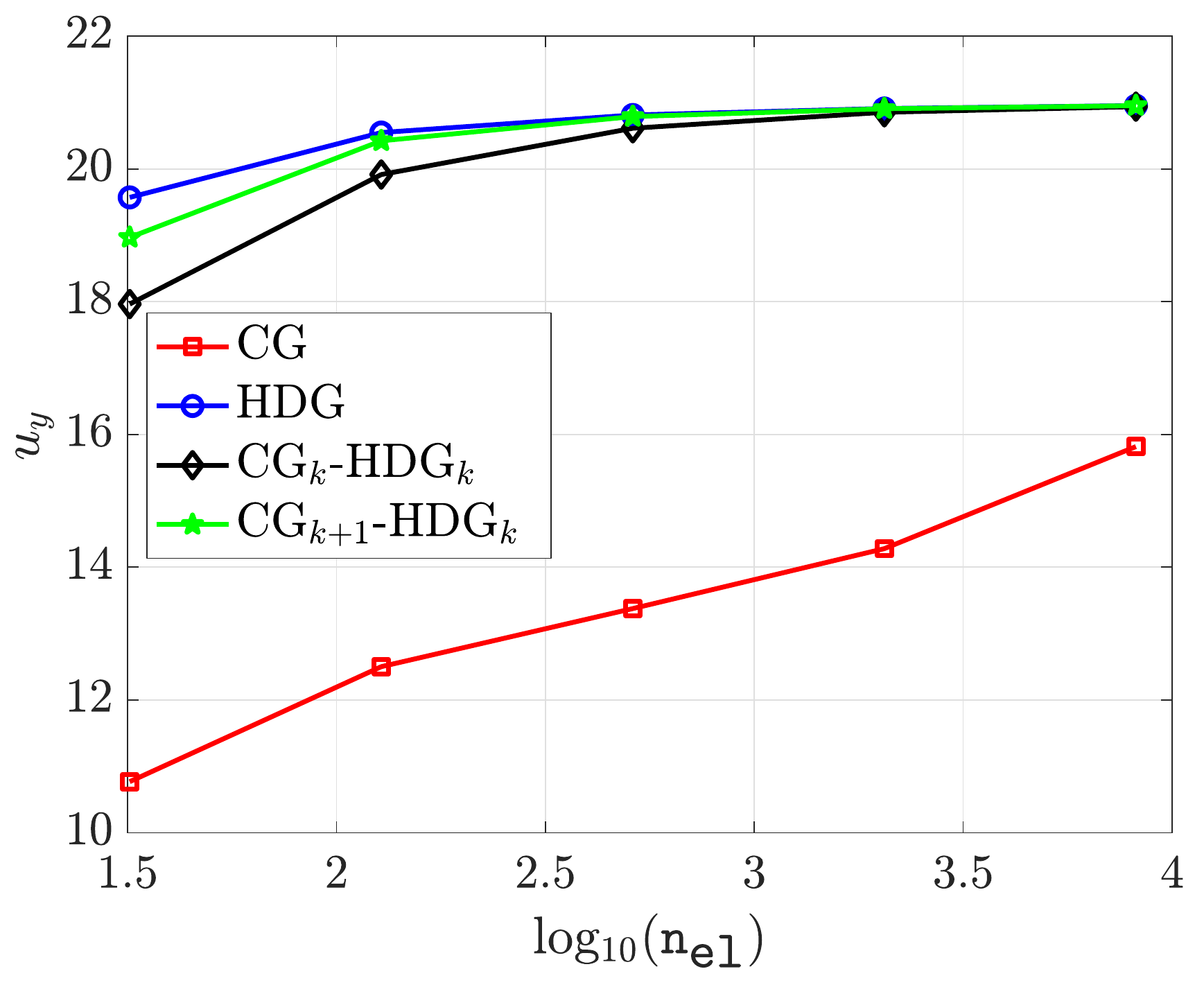}}
\hfill
\subfigure[$\nu {=} 0.499999$]{\includegraphics[width=0.32\columnwidth]{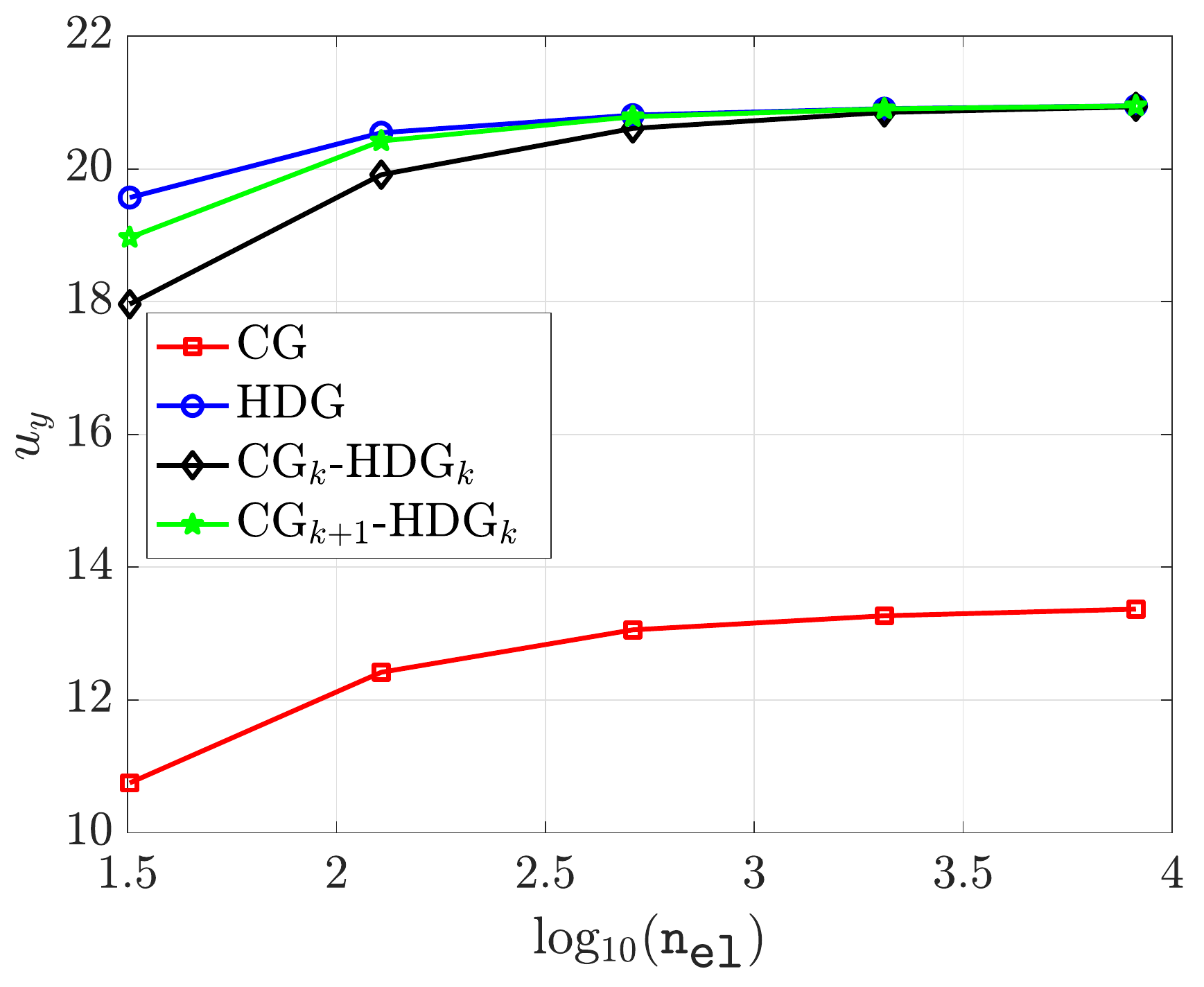}}
\caption{Vertical displacement of the top corner of the right end of the plate as a function of the number of mesh elements in the domain $\Omega$.}
\label{fig:cooksMembraneTip}
\end{figure}

The evolution of the vertical displacement of the top corner of the right end of the plate is displayed in Figure~\ref{fig:cooksMembraneTip} as a function of the number of mesh elements in the domain $\Omega$, for different values of the Poisson's ratio.
The CG approximation with polynomial functions of degree $k {=} 1$ leads to a significant underestimation of the displacement of the tip and this error increases as the Poisson's ratio tends to $0.5$ due to locking effects (Fig.~\ref{fig:cooksMembraneTip}c).
A good approximation of the displacement of the tip is recovered using both coupling strategies. As highlighted above, the $\KPdeg{CG}$-$\Kdeg{HDG}$ approach exploits the additional accuracy of the CG discretization in $\CG{\Omega}$ to construct a more reliable approximation, even when coarser meshes are utilized. Moreover, it is worth noticing that the CG-HDG couplings are locking-free and no loss of accuracy is experienced in the incompressible limit (Fig.~\ref{fig:cooksMembraneTip}c).

\subsection{Three-dimensional laminated composite beam}
\label{sec:3dlaminatedcompositebeam}

The last example considers a three-dimensional laminated composite beam to show the capability of the proposed coupling to treat problems of interest for engineering applications.
The beam $\Omega {=} [-1,1] {\times} [-1,1] {\times} [0,10]$ consists of four layers with alternating mechanical properties, namely a compressible isotropic material with Young's modulus $E=10$ and Poisson's ratio $\nu=0.3$ in $\CG{\Omega} {=} [-1,1] {\times} [-0.5,0] {\times} [0,10] \cup [-1,1] {\times} [0.5,1] {\times} [0,10]$ and a nearly incompressible one with $E=1$ and $\nu=0.49999$ in
$\HDG{\Omega} {=} \Omega \setminus \CG{\Omega}$.
The interfaces between the two materials are identified by the surfaces $\Gamma_I {=} \left\{ (x,y,z) \in \RR^3 : \right.$ $\left. y{=}{-}0.5, y{=}0 \text{ or } y{=}0.5 \right\}$.

\begin{figure}
\centering
\includegraphics[width=0.45\columnwidth]{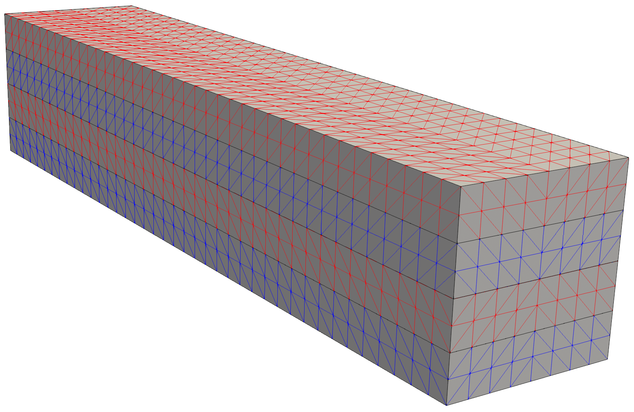}
\hfill
\includegraphics[width=0.45\columnwidth]{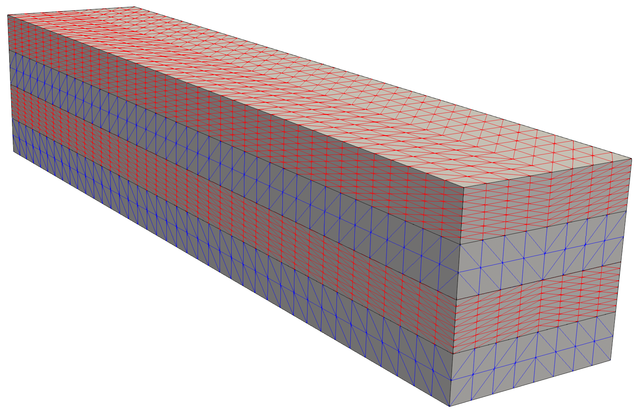}
\caption{Mesh configurations of the laminated composite beam. Coarse (left) and fine (right) configurations. The compressible region $\CG{\Omega}$ is displayed in red, the incompressible one $\HDG{\Omega}$ in blue and the interfaces $\Gamma_I$ in black.}
\label{fig:beamMesh}
\end{figure}
\begin{figure}
	\centering
	\begin{tabular}[c]{@{}c@{}c@{}c@{}}
		\parbox[c]{0.32\textwidth}{\includegraphics[width=0.32\textwidth]{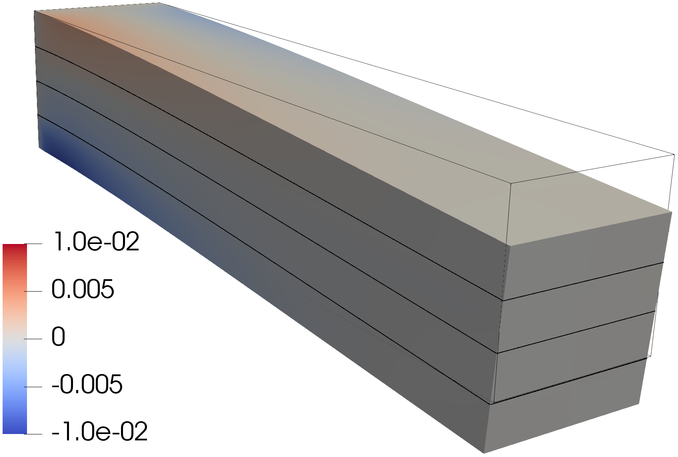}} & 
	 	\parbox[c]{0.32\textwidth}{\includegraphics[width=0.32\textwidth]{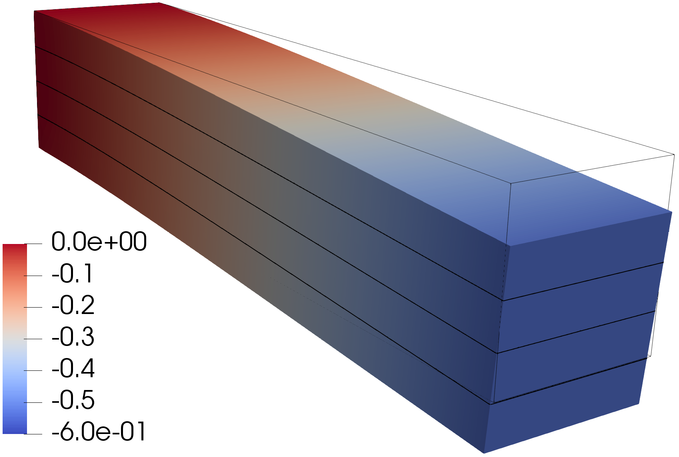}} &
		\parbox[c]{0.32\textwidth}{\includegraphics[width=0.32\textwidth]{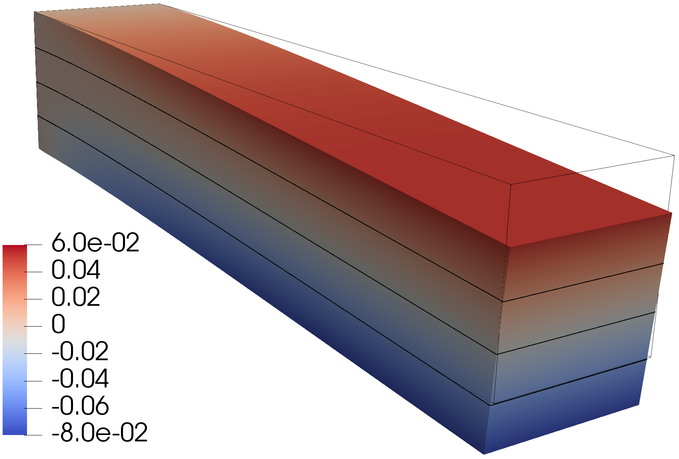}} \\[7em]
		\parbox[c]{0.32\textwidth}{\includegraphics[width=0.32\textwidth]{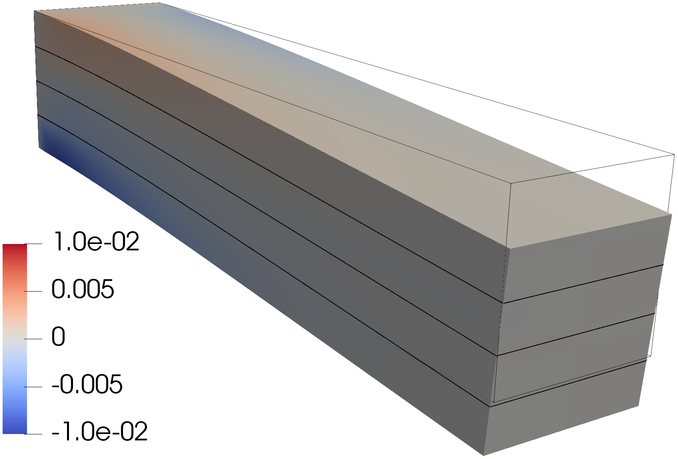}} & 
	 	\parbox[c]{0.32\textwidth}{\includegraphics[width=0.32\textwidth]{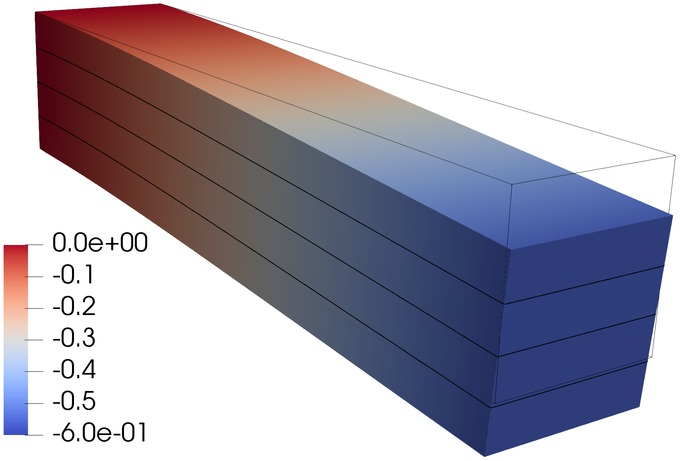}} &
		\parbox[c]{0.32\textwidth}{\includegraphics[width=0.32\textwidth]{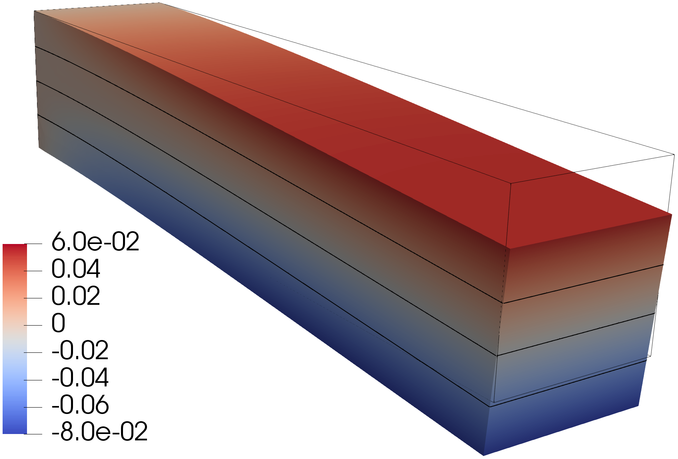}} \\[7em]
		$u_x$ & $u_y$ & $u_z$ \\
	\end{tabular}
\caption{Approximation of the displacement field in the laminated composite beam based on the $\Kdeg{CG}$-$\Kdeg{HDG}$ coupling and polynomial approximation of degree $k {=} 1$, with a coarse (top) and fine (bottom) mesh of the compressible region $\CG{\Omega}$. The interfaces $\Gamma_I$ are displayed in black.}
\label{fig:compositeBeamU}
\end{figure}

The cantilever beam is clamped on the left surface and a load $F {=} 1.3 {\times} 10^{-3}$ is uniformly distributed on the top surface of the beam, along the negative $y$-direction. The remaining surfaces are free boundaries on which a homogenous Neumann condition is imposed.

A reference solution is computed using HDG on a mesh featuring $414,720$ tetrahedral elements with a polynomial approximation of degree $k {=} 1$, for a total of $7,568,640$ global unknowns representing the displacement on the faces.
The HDG stabilization parameter $\tau$ is considered constant an equal to 10, whereas the Nitsche’s parameter $\gamma$ is set to $10^2$.

Contrary to the previous strategy based on nonuniform polynomial approximations, in this example the degree is set to $k {=} 1$ in all $\Omega$ and the mesh in $\CG{\Omega}$ is refined. 
The goal is to show that an improved accuracy of the discrete solution is achieved by increasing the number of degrees of freedom in the subdomain $\CG{\Omega}$, whereas accurate results are obtained by the corresponsing HDG discretization on the coarse grid.
More precisely, two mesh configurations are considered. The coarse mesh counts $7,680$ tetrahedra both in $\CG{\Omega}$ and $\HDG{\Omega}$, whereas the fine discretization respectively features $30,720$ and  $7,680$ elements in the two subdomains, as displayed in Figure~\ref{fig:beamMesh}. 
The subdomain $\CG{\Omega}$ is represented in red, $\HDG{\Omega}$ in blue and the interfaces $\Gamma_I$ in black.
%

\begin{figure}
\centering
\includegraphics[width=0.45\columnwidth]{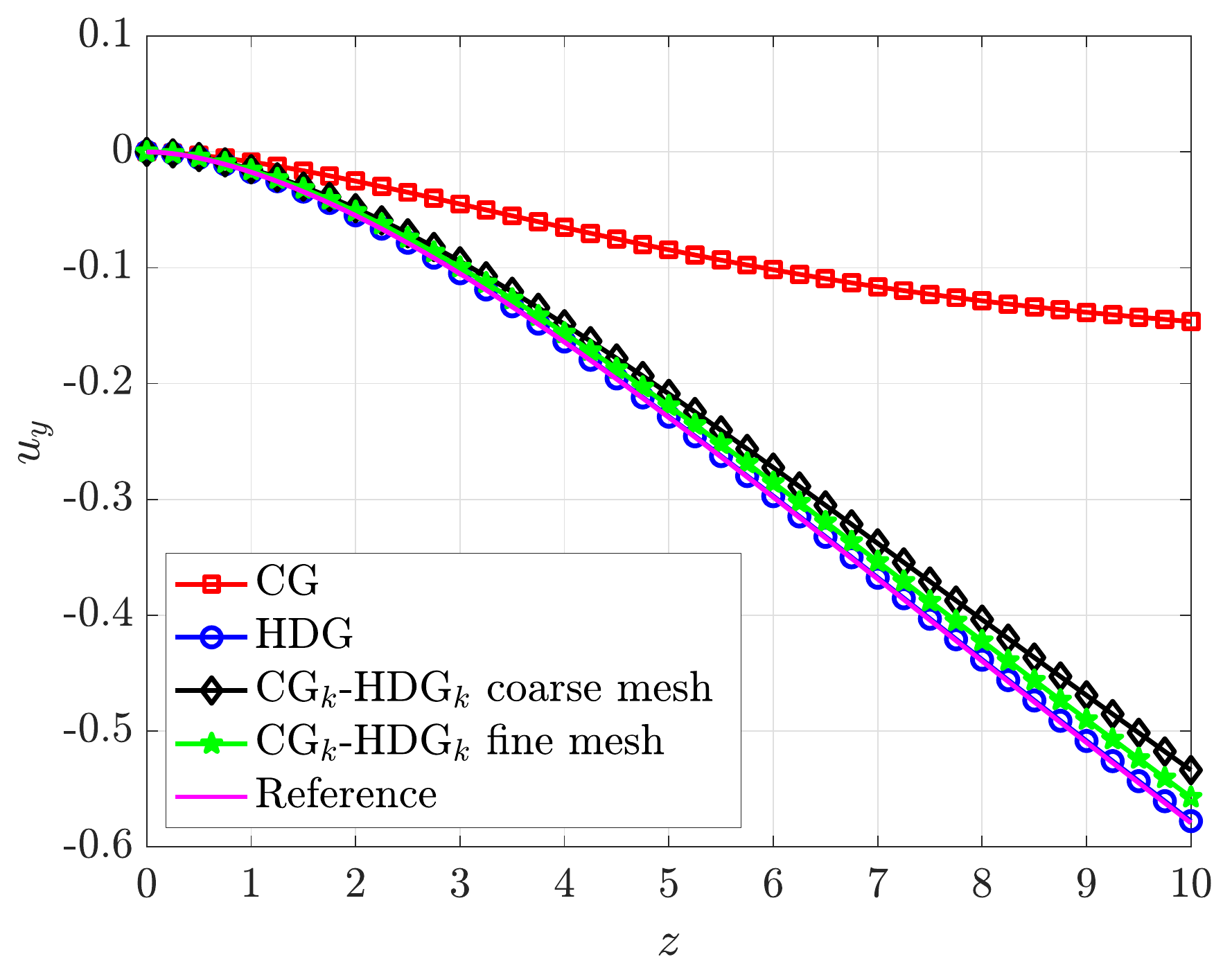}
\caption{Vertical displacement of the axis at $x {=} 0$ and $y {=} 0$ along the beam length.}
\label{fig:compositeBeamTip}
\end{figure}

Figure~\ref{fig:compositeBeamU} displays the comparison of the displacement field computed using the proposed $\Kdeg{CG}$-$\Kdeg{HDG}$ coupling strategy with a polynomial approximation of degree $k {=} 1$ on the meshes presented in Figure~\ref{fig:beamMesh}. 
The results are qualitatively comparable for both approaches. Nonetheless, a small discrepancy is observed in the quantitative evaluation of the vertical displacement of the axis at $x {=} 0$ and $y {=} 0$ along the beam length (Fig.~\ref{fig:compositeBeamTip}).
More precisely, the $\Kdeg{CG}$-$\Kdeg{HDG}$ coupling strategy applied to the configuration with a coarse mesh in the subdomain $\CG{\Omega}$ introduces a relative error for the vertical displacement and the bending moment at midspan of $8\%$, whereas using the fine mesh in $\CG{\Omega}$ they drop to $4\%$ and $2\%$, respectively.
For both quantities of interest, the HDG approximation using the coarse mesh of $\CG{\Omega}$ provides a relative error smaller than $1\%$.
On the contrary, the CG approximation severely underestimates the vertical displacement of the axis.
%

\section{Concluding remarks}
\label{sec:conclusion}

A novel strategy to couple CG and HDG discretizations based on Nitsche's method has been presented in the context of linear thermal and elastic problems.
Special emphasis has been devoted to the case of linear elasticity in which both the computational efficiency of CG in the simulation of compressible materials and the robustness of HDG in the incompressible limit are exploited.

The proposed coupling imposes continuity of the solution and of the trace of the normal numerical flux across the interface between the CG and HDG subdomains. 
The coupled degrees of freedom are only the ones associated with the HDG hybrid variable on the interface. The remaining terms in the formulation, namely the CG discrete matrix and the HDG local and global ones are not modified, leading to a minimally-intrusive implementation in existing CG and HDG libraries.

Numerical examples have been used to verify the optimal convergence of the proposed methodology in the global domain, as well as in the CG and HDG subdomains, separately.
Moreover, a strategy based on nonuniform polynomial degree approximation with a CG discretization of degree $k {+} 1$ and an HDG one of degree $k$ has been proved to lead to a global convergence of the flux/stress of order $k {+} 1$ and a superconvergence of the solution with order $k {+} 2$ by exploiting the inexpensive element-by-element HDG postprocess strategy proposed in~\cite{RS-SGKH:18,MG-GKSH:18}.

\subsection*{Acknowledgements}
This work is partially supported by the European Union's Horizon 2020 research and innovation programme under the Marie Sk\l odowska-Curie actions (Grant No. 675919), the Spanish Ministry of Economy and Competitiveness (Grant No. DPI2017-85139-C2-2-R) and the Generalitat de Catalunya (Grant No. 2017-SGR-1278). Andrea La Spina is supported by the European Education, Audiovisual and Culture Executive Agency (EACEA) under the Erasmus Mundus Joint Doctorate Simulation in Engineering and Entrepreneurship Development (SEED), FPA 2013-0043.

\bibliographystyle{unsrt}
\bibliography{Ref_couplingCG_HDG}

\end{document}